\documentclass[11pt,a4paper,twoside]{scrartcl}

\usepackage{amsfonts, amsmath, amssymb, amsthm}

\usepackage{graphicx}
\graphicspath{{pics/}}

\usepackage{enumerate}
\usepackage{verbatim}

\usepackage{
float,		subcaption, xfrac,		colortbl,
mathtools,  pifont,  	mathrsfs	}

\usepackage[symbol]{footmisc}
\usepackage{enumitem} 
\usepackage{hyperref}

\usepackage{pgfplots}
\pgfplotsset{compat=1.8}
\usepackage{tikz}
\usepackage{tikz-cd}

\DeclareMathOperator*{\diag}{diag}

\newtheorem{theorem}{Theorem}[section]

\newtheorem{remark}[theorem]{Remark}

\newtheorem{example}[theorem]{Example}
\newtheorem{corollary}[theorem]{Corollary}

\newtheorem{alg}[theorem]{Algorithm}

\newenvironment{Theorem}{\goodbreak \begin{theorem}\slshape}{\end{theorem}}
\newenvironment{Remark}{\goodbreak \begin{remark}\slshape}{\end{remark}}
\newenvironment{Example}{\goodbreak \begin{example}\slshape}{\end{example}}
\newenvironment{Corollary}{\goodbreak \begin{corollary}\slshape}{\end{corollary}}

\makeatletter
\def\imod#1{\allowbreak\mkern10mu({\operator@font mod}\,\,#1)}
\makeatother

\makeatletter 
\let\c@algorithm\c@theorem  \makeatother

\newenvironment{algorithm}[1]{\goodbreak~\begin{alg}[#1]~\vspace{-9pt}~\\
	\rule{\linewidth}{0.5pt}~\\}{\vspace{-9pt}~\\
	\rule{\linewidth}{0.5pt}~\end{alg}}

\numberwithin{equation}{section}
\numberwithin{table}{section}
\numberwithin{figure}{section}

\newcommand{\e}{\mathrm e}
\renewcommand{\i}{\mathrm i}
\newcommand{\sinc}{\mathrm{sinc}}
\renewcommand{\b}{\boldsymbol} 

\newcommand{\new}[1]{\textcolor{black}{ #1}}

\newcommand{\R}{\mathbb R}
\newcommand{\C}{\mathbb C}
\newcommand{\Z}{\mathbb Z}
\newcommand{\N}{\mathbb N}
\newcommand{\T}{\mathbb T}

\newcommand{\I}{\mathcal I}

\newcommand{\ex}{\hspace*{0ex} \hfill \hbox{\vrule height
	1.5ex\vbox{\hrule width 1.4ex \vskip 1.4ex\hrule  width 1.4ex}\vrule
	height 1.5ex}}

\long\def\symbolfootnote[#1]#2{\begingroup \def\thefootnote{\fnsymbol{footnote}}\footnote[#1]{#2}\endgroup}

\allowdisplaybreaks

\title{Fast and direct inversion methods for the multivariate nonequispaced fast Fourier transform}

\date{}
\author{Melanie Kircheis\footnotemark[1] \and Daniel Potts\footnotemark[3]}

\hypersetup{pdfauthor={Melanie Kircheis, Daniel Potts},
pdftitle={Fast and direct inversion methods for the multivariate nonequispaced fast Fourier transform},
pdfsubject={},
pdfcreator = {pdflatex},
plainpages=false,
pdfstartview=FitH,       
pdfview=FitH,            
pdfpagemode=UseOutlines, 
bookmarksnumbered=true, 
bookmarksopen=false,     
bookmarksopenlevel=0,   
colorlinks=true,       
linkcolor=black,
citecolor=black,
urlcolor=black}

\begin{document}	

\maketitle

\begin{abstract}
	The well-known discrete Fourier transform (DFT) can easily be generalized to arbitrary points in the spatial domain.
	The fast procedure for this generalization is referred to as nonequispaced fast Fourier transform (NFFT).
	Various applications such as MRI, solution of PDEs, etc., are interested in the inverse problem, i.\,e., computing Fourier coefficients from given nonequispaced data.
	In this paper we survey different kinds of approaches to tackle this problem.
	In contrast to iterative procedures, where multiple iteration steps are needed for computing a solution, we focus especially on so-called direct inversion methods.
	We review density compensation techniques and introduce a new scheme that leads to an exact reconstruction for trigonometric polynomials.
	In addition, we consider a matrix optimization approach using Frobenius norm minimization to obtain an inverse NFFT.

	\smallskip
	
	\noindent AMS \emph{Subject Classifications}: \text{
		65Txx, 
		65T50, 
		65F05. 
	}
\end{abstract}

\footnotetext[1]{Corresponding author: melanie.kircheis@math.tu-chemnitz.de, Chemnitz University of
	Technology, Faculty of Mathematics, D--09107 Chemnitz, Germany}
\footnotetext[3]{potts@math.tu-chemnitz.de, Chemnitz University of
	Technology, Faculty of Mathematics, D--09107 Chemnitz, Germany}


\section{Introduction} 

The NFFT, short hand for nonequispaced fast Fourier transform or nonuniform fast Fourier transform (NUFFT), respectively, is a fast algorithm to evaluate a trigonometric polynomial
\begin{equation}
	\label{eq:poly}
	f(\b x) = \sum_{\b k \in \I_{\b M}} \hat{f}_{\b k}\, \e^{2\pi\i \b k \b x}
\end{equation}
with given Fourier coefficients \mbox{$\hat f_{\b k}\in\C$}, \mbox{$\b k\in\I_{\b{M}}$},
at nonequispaced points \mbox{$\b x_j \in \left[-\frac 12,\frac 12\right)^d$}, \mbox{$j=1,\dots,N$}, \mbox{$N\in\N$}, 
where \mbox{$\I_{\b{M}} \coloneqq \Z^d \cap \left[-\tfrac{M}{2},\tfrac{M}{2}\right)^d$} with \mbox{$|\I_{\b M}| = M^d$}.
In case we are given equispaced points $\b x_j$ and $|\I_{\b{M}}|=N$, this evaluation can be realized by means of the well-known fast Fourier transform (FFT); an algorithm that is invertible.
However, various applications such as magnetic resonance imaging (MRI), cf.~\cite{MRI22, EgKiPo22}, solution of PDEs, cf.~\cite{fa07}, etc., need to perform an inverse nonequispaced fast
Fourier transform (iNFFT), i.\,e., compute the Fourier coefficients~$\hat f_{\b k}$ from given function evaluations~$f(\b x_j)$ of the trigonometric polynomial~\eqref{eq:poly}.
Hence, we are interested in an inversion also for nonequispaced data.

In general, the number~$N$ of points~$\b x_j$ is independent from the number~$|\I_{\b{M}}|$ of Fourier coefficients~$\hat f_{\b k}$ and hence the nonequispaced Fourier matrix
\begin{equation*}
	\b A \coloneqq \left( \e^{2\pi\i \b k \b x_j} \right)_{j=1,\,\b k \in \I_{\b M}}^{N} 
	\ \in \mathbb C ^{N\times |\I_{\b M}|},
\end{equation*}
which we would have to invert, is rectangular in most cases.
Considering the corresponding linear system \mbox{$\b A \b{\hat f} = \b f$} with \mbox{$\b f \coloneqq \left(f(\b x_j)\right)_{j=1}^N$} and \mbox{$\b{\hat f} \coloneqq (\hat f_{\b k})_{\b k\in\I_{\b{M}}}$}, this can either be overdetermined, if \mbox{$|\I_{\b{M}}| \leq N$}, or underdetermined, if \mbox{$|\I_{\b{M}}| > N$}. 
Generally, this forces us to look for a pseudoinverse solution. 
Moreover, we also require that the nonequispaced Fourier matrix $\b A$ has full rank. 
Eigenvalue estimates in \cite{FeGrSt95, BaGr03, BoPo06, KuPo06, KuNa18} indeed confirm that this condition is satisfied for sufficiently nice sampling sets.

In literature a variety of approaches for an inverse NFFT (iNFFT) can be found. 
This is why we give a short overview.

\paragraph{Iterative inversion methods}

We start surveying the iterative methods. 
For the one-dimensional setting \mbox{$d=1$} with \mbox{$|\I_{\b{M}}| = N$} an algorithm was published in \cite{RuTo18}, 
which is specially designed for jittered equispaced points and is based on the conjugate gradient (CG) method in connection with low rank approximation. 
An approach for the overdetermined case \mbox{$|\I_{\b{M}}| \leq N$} can be found in \cite{FeGrSt95}, where the Toeplitz structure of the matrix product \mbox{$\b A^* \b W \b A$} with a diagonal matrix \mbox{$\b W \coloneqq \diag(w_j)_{j=1}^N$} of Voronoi weights is used to compute the solution iteratively by means of the CG algorithm.

For higher dimensional problems with \mbox{$d \geq 1$} there are several approaches that compute a least squares approximation to the linear system \mbox{$\b A \b{\hat f} = \b f$}.
In the overdetermined case \mbox{$|\I_{\b{M}}| \leq N$}, the given data can typically only be approximated up to a residual \mbox{$\b r \coloneqq \b A \b{\hat f} - \b f$}.
Therefore, the weighted least squares problem 
\begin{align*}
	\underset{\b{\hat f} \in \C^{|\I_{\b{M}}|}}{\text{Minimize }} \, 
	\,\sum_{j=1}^N w_j \,\Bigg|\sum_{\b k \in \I_{\b M}} \hat{f}_{\b k}\, \e^{2\pi\i \b k \b x_j} – f(\b x_j)\Bigg|^2
\end{align*}
is considered, which is equivalent to solving the weighted normal equations of first kind
\mbox{$\b A^* \b W \b A \b{\hat f} = \b A^* \b W \b f$}
with the diagonal matrix \mbox{$\b W \coloneqq \diag(w_j)_{j=1}^N$} of weights in time domain. 
In \cite{SuFeNo01, fesu02, KnKuPo} these normal equations are solved iteratively by means of CG using the NFFT to realize fast matrix-vector multiplications involving $\b A$, whereas in \cite{PrWa01} a fast convolution is used.

In the consistent underdetermined case \mbox{$|\I_{\b{M}}| > N$} the data can be interpolated exactly and therefore one can choose a specific solution, e.\,g. the one that solves the constrained minimization problem
\begin{align*}
	\underset{\b{\hat f} \in \C^{|\I_{\b{M}}|}}{\text{Minimize }} \, 
	\,\sum_{\b k\in\I_{\b{M}}} \frac{|\hat f_{\b k}|^2}{\hat w_{\b k}} 
	\quad\text{ subject to }\quad
	\b A \b{\hat f} = \b f .
\end{align*} 
This interpolation problem is equivalent to the weighted normal equations of second kind
$\b A \b{\hat W} \b A^* \b y = \b f$, $\b{\hat f} = \b{\hat W} \b A^* \b y$
with the diagonal matrix \mbox{$\b{\hat W} \coloneqq \diag(\hat w_{\b k})_{\b k\in\I_{\b{M}}}$} of weights in frequency domain. 
In \cite{kupo04} the CG method was used in connection with the NFFT to iteratively compute a solution to this problem, see also \cite[Section~7.6.2]{PlPoStTa18}.

\paragraph{Regularization methods}

Moreover, there also exist several regularization techniques for the multidimensional setting \mbox{$d \geq 1$}.
For example, \cite{Gelb10, Gelb15, Gelb16} all solve the $\ell_1$-regularized problem
\begin{align*}
	\underset{\b{\hat f} \in \C^{|\I_{\b{M}}|}}{\text{Minimize }} \, 
	\tfrac 12 \|\b A \b{\hat f} - \b f\|_2^2 + \lambda \|\mathcal L^m \b{\hat f}\|_1
\end{align*}
with regularization parameter \mbox{$\lambda>0$} and the $m$-th order polynomial annihilation operator \mbox{$\mathcal L^m \in \R^{N \times |\I_{\b{M}}|}$} as sparsifying transform, see~\cite{Gelb05}.
Based on this, weighted $\ell_p$-schemes 
\begin{align*}
	\underset{\b{\hat f} \in \C^{|\I_{\b{M}}|}}{\text{Minimize }} \, 
	\tfrac 12 \|\b A \b{\hat f} - \b f\|_2^2 + \tfrac{1}{p} \|\b W \mathcal L^m \b{\hat f}\|_p^p 
\end{align*}
were introduced in \cite{CaWaBo08, ChYi08, DaDeFoG10, LiMa12}, which are designed to reduce the penalty at locations where \mbox{$\mathcal L^m \b{\hat f} $} is nonzero.
For instance, \cite{Gelb18, Gelb19} each state a two step method, that firstly uses edge detection to create a mask, i.\,e., a weighting matrix which indicates where non-zero entries are expected in the TV domain, and then targets weighted $\ell_2$-norm TV regularization appropriately to smooth regions of the function in a second minimization step. 

\paragraph{Direct inversion methods}

In contrast to these iterated procedures, there are also so-called direct methods that do not require multiple iteration steps.
Already in \cite{duro95} a direct method was explained for the setting \mbox{$d=1$} and \mbox{$|\I_{\b{M}}|=N$} which uses Lagrange interpolation in combination with fast multiple methods. 
Based on this, further methods were deduced for the same setting, which also use Lagrange interpolation, 
but additionally incorporate an imaginary shift in \cite{Selva18}, or utilize fast summation in \cite{KiPo19} for the fast evaluation of occurring sums, see also \cite[Section~7.6.1]{PlPoStTa18}. 

In the overdetermined setting \mbox{$|\I_{\b{M}}| \leq N$} another approach for computing an inverse NFFT can be obtained by using the fact that \mbox{$\b A^* \b A$} is of Toeplitz structure. 
To this end, the Gohberg-Semencul formula, see \cite{HeRo84}, can be used to solve the normal equations \mbox{$\b A^* \b A \b{\hat f} = \b A^* \b f$} exactly by analogy with \cite{AvShSh16}. 
Here the computation of the components of the Gohberg-Semencul formula can be viewed as a precomputational step.
In addition, also a frame-theoretical approach is known from \cite{GeSo14}, which provides a link between the adjoint NFFT and frame approximation, and could therefore be seen as a way to invert the NFFT.
Note that the method in \cite{GeSo14} is based only on optimizing a diagonal matrix (the matrix $\b D$ defined in \eqref{eq:matrix_D}), whereas in \cite{KiPo19} similar ideas were used to modify a sparse matrix (the matrix $\b B$ defined in \eqref{eq:matrix_B}).

For the multidimensional setting \mbox{$d > 1$} several methods have been developed that are tailored to the special structure of the linogram or pseudo-polar grid, respectively, see Figure~\ref{fig:polar_grids_linogram}, such that the inversion involves only one-dimensional FFTs and interpolations.
On the one hand, in \cite{AvCoDoElIs} a least squares solution is computed iteratively by using a fast multiplication technique with the matrix~$\b A$, which can be derived in case of the linogram grid.
On the other hand, \cite{AvCoDoIsSh08} is based on a fast resampling strategy, where a first step resamples the linogram grid onto a Cartesian grid, and the second phase recovers the image from these Cartesian samples.
However, these techniques are exclusively applicable for the special case of the linogram grid, see Figure~\ref{fig:polar_grids_linogram}, or the polar grid by another resampling, cf.~Figure~\ref{fig:polar_grids_mpolar}.
Since we are interested in more generally applicable methods, a brief introduction to direct inversion for general sampling patterns can be found in \cite{KiPo20}.

\paragraph{Current approach}

In this paper we focus on direct inversion methods that are applicable for general sampling patterns and present new methods for the computation of an iNFFT. 
Note that direct method in the context of the linear system \mbox{$\b A \b{\hat f} = \b f$} means, that for a fixed set of points $\b x_j$, \mbox{$j=1,\dots, N$}, the reconstruction of $\b{\hat f}$ from given $\b f$ can be realized with the same number of arithmetic operations as a single application of an adjoint NFFT (see Algorithm~\ref{alg:nfft*}).
To achieve this, a certain precomputational step is compulsory, since the adjoint NFFT does not yield an inversion of the NFFT per se, cf.~\eqref{eq:matrix_product_nonequi}. 
Although this precomputations might be rather costly, they need to be done only once for a given set of points 
$\b x_j$, \mbox{$j=1,\dots,N$}. 
Therefore, direct methods are especially beneficial in case of fixed points for several measurement vectors~$\b f$.

For this reason, the current paper is concerned with two different approaches of this kind.
Firstly, we consider the very well known approach of so-called sampling density compensation,
which can be written as \mbox{$\b{\hat f} \approx \b A^* \b W \b f$}
with a diagonal matrix \mbox{$\b W \coloneqq \diag(w_j)_{j=1}^N$} of weights.
The already mentioned precomputations then consist of computing suitable weights~$w_j$, while the actual reconstruction step includes only one adjoint NFFT applied to the scaled coefficient vector \mbox{$\b W \b f$}.
In this paper we examine several existing approaches for computing the weights~$w_j$ and introduce a new method, such that the iNFFT is exact for all trigonometric polynomials \eqref{eq:poly} of degree $\b M$.

As a second part, we reconsider and enhance our approach introduced in \cite{KiPo20}. 
Here the idea is using the matrix representation \mbox{$\b A \approx \b B \b F \b D$} of the NFFT to modify one of the matrix factors, such that an inversion is given by \mbox{$\b{\hat f} \approx \b D^* \b F^* \b B_{\mathrm{opt}}^* \b f$}.
Then the precomputational step consists of optimizing the matrix $\b B$, while the actual reconstruction step includes only one modified adjoint NFFT applied to the coefficient vector \mbox{$\b f$}.

\paragraph{Outline of this paper}

This paper is organized as follows.
In Section~\ref{sec:nfft} we introduce the already mentioned algorithm, the NFFT, as well as its adjoint version, the adjoint NFFT.
Secondly, in Section~\ref{sec:inv_density_comp} we introduce the inversion problem and deal with direct methods using so-called sampling density compensation.
We start our investigations with trigonometric polynomials in Section~\ref{sec:trig_poly}.
Here the main formula \eqref{eq:claim_quadrature_coefficients_double} yields exact reconstruction for all trigonometric polynomials of degree $\b M$.
We also discuss practical computation schemes for the overdetermined as well as the underdetermined setting.
Subsequently, in Section~\ref{sec:bandlim} we go on to bandlimited functions and show that the same numerical procedures as in Section~\ref{sec:trig_poly} can be used in this setting as well.
Section~\ref{sec:error_bound} then summarizes the previous findings by presenting a general error bound on density compensation factors computed by means of \eqref{eq:claim_quadrature_coefficients_double} in \autoref{Thm:error_est_lgs}.
In addition, this also yields an estimate on the condition number of a specific matrix product, as shown in \autoref{Thm:est_cond}.
In Section~\ref{sec:dcf_literature} we have a look at certain approaches from literature and their connection among each other as well as to the method presented in Section~\ref{sec:trig_poly}.
Afterwards, we examine another direct inversion method in Section~\ref{sec:opt_B}, where we aim to modify the adjoint NFFT by means of matrix optimization such that this yields an iNFFT.
Finally, in Section~\ref{sec:numerics} we show some numerical examples to investigate the accuracy of our approaches.

\section[Nonequispaced fast Fourier transform (NFFT)]{Nonequispaced fast Fourier transform \label{sec:nfft}}

Let
\begin{equation*}
	\label{eq:torus}
	\T^d \coloneqq \R^d\setminus\Z^d
	\cong 
	\left[-\tfrac 12, \tfrac 12\right)^d
	=
	\left\{ \b x\in\mathbb R^d \colon -\tfrac 12 \leq x_t < \tfrac 12,\, t=1,\dots,d \right\}
\end{equation*}
denote the $d$-dimensional torus with $d\in\N$. For $\b M \coloneqq (M,\dots,M)^T$, \mbox{$M\in 2\N$}, we define the multi-index set
\begin{equation*}
	\I_{\b M}
	\coloneqq
	\Z^d \cap \left[-\tfrac{M}{2},\tfrac{M}{2}\right)^d
	=
	\left\{ \b k \in \Z^d \colon -\tfrac{M}{2} \leq k_t < \tfrac{M}{2},\, t=1,\dots,d \right\}
\end{equation*}
with cardinality \mbox{$|\I_{\b M}| = M^d$}.
The inner product of two vectors shall be defined as usual as
\mbox{$\b k \b x \coloneqq k_1 x_1 + \dots + k_d x_d$}.
Additionally, we define the componentwise product as
\mbox{$\b x \odot \b y \coloneqq \left( x_1 y_1, \dots, x_d y_d \right)^T$},
the all ones vector \mbox{$\b 1_d \coloneqq (1,\dots,1)^T$} and 
the reciprocal of a vector $\b x$ with nonzero components shall be given by
\mbox{$\b x^{-1} \coloneqq \left( x_1^{-1},\dots,x_d^{-1} \right)^T$}.

We consider the Hilbert space $L_2(\T^d)$ of all $1$-periodic, complex-valued functions, which possesses the orthonormal basis \mbox{$\{\e^{2\pi\i \b k \b x} \colon \b k\in\Z^d\}$}.
It it known that every function \mbox{$f\in L_2(\T^d)$} is uniquely representable in the form
\begin{align}
	\label{eq:Fourier_series}
	f(\b x) = \sum_{\b k \in \Z^d} c_{\b k}(f) \,\e^{2\pi\i \b k \b x}
\end{align}
with the coefficients
\begin{align}
	\label{eq:Fourier_coeffs}
	c_{\b k}(f) \coloneqq \int_{\T^d} f(\b x) \,\e^{-2\pi\i\b k\b x} \,\mathrm d\b x,
	\quad \b k\in\Z^d,
\end{align}
where the sum in \eqref{eq:Fourier_series} converges to $f$ in the $L_2(\T^d)$-norm, cf. \cite[Thm.~4.5]{PlPoStTa18}.
A series of the form \eqref{eq:Fourier_series} is called \emph{Fourier series} with the \emph{Fourier coefficients}~\eqref{eq:Fourier_coeffs}.
Numerically, the Fourier coefficients \eqref{eq:Fourier_coeffs} are approximated on the uniform grid \mbox{$\{\b M^{-1} \odot \b\ell,\, \b\ell\in\I_{\b M}\}$} by the trapezoidal rule for numerical integration as
\begin{align}
	\label{eq:Fourier_coeffs_equi}
	c_{\b k}(f) 
	\approx 
	\frac{1}{|\I_{\b M}|} \sum_{\b\ell\in\I_{\b M}} f(\b M^{-1} \odot \b\ell) \,\e^{-2\pi\i\b k(\b M^{-1} \odot \b\ell)},
	\quad \b k\in\Z^d,
\end{align}
which is an acceptable approximation for \mbox{$\b k\in\I_{\b M}$}, see e.\,g.~\cite[p.~214]{PlPoStTa18}.
The fast evaluation of \eqref{eq:Fourier_coeffs_equi} can then be realized by means of the \emph{fast Fourier transform (FFT)}.
Moreover, it is know that this transformation is invertible and that the inverse problem of computing
\begin{equation}
	\label{eq:idft}
	f(\b M^{-1} \odot \b\ell) = \sum_{\b k \in \I_{\b M}} \hat{f}_{\b k}\, \e^{2\pi\i \b k (\b M^{-1} \odot \b\ell)},
	\quad \b\ell\in\I_{\b M},
\end{equation}
with \mbox{$\hat{f}_{\b k} \approx c_{\b k}(f)$}, \mbox{$\b k \in \I_{\b M}$}, can be realized by means of an \emph{inverse fast Fourier transform (iFFT)}, which is basically the same algorithm except for some reordering and scaling, cf.~\cite[Lem.~3.17]{PlPoStTa18}.

Now suppose we are given nonequispaced points \mbox{$\b x_j \in \T^d$},\, \mbox{$j=1,\dots,N$}, 
instead.
Then, we consider the computation of the sums 
\begin{equation}
	\label{eq:nfft}
	f_j \coloneqq f(\b x_j) = \sum_{\b k \in \I_{\b M}} \hat{f}_{\b k}\, \e^{2\pi\i \b k \b x_j}, \quad j=1,\dots,N,
\end{equation}
for given \mbox{$\hat{f}_{\b k} \in \C$}, \mbox{$\b k \in \I_{\b M}$}, as well as the adjoint problem of computing
\begin{equation}
	\label{eq:nfft*}
	h_{\b k} = \sum_{j=1}^N f_j\, \e^{-2\pi\i \b k \b x_j}, \quad \b k \in \I_{\b M},
\end{equation}
for given \mbox{$f_j \in \C$}, \mbox{$j=1,\dots,N$}.
By defining the nonequispaced Fourier matrix
\begin{equation}
	\label{eq:matrix_A}
	\b A = \b A_{|\I_{\b{M}}|} \coloneqq \left( \e^{2\pi\i \b k \b x_j} \right)_{j=1,\,\b k \in \I_{\b M}}^{N} 
	\ \in \mathbb C ^{N\times |\I_{\b M}|},
\end{equation}
as well as the vectors
\mbox{$\b f\coloneqq\left(f_j\right)_{j=1}^N$}, 
\mbox{$\b{\hat f}\coloneqq(\hat f_{\b k})_{\b k \in \I_{\b M}}$}
and
\mbox{$\b h \coloneqq (h_{\b k})_{\b k \in \I_{\b M}}$},
the computation of sums of the form \eqref{eq:nfft} and \eqref{eq:nfft*} can be written as
\mbox{$\b f = \b A \b{\hat f}$}
and
\mbox{$\b h = \b A^* \b f$},
where 
\mbox{$\b A^* \coloneqq \overline{\b A}^{T}$} 
denotes the adjoint matrix of~$\b A$.

Since the naive computation of \eqref{eq:nfft} and \eqref{eq:nfft*} is of complexity \mbox{$\mathcal O(N \cdot |\I_{\b M}|)$}, 
a fast approximate algorithm, the so-called \emph{nonequispaced fast Fourier transform (NFFT)}, is briefly explained below. For more information see \cite{duro93, bey95, st97, GrLe04, KeKuPo09} or \cite[pp.~377-381]{PlPoStTa18}.

\subsection{The NFFT\label{subsec:nfft}}

We firstly restrict our attention to problem~\eqref{eq:nfft}, 
which is equivalent to the evaluation of a trigonometric polynomial
\begin{equation}
	\label{eq:trig_poly_2d}
	f(\b x) = \sum_{\b k \in \I_{\b M}} \hat{f}_{\b k}\, \e^{2\pi\i \b k \b x}
\end{equation}
with given \mbox{$\hat{f}_{\b k} \in \C$}, \mbox{$\b k \in \I_{\b M}$},
at given nonequispaced points~\mbox{$\b x_j \in \T^d$}, \mbox{$j=1,\dots,N$}.
Let \mbox{$w\in L_2(\R^d)\cap L_1(\R^d)$} be a so-called \emph{window function}, which is well localized in space and frequency domain.
Now we define the \mbox{1-periodic} function~\mbox{$\tilde w (\b x) \coloneqq \sum_{\b r\in\Z^d} w(\b x+\b r)$} with absolute convergent Fourier series.
As a consequence, the Fourier coefficients of the periodization $\tilde w$ have the form
\begin{equation*}
	\label{eq:window_FT}
	c_{\b k}(\tilde w) 
	= 
	\int_{\T^d} \tilde w(\b x) \, \e^{-2\pi\i \b k \b x} \,\mathrm d \b x 
	= 
	\int_{\R^d} w(\b x) \, \e^{-2\pi\i \b k \b x} \,\mathrm d \b x 
	\eqqcolon 
	\hat w(\b k).
\end{equation*}
For a given oversampling factor~\mbox{$\sigma\geq 1$},
we define 
\mbox{$2\N \ni  M_{\sigma} \coloneqq 2 \lceil \lceil \sigma M \rceil / 2 \rceil$}
as well as
\mbox{$\b M_{\b\sigma} \coloneqq M_{\sigma} \cdot \b 1_{d}$},
and approximate~$f$ by a linear combination of translates of the periodized window function, i.\,e., 
\begin{equation}
	\label{eq:s1}
	f(\b x) \approx s_1(\b x) 
	\coloneqq 
	\sum_{\b\ell \in \I_{\b M_{\b\sigma}}} g_{\b\ell}\, \tilde w \hspace{-2pt} \left(\b x-\b M_{\b\sigma}^{-1} \odot \b\ell\right),
\end{equation}
where \mbox{$g_{\b\ell}\in\C$}, \mbox{$\b\ell\in\I_{\b M_{\b \sigma}}$}, are coefficients to be determined such that \eqref{eq:s1} yields a good approximation.
By means of the convolution theorem (see \cite[Lem.~4.1]{PlPoStTa18}), the approximant \mbox{$s_1\in L_2(\T^d)$} in \eqref{eq:s1} can be represented as
\begin{align}
	\label{eq:s1_aliasing}
	s_1(\b x) 
	&= 
	\sum_{\b k \in \Z^d} c_{\b k}(s_1)\, \e^{2\pi\i \b k \b x}
	\notag \\
	&= 
	\sum_{\b k \in \I_{\b M}} \hat g_{\b k}\;c_{\b k}(\tilde{w})\, \e^{2\pi\i \b k \b x} 
	+ 
	\sum_{\b r \in \Z^d\setminus\{\b 0\}} \sum_{\b k \in \I_{\b M}} \hat g_{\b k}\;c_{\b k+\b M_{\b\sigma}\odot\,\b r}(\tilde{w})\, \e^{2\pi\i (\b k+\b M_{\b\sigma}\odot\,\b r)\b x}, 
\end{align}
where the discrete Fourier transform of the coefficients $g_{\b\ell}$ is defined by
\begin{equation}
	\label{eq:discrete_FT}
	\hat g_{\b k}
	\coloneqq
	\sum_{\b\ell \in \I_{\b M_{\b\sigma}}} g_{\b\ell}\, \e^{-2\pi\i \b k (\b M_{\b\sigma}^{-1} \odot\, \b\ell)}, 
	\quad \b k \in \I_{\b M}.
\end{equation}
Comparing \eqref{eq:nfft} and \eqref{eq:s1_aliasing} then yields
\begin{equation*}
	\hat g_{\b k} =
	\left\{
	\begin{array}{cl}
		\dfrac{\hat f_{\b k}}{\hat w(\b k)} &: \b k \in \I_{\b M},\\
		0 &: \b k \in \I_{\b M_{\b\sigma}}\setminus\I_{\b M}.
	\end{array}
	\right.
\end{equation*}
Consequently, the coefficients $g_{\b\ell}$ in \eqref{eq:s1} can be obtained by inverting \eqref{eq:discrete_FT}, i.\,e., by the application of an iFFT.

Furthermore, we assume that $w$ is well localized such that it is small outside the square
\mbox{$\left[-\sfrac{m}{M_{\sigma}}, \sfrac{m}{M_{\sigma}}\right]^d$} with truncation parameter \mbox{$m \ll M_{\sigma}$}.
In this case, $w$ can be approximated by the compactly supported function
\begin{equation*}
	w_m(\b x)
	\coloneqq
	\left\{
	\begin{array}{cl}
		w(\b x) &: \b x \in \left[-\frac{m}{M_{\sigma}}, \frac{m}{M_{\sigma}}\right]^d,\\
		0 &: \text{otherwise}.
	\end{array}
	\right.
\end{equation*}
Thereby, we approximate $s_1$ by the short sums
\begin{align*}
	\label{eq:truncation}
	f(\b x_j) 
	\approx 
	s_1(\b x_j) 
	\approx 
	s(\b x_j) 
	&\coloneqq 
	\sum_{\b\ell \in \I_{\b M_{\b\sigma}}} g_{\b\ell}\, \tilde w_m \hspace{-2.5pt} \left(\b x_j-\b M_{\b\sigma}^{-1}\odot\b\ell\right)\\
	&= 
	\sum_{\b\ell \in \I_{\b M_{\b\sigma},m}(\b x_j)} g_{\b\ell}\, \tilde w_m \hspace{-2.5pt} \left(\b x_j-\b M_{\b\sigma}^{-1}\odot\b\ell\right),
\end{align*}
where the index set 
\begin{equation}
	\label{eq:indexset_x}
	\I_{\b M_{\b\sigma},m}(\b x_j)
	\coloneqq
	\left\{ \b\ell \in \I_{\b M_{\b\sigma}} \colon \exists \b z \in \Z^d \text{ with } -m \cdot \b 1_d \leq \b M_{\b\sigma} \odot \b x_j - \b\ell + \b z \leq m \cdot \b 1_d \right\}
\end{equation}
contains at most $(2m+1)^d$ entries for each fixed $\b x_j$.
Thus, the obtained algorithm can be summarized as follows.

\begin{algorithm}{NFFT}
	\label{alg:nfft}
	For \mbox{$d,N \in \N$} let \mbox{$\b x_j \in \T^d,\, j=1, \dots, N,$} be given points as well as \mbox{ $\hat f_{\b k} \in \C$}, \mbox{$\b k \in \I_{\b M}$}, given Fourier coefficients.
	Furthermore, we are given the oversampling factor \mbox{$\sigma\geq 1$}, \mbox{$2\N \ni  M_{\sigma} \coloneqq 2 \lceil \lceil \sigma M \rceil / 2 \rceil$}, \mbox{$\b M_{\b\sigma} \coloneqq M_{\sigma} \cdot \b 1_{d}$}, as well as the window function~$w$, the truncated function~$w_m$ with truncation parameter \mbox{$m \ll M_\sigma,$} and their 1-periodic versions~$\tilde w$ and $\tilde w_m$.
	\begin{enumerate}
		\item Set
		\begin{equation*}
			\hat g_{\b k}\coloneqq
			\left\{
			\begin{array}{cl}
				\frac{\hat f_{\b k}}{\hat w(\b k)} &: \b k \in \I_{\b M},\\
				0 &: \b k \in \I_{\b M_{\b\sigma}} \setminus \I_{\b M}.
			\end{array}
			\right.
		\end{equation*}
		\hfill $\mathcal O(|\I_{\b M}|)$
		\item Compute
		\begin{equation*}
			g_{\b\ell}
			\coloneqq
			\frac{1}{|\I_{\b M_{\b\sigma}}|} \sum_{\b k \in \I_{\b M}}
			\hat g_{\b k}\, \e^{2\pi\i \b k(\b M_{\b\sigma}^{-1}\odot\,\b\ell)},
			\quad \b\ell \in \I_{\b M_{\b\sigma}},
		\end{equation*}
		by means of a $d$-variate iFFT.
		\hfill $\mathcal O(|\I_{\b M}|\log(|\I_{\b M}|))$
		\item Compute the short sums
		\begin{equation*}
			\tilde f_j
			\coloneqq
			\sum_{\b\ell \in \I_{\b M_{\b\sigma},m}(\b x_j)} g_{\b\ell}\,\tilde w_m \hspace{-2.5pt}\left(\b x_j-\b M_{\b\sigma}^{-1}\odot\b\ell\right), \quad j = 1, \dots, N.
		\end{equation*}
		\hfill $\mathcal O(N)$
	\end{enumerate}
	\vspace{1ex} 
	\textnormal{\textbf{Output:}} $\tilde f_j \approx f_j$, $j = 1,\dots, N$, cf. \eqref{eq:nfft}. \hfill
	\textnormal{\textbf{Complexity:}} $\mathcal O(|\I_{\b M}|\log(|\I_{\b M}|) + N)$
\end{algorithm}

\begin{Remark} 
	\label{Remark:window_functions}
	Suitable window functions can be found e.\,g. in \cite{duro93, bey95, DuSc, Fou02, GrLe04, KeKuPo09, PoTa21b}.
	\ex
\end{Remark}

Next we give the matrix-vector representation of the NFFT.
To this end, we define
\begin{itemize}
	\item the diagonal matrix
	\begin{equation}
		\label{eq:matrix_D}
		\b D 
		\coloneqq 
		\text{diag} \left( \frac 1{|\I_{\b M_{\b\sigma}}|\cdot\hat{w}(\b k)} \right)_{\b k \in \I_{\b M}} 
		\ \in \mathbb C^{|\I_{\b M}|\times |\I_{\b M}|},
	\end{equation}
	\item the truncated Fourier matrix
	\begin{equation}
		\label{eq:matrix_F}
		\b F 
		\coloneqq 
		\left( \e^{2\pi\i \b k (\b M_{\b\sigma}^{-1}\odot\,\b\ell)} \right)_{\b\ell \in \I_{\b M_{\b\sigma}},\, \b k \in \I_{\b M}} 
		\ \in \mathbb C ^{|\I_{\b M_{\b\sigma}}|\times |\I_{\b M}|},
	\end{equation}
	\item and the sparse matrix
	\begin{equation}
		\label{eq:matrix_B}
		\b B 
		\coloneqq 
		\bigg( \tilde w_m \hspace{-2.5pt} \left(\b x_j - \b M_{\b\sigma}^{-1} \odot\b\ell\right) \bigg)_{j=1,\, \l \in \I_{\b M_{\b\sigma}}}^{N} 
		\ \in \mathbb R^{N\times |\I_{\b M_{\b\sigma}}|},
	\end{equation}
\end{itemize} 
where by definition \eqref{eq:indexset_x} each row of $\b B$ contains at most \mbox{$(2m+1)^d$} nonzeros.
In doing so, the NFFT in Algorithm~\ref{alg:nfft} can be formulated in matrix-vector notation such that we receive the approximation 
\mbox{$\b A \approx \b B \b F \b D$} of \eqref{eq:matrix_A}, cf. \cite[p.~383]{PlPoStTa18}.
In other words, the NFFT uses the approximation 
\begin{align*}
	\e^{2\pi\i \b k \b x_j} 
	\approx 
	\frac 1{|\I_{\b M_{\b\sigma}}|\cdot\hat{w}(\b k)} 
	\sum_{\b\ell \in \I_{\b M_{\b\sigma},m}(\b x_j)} 
	\e^{2\pi\i \b k(\b M_{\b\sigma}^{-1}\odot\,\b\ell)}
	\,\tilde w_m \hspace{-2.5pt}\left(\b x_j-\b M_{\b\sigma}^{-1}\odot\b\ell\right).
\end{align*}

\begin{Remark}
	It has to be pointed out that because of consistency the factor~$|\I_{\b M_{\b\sigma}}|^{-1}$ is here not located in the matrix~$\b F$ as usual but in the matrix~$\b D$.
	\ex
\end{Remark}

\subsection{The adjoint NFFT\label{subsec:nfft*}}

Now we proceed with the adjoint problem~\eqref{eq:nfft*}.
As already seen, this can be written as 
\mbox{$\b h = \b A^* \b f$} with the adjoint matrix \mbox{$\b A^*$} of \eqref{eq:matrix_A}.
Thus, using the matrices \eqref{eq:matrix_D}, \eqref{eq:matrix_F} and \eqref{eq:matrix_B}
we receive the approximation
\mbox{$\b A^* \approx \b D^* \b F^* \b B^*$},
such that a fast algorithm for the adjoint problem can be denoted as follows.

\begin{algorithm}{adjoint NFFT}
	\label{alg:nfft*}
	For \mbox{$d,N \in \mathbb N$} let \mbox{$\b x_j \in \T^d$}, \mbox{$j=1,\dots,N$}, be given points as well as \mbox{$f_j \in \mathbb C$} given coefficients.
	Further, we are given the oversampling factor~\mbox{$\sigma\geq 1$}, \mbox{$2\N \ni  M_{\sigma} \coloneqq 2 \lceil \lceil \sigma M \rceil / 2 \rceil$}, \mbox{$\b M_{\b\sigma} \coloneqq M_{\sigma} \cdot \b 1_{d}$}, as well as the window function~$w$, the truncated function~$w_m$ with truncation parameter \mbox{$m \ll M_\sigma,$} and their 1-periodic versions~$\tilde w$ and $\tilde w_m$.
	\begin{enumerate}
		\item Compute the sparse sums
		\begin{equation*}
			g_{\b\ell}
			\coloneqq
			\sum_{j=1}^{N} f_j\, \tilde w_m \hspace{-2.5pt} \left(\b x_j-\b M_{\b\sigma}^{-1}\odot\b\ell\right),
			\quad \b\ell \in \I_{\b M_{\b\sigma}}.
		\end{equation*}
		\hfill $\mathcal O(N)$
		\item Compute
		\begin{equation*}
			\hat g_{\b k}
			\coloneqq
			\frac{1}{|\I_{\b M_{\b \sigma}}|} 
			\sum_{\b\ell \in \I_{\b M_{\b\sigma}}} g_{\b\ell}\, \e^{-2\pi\i \b k(\b M_{\b \sigma}^{-1}\odot\,\b\ell)},
			\quad \b k \in \I_{\b M},
		\end{equation*}
		by means of a $d$-variate FFT.
		\hfill $\mathcal O(|\I_{\b M}|\log(|\I_{\b M}|))$
		\item Set
		\begin{equation*}
			\tilde h_{\b k}
			\coloneqq
			\frac{\hat g_{\b k}}{\hat w(\b k)} , \quad \b k \in \I_{\b M}.
		\end{equation*}
		\hfill $\mathcal O(|\I_{\b M}|)$	
	\end{enumerate}
	\vspace{1ex} 
	\textnormal{\textbf{Output:}} $\tilde h_{\b k} \approx h_{\b k}$, $\b k \in \I_{\b M}$, cf. \eqref{eq:nfft*}. \hfill
	\textnormal{\textbf{Complexity:}} $\mathcal O(|\I_{\b M}|\log(|\I_{\b M}|) + N)$
\end{algorithm}

The algorithms presented in this section (Algorithms \ref{alg:nfft} and \ref{alg:nfft*}) are part of the NFFT software~\cite{nfft3}. 
For algorithmic details we refer to~\cite{KeKuPo09}.

\section{Direct inversion using density compensation \label{sec:inv_density_comp}}

Having introduced the fast methods for nonequispaced data, we remark that various applications such as MRI, solution of PDEs, etc. are interested in the inverse problem, i.\,e., instead of the evaluation of \eqref{eq:nfft} the aim is computing the Fourier coefficients \mbox{$\hat{f}_{\b k}$}, \mbox{$\b k \in \I_{\b M}$}, from given nonequispaced data \mbox{$f(\b x_j)$}, \mbox{$j=1,\dots,N$}.
Therefore, this section shall be dedicated to this task.

To clarify the major dissimilarity between equispaced and nonequispaced data, we start considering the equispaced case.
When evaluating at the points \mbox{$x_{\b j} = \tfrac 1n \b j \in \T^d$}, \mbox{$\b j \in \I_{\b n}$},
with \mbox{$\b n \coloneqq n \cdot \b 1 _d$} and \mbox{$|\I_{\b n}|=N$},
the nonequispaced Fourier matrix \mbox{$\b A\in\C^{N \times |\I_{\b{M}}|}$} in \eqref{eq:matrix_A} turns into the equispaced Fourier matrix \mbox{$\b F \in \C^{|\I_{\b M_{\b\sigma}}| \times |\I_{\b{M}}|}$} from \eqref{eq:matrix_F}
with \mbox{$|\I_{\b M_{\b\sigma}}|=N$}.
Thereby, it results from the geometric sum formula that 
\begin{equation}
	\label{eq:matrix_product_F*F}
	\b F^* \b F 
	= 
	\Bigg( \sum_{\b j \in \I_{\b n}} \e^{2\pi\i (\b k-\b\ell) \b j/{n}} \Bigg)_{\b k,\b\ell\in\I_{\b{M}}}
	= 
	N \b I_{|\I_{\b M}|},
	\quad\text{ if }\,
	|\I_{\b M}| \leq N,
\end{equation}
as well as 
\begin{equation}
	\label{eq:matrix_product_FF*}
	\b F \b F^*
	=
	\Bigg( \sum_{\b k \in \I_{\b M}} \e^{2\pi\i \b k (\b j-\b h)/{n}} \Bigg)_{\b j,\b h\in\I_{\b{n}}}
	=
	|\I_{\b M}| \cdot \b I_{N},
	\quad\text{ if }\,
	|\I_{\b M}| \geq N 
\end{equation}
and \mbox{$|\I_{\b M}|$} is divisible by $N$.
Thus, in the equispaced setting a one-sided inverse is given by the (scaled) adjoint matrix.
However, when considering arbitrary points \mbox{$\b x_j \in \T^d$}, \mbox{$j=1,\dots,N$}, this property is lost, i.\,e., for the nonequispaced Fourier matrix we have
\begin{align}
	\label{eq:matrix_product_nonequi}
	\b A^* \b A \neq N \b I_{|\I_{\b M}|}
	\quad\text{ and }\quad
	\b A \b A^* \neq |\I_{\b M}| \cdot \b I_{N} .
\end{align}
Because of this, more effort is needed in the nonequispaced setting.

In general, we face the following two problems.
%
\begin{enumerate}
	\item[(1)] 
	Solve the linear system
	\begin{equation}
		\label{eq:problem_infft}
		\begin{split}
			\b A \b{\hat f} &= \b f, 
		\end{split}
	\end{equation}
	i.\,e., reconstruct the Fourier coefficients \mbox{$\b{\hat f} = (\hat f_{\b k})_{\b k \in \I_{\b M}}$} from given function values \mbox{$\b f = (f(\b x_j))_{j=1}^{N}$}.
	This problem is referred to as \textit{inverse NDFT (iNDFT)} and an efficient solver shall be called \textit{inverse NFFT (iNFFT)}.
	\item[(2)]
	Solve the linear system
	\begin{equation}
		\label{eq:problem_infft*}
		\begin{split}
			\b A^* \b f &= \b h, 
		\end{split}
	\end{equation}
	i.\,e., reconstruct the coefficients \mbox{$\b f = (f_j)_{j=1}^{N}$} from given data \mbox{$\b h = (h_{\b k})_{\b k \in \I_{\b M}}$}.
	This problem is referred to as \textit{inverse adjoint NDFT (iNDFT*)} and an efficient solver shall be called \textit{inverse adjoint NFFT (iNFFT*)}.
\end{enumerate}
Note that in both problems the numbers $|\I_{\b{M}}|$ and $N$ are independent, such that the nonequispaced Fourier matrix \mbox{$\b A \in \C^{N \times |\I_{\b{M}}|}$} in \eqref{eq:matrix_A} is generally rectangular.

At first, we restrict our attention to problem \eqref{eq:problem_infft}.
When considering iterative inversion procedures as those mentioned in the introduction, these methods require multiple iteration steps by definition.
Therefore, multiple matrix vector multiplications with the system matrix $\b A$, or rather multiple applications of the NFFT (see Algorithm~\ref{alg:nfft}), are needed to compute a solution.
To reduce the computational effort, we now proceed, in contrast to this iterated procedures, with so-called direct methods.
In the setting of problem \eqref{eq:problem_infft} we hereby mean methods, where for a fixed set of points 
$\b x_j$, \mbox{$j=1,\dots, N$}, the reconstruction of $\b{\hat f}$ from given $\b f$ can be realized with the same number of arithmetic operations as a single application of an adjoint NFFT (see Algorithm~\ref{alg:nfft*}).
To achieve this, a certain precomputational step is compulsory, since the adjoint NFFT does not yield an inversion of the NFFT per se, see~\eqref{eq:matrix_product_nonequi}. 
Although this precomputations might be rather costly, they need to be done only once for a given set of points 
$\b x_j$, \mbox{$j=1,\dots,N$}.
In fact, the actual reconstruction step is very efficient.
Therefore, direct methods are especially beneficial in case we are given fixed points for several measurement vectors~$\b f$.

In this section we focus on a direct inversion method for solving problem \eqref{eq:problem_infft} that utilizes so-called \emph{sampling density compensation}.
To this end, we consider the integral~\eqref{eq:Fourier_coeffs} and introduce a corresponding quadrature formula.
In contrast to the already known equispaced approximation \eqref{eq:Fourier_coeffs_equi} we now assume given arbitrary, nonequispaced points 
\mbox{$\b x_j\in\T^d$}, \mbox{$j=1,\dots,N$}.
Thereby, the Fourier coefficients \eqref{eq:Fourier_coeffs} are approximated by a general quadrature rule using 
quadrature weights \mbox{$w_j\in\C$}, \mbox{$j=1,\dots,N$}, 
which are needed for sampling density compensation due to the nonequispaced sampling.
Thus, for a trigonometric polynomial \eqref{eq:trig_poly_2d} we have 
\begin{equation}
	\label{eq:Fourier_coeffs_nonequi}
	\hat f_{\b k}
	=
	c_{\b k}(f)
	\approx
	h_{\b k}^{\mathrm{w}}
	\coloneqq 
	\sum_{j=1}^{N} 
	w_j\, f(\b x_j)\,\e^{-2\pi\i \b k \b x_j},
	\quad \b k\in\I_{\b M}.
\end{equation}
Using the nonequispaced Fourier matrix $\b A\in\C^{N\times |\I_{\b{M}}|}$ in \eqref{eq:matrix_A}, the diagonal matrix of weights
\mbox{$\b W \coloneqq \mathrm{diag} (w_j)_{j=1}^N\in \mathbb C ^{N\times N}$} as well as the vector \mbox{$\b{h}^{\mathrm{w}} \coloneqq (h_{\b k}^{\mathrm{w}})_{\b k\in\I_{\b M}}$},
the nonequispaced quadrature rule \eqref{eq:Fourier_coeffs_nonequi} 
can be written as \mbox{$\b{\hat f} \approx \b{h}^{\mathrm{w}} \coloneqq \b A^* \b W \b f$}.
For achieving a fast computation method we make use of the approximation of the adjoint NFFT, cf.~Section~\ref{subsec:nfft*}, i.\,e., the final approximation is given by
\begin{equation}
	\label{eq:reconstr_density_nfft}
	\b{\hat f} 
	\approx
	\b{\tilde h}^{\mathrm{w}}
	\coloneqq 
	\b D^* \b F^* \b B^* \b W \b f,
\end{equation}
with the matrices $\b D\in \mathbb C ^{|\I_{\b M}|\times |\I_{\b M}|}$, $\b F\in \mathbb C ^{|\I_{\b M_{\b\sigma}}|\times |\I_{\b M}|}$ and $\b B\in \mathbb R ^{N\times |\I_{\b M_{\b\sigma}}|}$ defined in \eqref{eq:matrix_D}, \eqref{eq:matrix_F} and \eqref{eq:matrix_B}.
In other words, for density compensation methods the already mentioned precomputations consist of computing the quadrature weights \mbox{$w_j\in\C$}, \mbox{$j=1,\dots,N$},
while the actual reconstruction step includes only one adjoint NFFT (see Algorithm~\ref{alg:nfft*}) applied to the scaled measurement vector \mbox{$\b W \b f$}.

The aim of all density compensation techniques is then to choose appropriate weights \mbox{$w_j\in\C$}, \mbox{$j=1,\dots,N$}, 
such that the underlying quadrature \eqref{eq:Fourier_coeffs_nonequi} is preferably exact.
In the following we have a look at the specific choice of the so-called \emph{density compensation factors} $w_j$.

An intuitive approach for density compensation is based on geometry, where each sample is considered as representative of a certain surrounding area, as in numerical integration. 
The weights for each sample can be obtained for instance by constructing a Voronoi diagram and calculating the area of each cell, see e.\,g. \cite{RaPrSiBoEg99}.	
This approach of \emph{Voronoi weights} is well-known and widely used in practice.
However, it does not necessarily yield a good approximation \eqref{eq:reconstr_density_nfft}, which is why we examine some more sophisticated approaches in the remainder of this section.

To this end, this section is organized as follows. 
Firstly, in Section~\ref{sec:trig_poly} we introduce density compensation factors $w_j$, \mbox{$j=1,\dots,N,$} that lead to an exact reconstruction formula \eqref{eq:Fourier_coeffs_nonequi} for all trigonometric polynomials \eqref{eq:trig_poly_2d} of degree~$\b M$.
In addition to the theoretical results, we also discuss methods for the numerical computation.
Secondly, in Section~\ref{sec:bandlim} we show that it is reasonable to consider the inversion problem \eqref{eq:problem_infft} and density compensation via \eqref{eq:reconstr_density_nfft} for bandlimited functions \mbox{$f\in L_1(\R^d)\cap C_0(\R^d)$} as well.
Subsequently, we summarize our previous findings by presenting a general error bound on density compensation factors in Section~\ref{sec:error_bound}.
Finally, in Section~\ref{sec:dcf_literature} we reconsider certain approaches from literature and illustrate their connection among each other as well as to the method introduced in Section~\ref{sec:trig_poly}.

\begin{Remark}
	\label{Remark:opt_density}
	Recapitulating, we have a closer look at some possible interpretation perspectives on the reconstruction~\eqref{eq:reconstr_density_nfft}.
	\begin{enumerate}
		\item[(i)]
		If we define \mbox{$\b g \coloneqq \b W \b f$}, i.\,e., 
		each entry of $\b f$ is scaled with respect to the points~$\b x_j$, \mbox{$j=1,\dots,N$},
		the approximation \eqref{eq:reconstr_density_nfft} can be written as
		\mbox{$\b{\hat f} \approx \b D^* \b F^* \b B^* \b g$}.
		As mentioned before, this coincides with an ordinary adjoint NFFT applied to a modified coefficient vector~$\b g$.
		\item[(ii)] 
		By defining the matrix \mbox{$\b{\tilde B} \coloneqq \b W^* \b B$}, i.\,e.,
		scaling the rows of $\b B$ with respect to the points~$\b x_j$, \mbox{$j=1,\dots,N$},
		the approximation \eqref{eq:reconstr_density_nfft} can be written as
		\mbox{$\b{\hat f} \approx \b D^* \b F^* \b{\tilde B}^* \b f$}.		
		In this sense, density compensation can also be seen as a modification of the adjoint NFFT and its application to the original coefficient vector.
	\end{enumerate}
	Note that (i) is the common viewpoint.
	However, we keep (ii) in mind, since this allows treating density compensation methods as an optimization of the sparse matrix~\mbox{$\b B\in\R^{N \times |\I_{\b M_{\b\sigma}}|}$} in \eqref{eq:matrix_B}, as it shall be done in Section~\ref{sec:opt_B}.
	We remark that density compensation methods allow only $N$ degrees of freedom.
	\ex
\end{Remark}

\subsection{Exact quadrature weights for trigonometric polynomials \label{sec:trig_poly}}

Similar to \cite{GrKuPo09}, we aim to introduce density compensation factors $w_j$, \mbox{$j=1,\dots,N,$} that lead to an exact reconstruction formula \eqref{eq:Fourier_coeffs_nonequi} for all trigonometric polynomials \eqref{eq:trig_poly_2d} of degree~$\b M$.
To this end, we firstly examine certain properties that arise from \eqref{eq:Fourier_coeffs_nonequi} being exact.

\begin{Theorem}
	\label{Thm:exact_trig_poly}
	Let a polynomial degree \mbox{$\b M \in (2\N)^d$}, nonequispaced points \mbox{$\b x_j\in\T^d$}, \mbox{$j=1,\dots, N,$} and quadrature weights \mbox{$w_j\in\C$} be given. 
	Then an exact reconstruction formula \eqref{eq:Fourier_coeffs_nonequi} for trigonometric polynomials \eqref{eq:trig_poly_2d} with maximum degree $\b M$ satisfying 
	\begin{align}
		\label{eq:exact_reconstr}
		\hat f_{\b k} = c_{\b k}(f) = h_{\b k}^{\mathrm{w}}, \quad \b k\in\I_{\b M},
	\end{align}
	implies the following equivalent statements.
	\begin{enumerate}[label=(\roman*), resume]
		\item The quadrature rule
		\label{claim_quadrature}
		\begin{equation}
			\label{eq:claim_quadrature}
			\int_{\T^d} f(\b x) \,\mathrm d \b x = \sum_{j=1}^N w_j f(\b x_j)
		\end{equation}
		is exact for all trigonometric polynomials \eqref{eq:trig_poly_2d} with maximum degree $\b M$. 
		\item The linear system of equations
		\label{claim_quadrature_coefficients}
		\begin{align}
			\label{eq:claim_quadrature_coefficients}
			\left[ \b A^T \b w \right]_{\b k}
			=
			\sum_{j=1}^N w_j \,\e^{2\pi\i\b k\b x_j}
			=
			\delta_{\b 0,\b k} \,
			\new{=
				\left\{
				\begin{array}{ll}
					1 &\colon \b k = \b 0 \\
					0 &\colon \text{otherwise}
				\end{array}
				\right\}},
			\quad\b k\in{\I_{\b M}},
		\end{align}
		is fulfilled with the matrix \mbox{$\b A\in\C^{N\times |\I_{\b{M}}|}$} in \eqref{eq:matrix_A} and \mbox{$\b w \coloneqq \left( w_j \right)_{j=1}^N$}.
	\end{enumerate}
\end{Theorem}

\begin{proof} \phantom{text}
	
	\eqref{eq:exact_reconstr} $\Rightarrow$ \ref{claim_quadrature}: 
	By inserting the definition \eqref{eq:trig_poly_2d} of a trigonometric polynomial of degree~$\b M$ into the integral considered in \eqref{eq:claim_quadrature} we have
	\begin{align}
		\label{eq:approx_quadrature}
		\int_{\T^d} f(\b x) \,\mathrm d \b x
		&=
		\sum_{\b k\in{\I_{\b M}}} \hat{f}_{\b k} \cdot \int_{\T^d} \e^{2\pi\i\b k\b x} \,\mathrm d \b x 
		=
		\sum_{\b k\in{\I_{\b M}}} \hat{f}_{\b k} \cdot \delta_{\b 0,\b k}
		=
		\hat{f}_{\b 0} ,
	\end{align}
	\new{with the Kronecker delta \mbox{$\delta_{\b 0,\b k}$}.}
	Now using the property \eqref{eq:exact_reconstr} as well as definition~\eqref{eq:Fourier_coeffs_nonequi} of $h_{\b k}^{\mathrm{w}}$ we proceed with
	\begin{align}
		\label{eq:claim_quadrature_leftside}
		\hat{f}_{\b 0}
		=
		h_{\b 0}^{\mathrm{w}}
		&=
		\sum_{j=1}^N w_j f(\b x_j)
		\sum_{\b k\in{\I_{\b M}}} 
		\e^{0}
		=
		\sum_{j=1}^N w_j f(\b x_j) ,
	\end{align}
	such that \eqref{eq:approx_quadrature} combined with \eqref{eq:claim_quadrature_leftside} yields the assertion \eqref{eq:claim_quadrature}.

	\ref{claim_quadrature} $\Rightarrow$ \ref{claim_quadrature_coefficients}: Inserting the definition \eqref{eq:trig_poly_2d} of a trigonometric polynomial of degree~$\b M$ into the right-hand side of \eqref{eq:claim_quadrature} implies
	\begin{align}
		\label{eq:claim_quadrature_rightside}
		\sum_{j=1}^N w_j f(\b x_j)
		&=
		\sum_{j=1}^N w_j \sum_{\b k\in\I_{\b M}} \hat{f}_{\b k} \,\e^{2\pi\i\b k\b x_j}
		=
		\sum_{\b k\in\I_{\b M}} \hat{f}_{\b k} \,\sum_{j=1}^N w_j \,\e^{2\pi\i\b k\b x_j} .
	\end{align}
	This together with the property \ref{claim_quadrature} and \eqref{eq:approx_quadrature} leads to
	\begin{align*}
		\hat{f}_{\b 0}
		=
		\sum_{\b k\in\I_{\b M}} \hat{f}_{\b k} \,\sum_{j=1}^N w_j \,\e^{2\pi\i\b k\b x_j} 
	\end{align*}
	and thus to assertion \eqref{eq:claim_quadrature_coefficients}.

	\ref{claim_quadrature_coefficients} $\Rightarrow$ \ref{claim_quadrature}:
	Combining \eqref{eq:approx_quadrature}, \eqref{eq:claim_quadrature_coefficients} and \eqref{eq:claim_quadrature_rightside} yields the assertion via
	\begin{align*}
		\int_{\T^d} f(\b x) \,\mathrm d \b x
		&=
		\sum_{\b k\in{\I_{\b M}}} \hat{f}_{\b k} \cdot \delta_{\b 0,\b k}
		=
		\sum_{\b k\in\I_{\b M}} \hat{f}_{\b k} \,\sum_{j=1}^N w_j \,\e^{2\pi\i\b k\b x_j}
		=
		\sum_{j=1}^N w_j f(\b x_j) .
	\end{align*}
\end{proof}

\begin{Remark}
	Comparable results can also be found in literature.
	A fundamental theorem in numerical integration, see \cite{Tch57}, states that for any integral
	\mbox{$\int_{\T^d} f(\b x) \,\mathrm d\b x$}
	there exists an exact quadrature rule \eqref{eq:claim_quadrature}, i.\,e., optimal points $\b x_j\in\T^d$ and weights $w_j\in\C$, \mbox{$j=1,\dots,N$}, such that \eqref{eq:claim_quadrature} is fulfilled.
	In \cite[Lemma~2.6]{Groechenig20} it was shown that for given points \mbox{$\b x_j\in\T^d$}, \mbox{$j=1,\dots, N,$} certain quadrature weights $w_j$ can be stated by means of frame theoretical considerations which lead to an exact quadrature rule \eqref{eq:claim_quadrature} by definition.
	Moreover, it was shown (cf.~\cite[Lemma~3.6]{Groechenig20}) that these weights 
	are the ones with minimal (weighted) $\ell_2$-norm, which are already known under the name ``least squares quadrature'', see~\cite{Hu09}.
	According to \cite[Sec.~2.1]{Hu09} these quadrature weights $w_j$, \mbox{$j=1,\dots, N,$} can be found by solving a linear system of equations \mbox{$\b\Phi \b w = \b v$}, where \mbox{$\Phi_{\b k,j} = \phi_{\b k}(\b x_j)$} and \mbox{$v_{\b k} = \int_{\T^d} \phi_{\b k}(\b x) \,\mathrm d\b x$} for a given set of basis functions \mbox{$\{\phi_{\b k}\}_{\b k \in\I_{\b M}}$}.
	In our setting we have \mbox{$\phi_{\b k}(\b x)=\e^{2\pi\i\b k \b x}$} and therefore 
	\begin{align*}
		\b\Phi = \left(\e^{2\pi\i\b k \b x_j}\right)_{\b k,j} = \b A^T
		\quad\text{ and }\quad
		v_{\b k} = \int_{\T^d} 1 \cdot \e^{2\pi\i\b k \b x} \,\mathrm d\b x = \delta_{\b 0,\b k},
	\end{align*}
	i.\,e., the same linear system of equations as in \eqref{eq:claim_quadrature_coefficients}.
	We remark that both \cite{Groechenig20} and \cite{Hu09} state the results in the case \mbox{$d=1$}, a generalization to \mbox{$d>1$}, however, is straight-forward.
	\ex
\end{Remark}

By means of Theorem~\ref{Thm:exact_trig_poly} we can now give a condition that guaranties \eqref{eq:Fourier_coeffs_nonequi} being exact for all trigonometric polynomials \eqref{eq:trig_poly_2d} with maximum degree $\b M$.

\begin{Corollary}
	\label{Corollary:cond_exact}
	The two statements \ref{claim_quadrature} and \ref{claim_quadrature_coefficients} in Theorem~\ref{Thm:exact_trig_poly} are not equivalent to property \eqref{eq:exact_reconstr}, since \eqref{eq:claim_quadrature_coefficients} does not imply an exact reconstruction in \eqref{eq:Fourier_coeffs_nonequi}.
	
	However, an augmented variant of \eqref{eq:claim_quadrature_coefficients}, namely
	\begin{align}
		\label{eq:claim_quadrature_coefficients_double}
		\sum_{j=1}^N w_j \,\e^{2\pi\i\b k\b x_j} = \delta_{\b 0,\b k},
		\quad\b k\in{\I_{\b {2M}}},
	\end{align}
	yields an exact reconstruction \mbox{$\hat f_{\b k} = h_{\b k}^{\mathrm{w}}$} in \eqref{eq:Fourier_coeffs_nonequi} for trigonometric polynomials \eqref{eq:trig_poly_2d} with maximum degree $\b M$.
	Additionally, \eqref{eq:claim_quadrature_coefficients_double} implies the matrix equation \mbox{$\b A^* \b W \b A = \b I_{|\I_{\b{M}}|}$}
	with \mbox{$\b A \in \mathbb C ^{N\times |\I_{\b M}|}$} in~\eqref{eq:matrix_A} and the identity matrix \mbox{$\b I_{|\I_{\b{M}}|}$} of size $|\I_{\b{M}}|$.
\end{Corollary}

\begin{proof}
	Utilizing definitions \eqref{eq:Fourier_coeffs_nonequi} and \eqref{eq:trig_poly_2d} we have
	\begin{align*}
		h_{\b k}^{\mathrm{w}}
		&=
		\sum_{j=1}^{N} w_j\, \Bigg( \sum_{\b\ell\in\I_{\b M}} \hat{f}_{\b\ell} \,\e^{2\pi\i\b\ell\b x_j} \Bigg) \,\e^{-2\pi\i \b k \b x_j}
		= 
		\sum_{\b\ell\in\I_{\b M}} \hat{f}_{\b\ell} \,\sum_{j=1}^{N} w_j\, \e^{2\pi\i(\b\ell-\b k)\b x_j} \notag \\
		&= 
		\sum_{\substack{\b\ell\in\I_{\b M} \\ (\b\ell-\b k) \in\I_{\b{M}}}}
		\hat{f}_{\b\ell} \,\sum_{j=1}^{N} w_j\, \e^{2\pi\i(\b\ell-\b k)\b x_j}
		+
		\sum_{\substack{\b\ell\in\I_{\b M} \\ (\b\ell-\b k) \notin\I_{\b{M}}}}
		\hat{f}_{\b\ell} \,\sum_{j=1}^{N} w_j\, \e^{2\pi\i(\b\ell-\b k)\b x_j} ,
		\quad \b k\in\I_{\b M} .
	\end{align*}
	Since \eqref{eq:claim_quadrature_coefficients} only holds for \mbox{$\b k,\b\ell\in\I_{\b M}$} with \mbox{$(\b\ell-\b k) \in\I_{\b{M}}$}, this implies
	\begin{align*}
		h_{\b k}^{\mathrm{w}}
		&= 
		\hat f_{\b k}
		+
		\sum_{\substack{\b\ell\in\I_{\b M} \\ (\b\ell-\b k) \notin\I_{\b{M}}}}
		\hat{f}_{\b\ell} \,\sum_{j=1}^{N} w_j\, \e^{2\pi\i(\b\ell-\b k)\b x_j} ,
		\quad \b k\in\I_{\b M} ,
	\end{align*}
	%
	where for all \mbox{$\b k\in\I_{\b M}\setminus\{\b 0\}$} there exists an \mbox{$\b\ell\in\I_{\b M}$} with \mbox{$(\b\ell-\b k) \in\I_{\b{2M}}\setminus\I_{\b{M}}$}.

	As \mbox{$(\b\ell-\b k) \in\I_{\b{2M}}$} for \mbox{$\b k,\b\ell\in\I_{\b M}$}, the augmented variant \eqref{eq:claim_quadrature_coefficients_double} yields
	\begin{align*}
		h_{\b k}^{\mathrm{w}}
		&= 
		\sum_{\b\ell\in\I_{\b M}} \hat{f}_{\b\ell} \,\sum_{j=1}^{N} w_j\, \e^{2\pi\i(\b\ell-\b k)\b x_j} 
		= 
		\sum_{\b\ell\in\I_{\b M}} \hat{f}_{\b k} \cdot \delta_{\b 0,\b k}
		= 
		\hat f_{\b k},
		\quad \b k\in\I_{\b M}.
	\end{align*}
	Moreover, since \mbox{$\delta_{(\b\ell-\b k),\b 0} = \delta_{\b k,\b\ell}$}, the condition \eqref{eq:claim_quadrature_coefficients_double} implies
	\begin{align*}
		\delta_{\b k,\b\ell}
		=
		\sum_{j=1}^{N} w_j\, \e^{2\pi\i(\b\ell-\b k)\b x_j} 
		= 
		\sum_{j=1}^N \,\e^{-2\pi\i\b k\b x_j} \left( w_j \,\e^{2\pi\i\b\ell\b x_j} \right),
		\quad \b k, \b\ell\in\I_{\b{M}}.
	\end{align*}
	In matrix-vector notation this can be written as \mbox{$\b A^* \b W \b A = \b I_{|\I_{\b{M}}|}$}
	with \mbox{$\b A \in \mathbb C ^{N\times |\I_{\b M}|}$} in~\eqref{eq:matrix_A} and the identity matrix \mbox{$\b I_{|\I_{\b{M}}|}$} of size $|\I_{\b{M}}|$.
	We remark that this matrix equation immediately shows that we have an exact reconstruction of the form~\eqref{eq:exact_reconstr}, since if \mbox{$\b A^* \b W \b A = \b I_{|\I_{\b{M}}|}$} is fulfilled, \eqref{eq:problem_infft} implies that
	\mbox{$\b{\hat f} = \b A^* \b W \b A \b{\hat f} = \b A^* \b W \b f$}.
\end{proof}

\begin{Remark}
	Let \mbox{$f\in L_2(\T^d)$} be an arbitrary \mbox{1-periodic} function \eqref{eq:Fourier_series}.
	Then \eqref{eq:claim_quadrature_coefficients_double} yields
	\begin{align*}
		h_{\b k}^{\mathrm{w}}
		&=
		\sum_{\b\ell\in\Z^d} c_{\b\ell}(f) \,\sum_{j=1}^{N} w_j\, \e^{2\pi\i(\b\ell-\b k)\b x_j} \notag \\
		&= 
		\sum_{\substack{\b\ell\in\Z^d \\ (\b\ell-\b k) \in\I_{\b{2M}}}}
		c_{\b\ell}(f) \,\sum_{j=1}^{N} w_j\, \e^{2\pi\i(\b\ell-\b k)\b x_j}
		+
		\sum_{\substack{\b\ell\in\Z^d \\ (\b\ell-\b k) \notin\I_{\b{2M}}}}
		c_{\b\ell}(f) \,\sum_{j=1}^{N} w_j\, \e^{2\pi\i(\b\ell-\b k)\b x_j} \notag \\
		&= 
		c_{\b k}(f)
		+
		\sum_{\substack{\b\ell\in\Z^d \\ (\b\ell-\b k) \notin\I_{\b{2M}}}}
		c_{\b\ell}(f) \,\sum_{j=1}^{N} w_j\, \e^{2\pi\i(\b\ell-\b k)\b x_j},
		\quad \b k\in\I_{\b M},
	\end{align*}
	i.\,e., for a function \mbox{$f\in L_2(\T^d)$} we only have a good approximation in case the coefficients~$c_{\b\ell}(f)$ are small for \mbox{$\b\ell\notin\I_{\b M}$}, whereas this reconstruction can only be exact for $f$ being a trigonometric polynomial \eqref{eq:trig_poly_2d}.
	\ex
\end{Remark}

\subsubsection{Practical computation in the underdetermined setting \mbox{$|\I_{\b {2M}}| \leq N$} \label{sec:underdet}}

So far, we have seen in Corollary~\ref{Corollary:cond_exact} that an exact solution \mbox{$\b w = (w_j)_{j=1}^N$} to the linear system \eqref{eq:claim_quadrature_coefficients_double} leads to an exact reconstruction formula~\eqref{eq:Fourier_coeffs_nonequi} for all trigonometric polynomials \eqref{eq:trig_poly_2d} with maximum degree~$\b M$.
Therefore, we aim to use this condition \eqref{eq:claim_quadrature_coefficients_double} to numerically find optimal density compensation factors \mbox{$w_j\in\C$}, \mbox{$j=1,\dots,N$}.

Having a closer look at the condition \eqref{eq:claim_quadrature_coefficients_double} we recognize that it can be written as the linear system of equations
\mbox{$\b A_{|\I_{\b {2M}}|}^T \,\b w = \b e_{\b 0}$} with
the matrix \mbox{$\b A_{|\I_{\b{2M}}|} \in \C^{N\times|\I_{\b{2M}}|}$}, cf.~\eqref{eq:matrix_A},
and right side \mbox{$\b e_{\b 0} \coloneqq \left( \delta_{\b 0,\b k} \right)_{\b k\in{\I_{\b{2M}}}}$}.
We remark that in contrast to \mbox{$\b A\in\C^{N \times |\I_{\b{M}}|}$} we now deal with the enlarged matrix \mbox{$\b A_{|\I_{\b{2M}}|} \in \C^{N\times|\I_{\b{2M}}|}$},
such that single matrix operations are more costly.
Nevertheless, Corollary~\ref{Corollary:cond_exact} yields a direct inversion method for \eqref{eq:problem_infft}, where the system \mbox{$\b A_{|\I_{\b {2M}}|}^T \,\b w = \b e_{\b 0}$} needs to be solved only once for fixed points \mbox{$\b x_j\in\T^d$}, \mbox{$j=1,\dots, N$}.
Its solution $\b w$ can then be used to efficiently approximate $\b{\hat f}$ for multiple measurement vectors $\b f$, whereas iterative methods for \eqref{eq:problem_infft} need to solve \mbox{$\b A \b{\hat f} = \b f$} each time.

As already mentioned in~\cite[Sec.~3.1]{Hu09} an exact solution to \eqref{eq:claim_quadrature_coefficients_double} can only be found if \mbox{$|\I_{\b {2M}}| \leq N$}, 
i.\,e., in case \mbox{$\b A_{|\I_{\b {2M}}|}^T \,\b w = \b e_{\b 0}$} is an underdetermined system of equations.
By \cite[Lem.~3.1]{Hu09} this system has at least one solution, which is why we may choose the one with minimal $\ell_2$-norm.
If \mbox{$\mathrm{rank}(\b A_{|\I_{\b {2M}}|}) = {|\I_{\b {2M}}|}$}, then the system \mbox{$\b A_{|\I_{\b {2M}}|}^T \,\b w = \b e_{\b 0}$} is consistent and the unique solution is given by the normal equations of second kind
\begin{align}
	\label{eq:normal_equations_second_kind_double}
	\b A_{|\I_{\b {2M}}|}^T \overline{\b A_{|\I_{\b {2M}}|}} \,\b v = \b e_{\b 0},
	\quad 
	\overline{\b A_{|\I_{\b {2M}}|}} \,\b v = \b w.
\end{align}
More precisely, we may compute the vector $\b v$ using an iterative procedure such as the CG algorithm, such that only matrix multiplications with $\b A_{|\I_{\b {2M}}|}^T$ and $\overline{\b A_{|\I_{\b {2M}}|}}$ are needed.
Since fast multiplication with $\b A_{|\I_{\b {2M}}|}^T$ and $\overline{\b A_{|\I_{\b {2M}}|}}$ can easily be realized by means of an adjoint NFFT (see Algorithm~\ref{alg:nfft*}) and an NFFT (see Algorithm~\ref{alg:nfft}), respectively, computing the solution $\b w$ to \eqref{eq:normal_equations_second_kind_double} is of complexity \mbox{$\mathcal{O}(|\I_{\b {2M}}|\log(|\I_{\b {2M}}|)+N)$}, where
\begin{equation*}
	|\I_{\b {2M}}| = (2M)^d = 2^d M^d = 2^d \, |\I_{\b{M}}|.
\end{equation*}
Thus, in order to receive exact quadrature weights $w_j\in\C$, \mbox{$j=1,\dots, N,$} via \eqref{eq:normal_equations_second_kind_double} we need to satisfy the full rank condition \mbox{$\mathrm{rank}(\b A_{|\I_{\b {2M}}|}) = {|\I_{\b {2M}}|}$}.
In case of a low rank matrix \mbox{$\b A_{|\I_{\b{2M}}|} \in \C^{N\times|\I_{\b{2M}}|}$} for \mbox{$|\I_{\b {2M}}| \leq N$}, we may still use \eqref{eq:normal_equations_second_kind_double} to obtain a least squares approximation to \eqref{eq:claim_quadrature_coefficients_double}.

\subsubsection{Practical computation in the overdetermined setting \mbox{$|\I_{\b {2M}}| > N$}}
In the setting \mbox{$|\I_{\b {2M}}| > N$}, we cannot expect finding an exact solution $\b w$ to \eqref{eq:claim_quadrature_coefficients_double}, 
since we have to deal with an overdetermined system possessing more conditions than variables.
However, we still aim to numerically find optimal density compensation factors \mbox{$w_j\in\C$}, \mbox{$j=1,\dots,N$}, 
by considering a least squares approximation to \eqref{eq:claim_quadrature_coefficients_double} that minimizes
\mbox{$\big\| \b A_{|\I_{\b {2M}}|}^T \,\b w - \b e_{\b 0} \big\|_2$}.
In \cite[Thm.~1.1.2]{Bj96} it was shown that every least squares solution satisfies the normal equations of first kind
\begin{align}
	\label{eq:normal_equations_first_kind_double}
	\overline{\b A_{|\I_{\b {2M}}|}} \,\b A_{|\I_{\b {2M}}|}^T \,\b w = \overline{\b A_{|\I_{\b {2M}}|}} \,\b e_{\b 0} .
\end{align}
By means of definitions of \mbox{$\b A_{|\I_{\b {2M}}|}\in\C^{N \times |\I_{\b{2M}}|}$}, cf.~\eqref{eq:matrix_A}, and \mbox{$\b e_{\b 0} = \left( \delta_{\b 0,\b k} \right)_{\b k\in{\I_{\b{2M}}}}$}
we simplify the right hand side via
\begin{align*}
	\overline{\b A_{|\I_{\b {2M}}|}} \,\b e_{\b 0}
	= 
	\Bigg( \sum_{\b k \in \I_{\b{2M}}} \delta_{\b 0,\b k} \,\e^{-2\pi\i \b k \b x_j} \Bigg)_{j=1}^N
	=
	\b 1_N .
\end{align*}
Since fast multiplication with $\b A_{|\I_{\b {2M}}|}^T$ and $\overline{\b A_{|\I_{\b {2M}}|}}$ can easily be realized by means of an adjoint NFFT (see Algorithm~\ref{alg:nfft*}) and an NFFT (see Algorithm~\ref{alg:nfft}), respectively, the solution $\b w$ to \eqref{eq:normal_equations_first_kind_double} can be computed iteratively by means of the CG algorithm in \mbox{$\mathcal{O}(|\I_{\b {2M}}|\log(|\I_{\b {2M}}|)+N)$} arithmetic operations.
Note that the solution to \eqref{eq:normal_equations_first_kind_double} is only unique if the full rank condition \mbox{$\mathrm{rank}(\b A_{|\I_{\b {2M}}|})=N$} is satisfied, cf.~\cite[p.~7]{Bj96}.
We remark that the computed weight matrix \mbox{$\b W=\diag(\b w)$} can further be used in an iterative procedure as in \cite[Alg.~7.27]{PlPoStTa18} to improve the approximation of $\b{\hat f}$.

The previous considerations lead to the following algorithms.

\begin{algorithm}{Computation of the optimal density compensation factors}
	\label{alg:precompute_density}
	For \mbox{$d,N \in \N$} let \mbox{$\b x_j \in \T^d$}, \mbox{$j=1,\dots,N$}, as well as \mbox{$M\in 2\N$} and \mbox{$\b M \coloneqq M \cdot \b 1_{d}$} be given.
	\begin{enumerate}	
		\item
		Compute \mbox{$|\I_{\b{2M}}| = (2M)^d$}. \hfill \mbox{$\mathcal{O}(1)$}
		\item If \mbox{$|\I_{\b {2M}}| \leq N$}
		\begin{itemize}
			\item[] \hspace{-0.6cm} Compute the solution $\b v$ to 
			\eqref{eq:normal_equations_second_kind_double} iteratively using the NFFT. \\ \phantom{.}
			\hfill \mbox{$\mathcal{O}(|\I_{\b {2M}}|\log(|\I_{\b {2M}}|)+N)$}
			\item[] \hspace{-0.6cm} Compute the solution \mbox{$\b w = \overline{\b A_{|\I_{\b {2M}}|}} \,\b v$}, see~\eqref{eq:normal_equations_second_kind_double}, using an NFFT. \\ \phantom{.}
			\hfill \mbox{$\mathcal{O}(|\I_{\b {2M}}|\log(|\I_{\b {2M}}|)+N)$}
		\end{itemize}
		elseif \mbox{$|\I_{\b {2M}}| > N$}
		\begin{itemize}
			\item[] \hspace{-0.6cm} Compute the solution $\b w$ to \eqref{eq:normal_equations_first_kind_double} iteratively using the NFFT. \\ \phantom{.}
			\hfill \mbox{$\mathcal{O}(|\I_{\b {2M}}|\log(|\I_{\b {2M}}|)+N)$}
		\end{itemize}
		\item Compose $\b W  = \diag(\b w)\in \mathbb C^{N\times N}$.
		\hfill $\mathcal O(N)$
	\end{enumerate}
	\vspace{1ex} 
	\textnormal{\textbf{Output:}} weights matrix $\b W$ \hfill
	\textnormal{\textbf{Complexity:}} \mbox{$\mathcal{O}(|\I_{\b {2M}}|\log(|\I_{\b {2M}}|)+N)$}
\end{algorithm}

\begin{algorithm}{iNFFT -- density compensation approach}
	\label{alg:infft_density}
	For \mbox{$d,N \in \N$} let \mbox{$\b x_j \in \T^d$}, \mbox{$j=1,\dots,N$}, as well as \mbox{$\b f \in \C^N$}, \mbox{$M\in 2\N$} and \mbox{$\b M \coloneqq M \cdot \b 1_{d}$} be given.
	\begin{enumerate}
		\setcounter{enumi}{-1}
		\item Precompute the weights matrix $\b W$ using Algorithm~\ref{alg:precompute_density}.
		\item Compute 
		\mbox{$\b{\tilde h}^{\mathrm{w}} \coloneqq \b D^* \b F^* \b B^* \b W \b f$}, cf. \eqref{eq:reconstr_density_nfft},
		by means of an adjoint NFFT.
	\end{enumerate}
	\vspace{1ex} 
	\textnormal{\textbf{Output:}} $\b{\tilde h}^{\mathrm{w}} \approx \b{\hat f}\in\C^{|\I_{\b{M}}|}$, cf. \eqref{eq:problem_infft}. \hfill
	\textnormal{\textbf{Complexity:}} $\mathcal O(|\I_{\b M}|\log(|\I_{\b M}|) + N)$
\end{algorithm}

\subsection{Bandlimited functions \label{sec:bandlim}}

In some numerical examples, such as in MRI, we are concerned with bandlimited functions \mbox{$f\in L_1(\R^d)\cap C_0(\R^d)$} instead of trigonometric polynomials \mbox{$f\in L_2(\T^d)$} in \eqref{eq:trig_poly_2d}, cf.~\cite{MRI22}.
In the following we show that it is reasonable to consider the inversion problem \eqref{eq:problem_infft} as well as density compensation via \eqref{eq:reconstr_density_nfft} for bandlimited functions as well.

To this end, let \mbox{$f\in L_1(\R^d)\cap C_0(\R^d)$} be a bandlimited function with bandwidth $\b M$, i.\,e., its (continuous) Fourier transform 
\begin{equation}
	\label{eq:inverse_integral}
	\hat f(\b v)
	\coloneqq 
	\int\limits_{\R^d} 
	f(\b x)\,\e^{-2\pi\i \b v\b x}\,\mathrm d\b x,
	\quad \b v\in\R^d,
\end{equation}
is supported on \mbox{$\left[-\tfrac M2,\tfrac M2\right)^d$}.
Utilizing this fact, we have \mbox{$\hat f \in L_1(\R^d)$}
and thus by the Fourier inversion theorem \cite[Thm.~2.10]{PlPoStTa18} the inverse Fourier transform of $f$ can be written as
\begin{equation}
	\label{eq:forward_integral}
	f(\b x) 
	= 
	\int\limits_{\R^d} \hat f(\b v)\,\e^{2\pi\i \b v\b x}\,\mathrm d\b v
	= 
	\int\limits_{\left[-\frac M2,\frac M2\right)^d} \hat f(\b v)\,\e^{2\pi\i \b v\b x}\,\mathrm d\b v,
	\quad \b x\in \R^d.
\end{equation}
Analogous to \eqref{eq:Fourier_coeffs_equi}, the approximation using equispaced quadrature points \mbox{$\b k \in \I_{\b M}$} yields
\begin{equation}
	\label{eq:quadrature_forward}
	f(\b x) 
	\approx
	\frac{\left(\frac M2-(-\frac M2)\right)^d}{|\I_{\b M}|} \sum_{\b k \in \I_{\b M}} \hat f(\b k)\,\e^{2\pi\i \b k\b x}
	=
	\sum_{\b k \in \I_{\b M}} \hat f(\b k)\,\e^{2\pi\i \b k\b x},
	\quad \b x\in \R^d,
\end{equation}
such that evaluation at the given nonequispaced points \mbox{$\b x_j\in \left[ -\tfrac 12, \tfrac 12 \right)^d$}, \mbox{$j=1,\dots, N$}, leads to
\begin{equation*}
	f(\b x_j) 
	\approx
	\sum_{\b k \in \I_{\b M}} \hat f(\b k)\,\e^{2\pi\i \b k\b x_j}.
\end{equation*}
By means of the definition \eqref{eq:matrix_A} of the matrix \mbox{$\b A\in\C^{N\times|\I_{\b{M}}|}$} this can be written as
\mbox{$\b f \approx \b A \b{\hat f}$},
where we used the notation \mbox{$\b{\hat f} \coloneqq (\hat f(\b k))_{\b k\in\I_{\b M}}$} in this setting.
Thus, also for band\-limited functions $f$ its evaluations at points \mbox{$\b x_j$} can be approximated in the form~\eqref{eq:nfft}, such that it is reasonable to consider the inversion problem \eqref{eq:problem_infft} for bandlimited functions as well.

Considering \eqref{eq:inverse_integral} we are given an exact formula for the evaluation of the Fourier transform $\hat f$.
However, in practical applications, such as MRI, this is only a hypothetical case, since $f$ cannot be sampled on whole $\R^d$, cf.~\cite{MRI22}.
Due to a limited coverage of space by the acquisition, the function $f$ is typically only known on a bounded domain, w.l.o.g. for $\b x\in \left[ -\tfrac 12, \tfrac 12 \right)^d$.
Thus, we have to deal with the approximation
\begin{equation}
	\label{eq:approx_inverse_integral}
	\hat f(\b v)
	\approx 
	\int\limits_{\left[ -\tfrac 12, \tfrac 12 \right)^d} 
	f(\b x)\,\e^{-2\pi\i \b v\b x}\,\mathrm d\b x,
	\quad \b v\in\left[-\tfrac M2,\tfrac M2\right)^d.
\end{equation}
Using the nonequispaced quadrature rule in \eqref{eq:Fourier_coeffs_nonequi}, we find that evaluation at uniform grid points \mbox{$\b k\in\I_{\b M}$} can be approximated via
\begin{equation*}
	\hat f(\b k)
	\approx
	\tilde h(\b k) 
	\coloneqq 
	\sum_{j=1}^{N} 
	w_j\, f(\b x_j)\,\e^{-2\pi\i \b k \b x_j},
	\quad \b k\in\I_{\b M} .
\end{equation*}
This is to say, equispaced samples of the Fourier transform of a bandlimited function may be approximated in the same form \mbox{$(\hat f(\b k))_{\b k\in\I_{\b M}} = \b{\hat f} \approx \b{h}^{\mathrm{w}} \coloneqq \b A^* \b W \b f$} as in \eqref{eq:Fourier_coeffs_nonequi},
where we used the notation \mbox{$\b{h}^{\mathrm{w}} \coloneqq (\tilde h(\b k))_{\b k\in\I_{\b M}}$} in this setting.
Moreover, we extend this approximation onto the whole interval 
\mbox{$\left[-\tfrac M2,\tfrac M2\right)^d$}, i.\,e., we consider 
\begin{equation}
	\label{eq:quadrature_inverse}
	\hat f(\b v)
	\approx
	\tilde h(\b v) 
	\coloneqq 
	\sum_{j=1}^{N} 
	w_j\, f(\b x_j)\,\e^{-2\pi\i \b v \b x_j},
	\quad \b v\in\left[-\tfrac M2,\tfrac M2\right)^d.
\end{equation}

So all in all, we have seen that it is reasonable to study the inversion problem \eqref{eq:problem_infft} and the associated density compensation via \eqref{eq:reconstr_density_nfft} for bandlimited functions as well.
Analogous to Section~\ref{sec:trig_poly} we now aim to find a numerical method for computing suitable weights \mbox{$w_j\in\C$}, \mbox{$j=1,\dots,N$}, such that the reconstruction formula~\eqref{eq:quadrature_inverse} is preferably exact.
To this end, we have a closer look at \eqref{eq:quadrature_inverse} being exact and start with analogous considerations as in \autoref{Thm:exact_trig_poly}.

\begin{Theorem}
	\label{Thm:exact_bandlim}
	Let a bandwidth \mbox{$\b M\in\N^d$}, nonequispaced points \mbox{$\b x_j \in \left[-\tfrac 12,\tfrac 12\right)^d$} as well as quadrature weights \mbox{$w_j\in\C$}, \mbox{$j=1,\dots, N,$} be given. 
	Then an exact reconstruction formula~\eqref{eq:quadrature_inverse} for bandlimited functions \mbox{$f\in L_1(\R^d)\cap C_0(\R^d)$} with bandwidth $\b M$, i.\,e., 
	\begin{equation}
		\label{eq:exact_reconstr_bandlim}
		\hat f(\b v) = \tilde{h}(\b v) = \sum_{j=1}^N w_j f(\b x_j) \,\e^{-2\pi\i\b v \b x_j},
		\quad \b v\in\left[-\tfrac M2,\tfrac M2\right)^d ,
	\end{equation}
	implies that the quadrature rule
	\begin{align*}
		\int_{\R^d} f(\b x) \,\mathrm d \b x = \sum_{j=1}^N w_j f(\b x_j)
	\end{align*}
	is exact for all bandlimited functions \mbox{$f\in L_1(\R^d)\cap C_0(\R^d)$} with bandwidth $\b M$.
\end{Theorem}

\begin{proof} 
	By \eqref{eq:inverse_integral} the assumption \eqref{eq:exact_reconstr_bandlim} can be written as 
	\begin{equation}
		\label{eq:reconstr_bandlim}
		\int\limits_{\R^d} f(\b x)\,\e^{-2\pi\i \b v\b x}\,\mathrm d\b x
		= \hat f(\b v) 
		= \sum_{j=1}^N w_j f(\b x_j) \,\e^{-2\pi\i\b v \b x_j},
		\quad \b v\in\left[-\tfrac M2,\tfrac M2\right)^d .
	\end{equation}
	Especially, for \mbox{$\b v = \b 0$} evaluation of \eqref{eq:reconstr_bandlim} yields the assertion 
	\begin{equation*}
		\int\limits_{\R^d} f(\b x)\,\mathrm d\b x
		= \hat f(\b 0) 
		= \sum_{j=1}^N w_j f(\b x_j) .
	\end{equation*}
\end{proof}

However, 
in contrast to \autoref{Thm:exact_trig_poly}, this \autoref{Thm:exact_bandlim} does not yield an explicit condition for computing suitable weights \mbox{$w_j\in\C$}, \mbox{$j=1,\dots,N$}.
To derive a numerical procedure anyway, we generalize the notion of an exact reconstruction $\tilde h$ of $f$ and have a look at the theory of tempered distributions.
To this end, let \mbox{$\mathscr S(\R^d)$} be the Schwartz space of rapidly decaying functions, cf.~\cite[Sec.~4.2.1]{PlPoStTa18}.
The tempered Dirac distribution $\delta$ shall be defined by 
\mbox{$\langle \delta, \varphi \rangle \coloneqq \int_{\R^d} \varphi(\b v) \,\delta(\b v) \,\mathrm d\b v = \varphi(\b 0)$} for all \mbox{$\varphi\in\mathscr S(\R^d)$}, cf.~\cite[Ex.~4.36]{PlPoStTa18}.
For a slowly increasing function \mbox{$f \colon \R^d \to \C$} satisfying \mbox{$|f(\b x)| \leq c (1+\|\b x\|_2)^n$} almost everywhere with \mbox{$c>0$} and \mbox{$n\in\N_0$}, the induced distribution $T_f$ shall be defined by 
\mbox{$\langle T_{f}, \varphi \rangle \coloneqq \int_{\R^d} \varphi(\b x) \,f(\b x) \,\mathrm d\b x$} for all \mbox{$\varphi\in\mathscr S(\R^d)$}.
For a detailed introduction to the topic we refer to \cite[Sec.~4.2.1 and Sec.~4.3]{PlPoStTa18}.

Then the following property can be shown.

\begin{Theorem}
	\label{Thm:exact_delta}
	Let nonequispaced points \mbox{$\b x_j \in \left[-\frac 12, \frac 12\right)^d$}, \mbox{$j=1,\dots, N,$} and quadrature weights $w_j \in \C$ be given. 
	Further let $T_f$ be the distribution induced by some bandlimited function \mbox{$f\in L_1(\R^d)\cap C_0(\R^d)$} with bandwidth $\b M$.
	Then 
	\begin{align}
		\label{claim_delta}
		\langle \delta, \varphi \rangle = \langle T_\xi, \varphi \rangle,
		\quad 
		\varphi\in\mathscr S(\R^d),
	\end{align}
	%
	with
	%
	\begin{align}
		\label{eq:claim_delta}
		\xi(\b v) \coloneqq \sum_{j=1}^N w_j \,\e^{2\pi\i\b v\b x_j},
		\quad \b v\in\R^d,
	\end{align}
	%
	%
	implies
	\begin{align}
		\label{claim_exact_reconstr_delta}
		\langle \hat T_{f}, \varphi \rangle = \langle T_{\tilde h}, \varphi \rangle,
		\quad 
		\varphi\in\mathscr S(\R^d),
	\end{align}
	%
	with the function $\tilde h$ defined in \eqref{eq:quadrature_inverse}.
	%
	%
	%
\end{Theorem}

\begin{proof} 

	Using the definition of the function $\tilde h$ in \eqref{eq:quadrature_inverse} as well as the fact that the inversion formula~\eqref{eq:forward_integral} holds for all \mbox{$\b x \in\R^d$}, we have
	\begin{align*}
		\langle T_{\tilde h}, \varphi \rangle
		&=
		\int_{\R^d} \varphi(\b v) \,\sum_{j=1}^N w_j  f(\b x_j) \, \e^{-2\pi\i\b v\b x_j} \,\mathrm{d}\b v \\
		&=
		\int_{\R^d} \varphi(\b v) \,\sum_{j=1}^N w_j  
		\left( \int_{\R^d} \hat f(\b u)\,\e^{2\pi\i \b u\b x_j}\,\mathrm d\b u \right) 
		\, \e^{-2\pi\i\b v\b x_j} \,\mathrm{d}\b v \\
		&=
		-\int_{\R^d} \hat f(\b u)
		\int_{\R^d} \varphi(\b u-\b v) \,\sum_{j=1}^N w_j \,\e^{2\pi\i \b v\b x_j} \,\mathrm{d}\b v \,\mathrm d\b u .
	\end{align*}
	Hence, by \eqref{claim_delta} this implies
	\begin{align*}
		\langle T_{\tilde h}, \varphi \rangle
		&=
		-\int_{\R^d} \hat f(\b u)
		\int_{\R^d} \varphi(\b u-\b v) \,\delta(\b v) \,\mathrm{d}\b v \,\mathrm d\b u  \\
		&=
		\int_{\R^d} \hat f(\b u)
		\int_{\R^d} \varphi(\b v) \,\delta(\b u-\b v) \,\mathrm{d}\b v \,\mathrm d\b u 
		=
		\int_{\R^d} \hat f(\b u) \, \varphi(\b u) \,\mathrm d\b u .
	\end{align*}
\end{proof}

Considering the property \eqref{claim_exact_reconstr_delta}, we remark that this indeed states an exact reconstruction 
\mbox{$\hat f = \tilde h$} in the sense of tempered distributions, as distinct from \eqref{eq:exact_reconstr_bandlim}.
Since it is known by Corollary~\ref{Corollary:cond_exact} that the condition \eqref{eq:claim_quadrature_coefficients_double} yields an exact reconstruction for trigonometric polynomials,
we aim to use this result to compute suitable weights \mbox{$w_j\in\C$}, \mbox{$j=1,\dots,N$}, for bandlimited functions as well.
To this end, suppose we have \eqref{eq:claim_quadrature_coefficients_double}, i.\,e.,
\begin{align*}
	\sum_{j=1}^N w_j \,\e^{2\pi\i\b k\b x_j} 
	=
	\delta_{\b 0,\b k},
	\quad \b k \in \I_{\b{2M}}.
\end{align*}
Then this yields
\begin{align}
	\label{eq:claim_delta_double_discretized}
	\varphi(\b 0) 
	=
	\sum_{\b k \in \I_{\b{2M}}} \varphi(\b k) \,\sum_{j=1}^N w_j \,\e^{2\pi\i\b k\b x_j} ,
	\quad \varphi\in\mathscr S(\R^d) .
\end{align}
Having a look at \autoref{Thm:exact_delta}, an exact reconstruction \eqref{claim_exact_reconstr_delta} is implied by \eqref{claim_delta}, i.\,e.,
\begin{align*}
	\varphi(\b 0) = \langle \delta, \varphi \rangle = \langle T_\xi, \varphi \rangle
	=
	\int_{\R^d} \varphi(\b v) \,\sum_{j=1}^N w_j \,\e^{2\pi\i\b v\b x_j} \,\mathrm{d}\b v ,
	\quad \varphi\in\mathscr S(\R^d) .
\end{align*}
Thus, the property \eqref{eq:claim_delta_double_discretized} that is fulfilled by \eqref{eq:claim_quadrature_coefficients_double} could be interpreted as an equispaced quadrature of \eqref{claim_delta} at integer frequencies \mbox{$\b k\in\I_{\b{2M}}$}. 

\begin{Remark}
	\label{Remark:implicit_truncation}
	We remark that for deriving the quadrature rule \eqref{eq:claim_delta_double_discretized} from \eqref{claim_delta}, we implicitly truncate the integral bounds in \eqref{claim_delta} as
	\begin{align*}
		\langle T_\xi, \varphi \rangle
		=
		\int_{\R^d} \varphi(\b v) \,\sum_{j=1}^N w_j \,\e^{2\pi\i\b v\b x_j} \,\mathrm{d}\b v
		&\approx
		\int\limits_{[-M,M)^d} \varphi(\b v) \,\sum_{j=1}^N w_j \,\e^{2\pi\i\b v\b x_j} \,\mathrm{d}\b v \notag \\
		&=
		\int_{\R^d} \varphi(\b v) \,\sum_{j=1}^N w_j \,\e^{2\pi\i\b v\b x_j} \,\chi_{\left[-M,M\right)^d}(\b v) \,\mathrm{d}\b v ,
	\end{align*}
	i.\,e., instead of \eqref{eq:claim_delta} we rather deal with a distribution induced by
	\begin{align}
		\label{eq:claim_delta_bandlim}
		\tilde\xi(\b v) \coloneqq \xi(\b v) \,\chi_{\left[-M,M\right)^d}(\b v), \quad \b v\in\R^d.
	\end{align}
	However, an analogous implication as in \autoref{Thm:exact_delta} cannot be shown when using the function \eqref{eq:claim_delta_bandlim} instead of \eqref{eq:claim_delta}.
	\ex
\end{Remark}

As seen in Remark~\ref{Remark:implicit_truncation}, our numerical method for computing suitable weights \mbox{$w_j\in\C$}, \mbox{$j=1,\dots,N$}, can also be derived by means of a quadrature formula applied to the property
\mbox{$\langle \delta, \hat\varphi \rangle = \langle T_{\tilde\xi}, \hat\varphi \rangle$} 
with $\tilde\xi$ defined in \eqref{eq:claim_delta_bandlim}.
Having a closer look at this property, the following equivalent characterization can be shown.

\begin{Theorem}
	\label{Thm:exact_sinc_bandlim}
	Let a bandwidth $\b M \in \N^d$, nonequispaced points \mbox{$\b x_j \in \left[-\frac 12, \frac 12\right)^d$} and quadrature weights $w_j \in \C$, \mbox{$j=1,\dots, N,$} be given. 
	Then the following two statements are equivalent. 
	%
	\begin{enumerate}[label=(\roman*), resume]
		\item For all \mbox{$\varphi\in\mathscr S(\R^d)$} we have \mbox{$\langle \delta, \hat\varphi \rangle = \langle T_{\tilde\xi}, \hat\varphi \rangle$} with $\tilde\xi$ defined in \eqref{eq:claim_delta_bandlim}.
		\label{claim_delta_bandlim}
		\item We have \mbox{$\langle 1, \varphi \rangle = \langle T_{\psi}, \varphi \rangle$} for all \mbox{$\varphi\in\mathscr S(\R^d)$}, where
		\label{claim_sinc_bandlim}
		\begin{align*}
			\psi(\b x) \coloneqq \sum_{j=1}^N w_j \cdot |\I_{\b{2M}}| \,\sinc\left(2M\pi\left(\b x_j - \b x\right)\right),
			\quad \b x\in\R^d ,
		\end{align*}
		with the $d$-variate \mbox{$\sinc$} function
		\mbox{$\sinc(\b x) \coloneqq \prod_{t=1}^d \mathrm{sinc}(x_t)$} and
		$$
		\mathrm{sinc}(x) \coloneqq \left\{ \begin{array}{ll}  \frac{\sin x}{x} & \quad x \in \mathbb R \setminus \{0\}\,, \\ [1ex]
			1 & \quad x = 0\,. \end{array} \right.
		$$
	\end{enumerate}
\end{Theorem}

\begin{proof} 

	By the definition \mbox{$\langle \hat T, \varphi \rangle = \langle T, \hat\varphi \rangle$}
	of the Fourier transform of a tempered distribution \mbox{$T\in\mathscr{S}'(\R^d)$} we have 
	\mbox{$\langle 1, \varphi \rangle = \langle \delta, \hat\varphi \rangle$}, cf.~\cite[Ex.~4.46]{PlPoStTa18}.
	Moreover, the distribution induced by \eqref{eq:claim_delta_bandlim} can be rewritten using the Fourier transform \eqref{eq:inverse_integral} as
	\begin{align}
		\label{eq:distrib_delta_bandlim}
		\langle T_{\tilde\xi}, \hat\varphi \rangle
		&=
		\int_{\R^d} \hat\varphi(\b v) \,\tilde\xi(\b v) \,\mathrm{d}\b v
		=
		\sum_{j=1}^N w_j \,\int\limits_{\left[-M,M\right)^d} \hat\varphi(\b v) \,\e^{2\pi\i\b v\b x_j} \,\mathrm{d}\b v \notag \\
		&=
		\sum_{j=1}^N w_j \int\limits_{\left[-M,M\right)^d} 
		\left( \int_{\R^d} \varphi(\b x) \,\e^{-2\pi\i\b v\b x} \,\mathrm{d}\b x \right)
		\e^{2\pi\i\b v\b x_j} \,\mathrm{d}\b v \notag \\
		&=
		\sum_{j=1}^N w_j \int_{\R^d} \varphi(\b x) 
		\int\limits_{\left[-M,M\right)^d} \e^{2\pi\i\b v(\b x_j-\b x)} \,\mathrm{d}\b v
		\,\mathrm{d}\b x . 
	\end{align}
	The inner integral can be determined by
	\begin{align*}
		\int_{\left[-M,M\right)^d} \e^{2\pi\i\b v(\b x_j-\b x)} \,\mathrm{d}\b v = |\I_{\b{2M}}| \,\sinc\left(2M\pi\left(\b x_j - \b x\right)\right), 
	\end{align*}
	such that \eqref{eq:distrib_delta_bandlim} shows the equality
	\mbox{$\langle T_{\tilde\xi}, \hat\varphi \rangle = \langle T_{\psi}, \varphi \rangle$}.
	Hence, the assertions \ref{claim_delta_bandlim} and \ref{claim_sinc_bandlim} are equivalent.
\end{proof}

Thus, since the statements \ref{claim_delta_bandlim} and \ref{claim_sinc_bandlim} of Theorem~\ref{Thm:exact_sinc_bandlim} are equivalent,
one could also consider the property \ref{claim_sinc_bandlim} for deriving a numerical method to compute suitable weights \mbox{$w_j\in\C$}, \mbox{$j=1,\dots,N$}.
To this end, we have a closer look at 
\begin{align}
	\label{eq:claim_sinc_bandlim_double}
	\int_{\R^d} \varphi(\b x) \,\mathrm{d}\b x
	=
	\langle 1, \varphi \rangle = \langle T_{\psi}, \varphi \rangle
	=
	\int_{\R^d} \varphi(\b x) \,\sum_{j=1}^N w_j  \, \cdot |\I_{\b{2M}}| \,\sinc\left(2M\pi\left(\b x_j - \b x\right)\right) \,\mathrm{d}\b x 
\end{align}
for all \mbox{$\varphi\in\mathscr S(\R^d)$}.
Due to the integrals on both sides of \eqref{eq:claim_sinc_bandlim_double}, we need to discretize twice
and therefore use the same quadrature rule on both sides of \eqref{eq:claim_sinc_bandlim_double}.
For better comparability to \eqref{eq:claim_delta_double_discretized} we utilize the same number \mbox{$|\I_{\b{2M}}|$} of equispaced quadrature points \mbox{$\b y_{\b \ell} \coloneqq (\b{2M})^{-1} \odot\b\ell$}, \mbox{$\b\ell \in \I_{\b{2M}}$}, as in \eqref{eq:claim_delta_double_discretized}, i.\,e., we consider
\begin{align*}
	\sum_{\b\ell \in \I_{\b{2M}}} \varphi(\b y_{\b \ell})
	=
	\sum_{\b\ell \in \I_{\b{2M}}} \varphi(\b y_{\b \ell}) \,\sum_{j=1}^N w_j  \cdot |\I_{\b{2M}}| \,\sinc\left(2M\pi\left(\b x_j - \b y_{\b \ell}\right)\right) .
\end{align*}
In order that this applies for all \mbox{$\varphi\in\mathscr S(\R^d)$}, we need to satisfy 
\begin{align}
	\label{eq:lgs_claim_sinc_bandlim}
	1
	=
	\sum_{j=1}^N w_j \cdot |\I_{\b{2M}}| \,\sinc\left(2M\pi\left(\b x_j - \b y_{\b \ell}\right)\right) ,
	\quad \b\ell \in \I_{\b{2M}} ,
\end{align}
i.\,e., one could also compute weights \mbox{$w_j\in\C$}, \mbox{$j=1,\dots,N$}, as a least squares solution to the linear system of equations \eqref{eq:lgs_claim_sinc_bandlim}.
Hence, it merely remains the comparison of the two computation schemes.

\begin{Remark}
	Since we derived discretizations out of both statements of \autoref{Thm:exact_sinc_bandlim}, we examine if also the two linear systems \eqref{eq:claim_quadrature_coefficients_double} and \eqref{eq:lgs_claim_sinc_bandlim} are related. 
	Considering the statements in \autoref{Thm:exact_sinc_bandlim} we notice that in some sense they are the Fourier transformed versions of each other.
	To this end, we need to Fourier transform one of the linear systems for better comparability.
	More precisely, we apply an iFFT of length~\mbox{$|\I_{\b{2M}}|$}, cf.~\eqref{eq:idft}, to both sides of equation~\eqref{eq:lgs_claim_sinc_bandlim}.
	Since the left side transforms to
	\begin{align*}
		\sum_{\b\ell\in\I_{\b{2M}}} 1 \cdot \e^{2\pi\i\b k\b y_{\b \ell}}
		=
		|\I_{\b{2M}}| \cdot \delta_{\b 0,\b k},
		\quad \b k\in\I_{\b{2M}},
	\end{align*}
	we obtain the transformed system
	\begin{align}
		\label{eq:lgs_sinc_transformed}
		\delta_{\b 0,\b k}
		=
		\sum_{j=1}^N w_j 
		\sum_{\b\ell\in\I_{\b{2M}}} \,\sinc\left(2M\pi\left(\b x_j - \b y_{\b \ell}\right)\right)
		\,\e^{2\pi\i\b k\b y_{\b \ell}},
		\quad \b k\in\I_{\b{2M}}.
	\end{align}
	Comparing this linear system of equations to \eqref{eq:claim_quadrature_coefficients_double}, we recognize an identical structure. 
	Hence, we have a closer look at the connection between the expressions \mbox{$\e^{2\pi\i\b k \b x_j}$} and \mbox{$\sum_{\b\ell\in\I_{\b{2M}}} \sinc\left(2M\pi\left(\b x_j - \b y_{\b \ell}\right)\right)\,\e^{2\pi\i\b k\b y_{\b \ell}}$}.
	
	For this purpose, we consider the function \mbox{$f(\b t) = \e^{2\pi\i \b t \b x}$}, \mbox{$\b t\in \left[-M,M\right)^d$}, for fixed \mbox{$\b x\in\C^d$}.
	By means of \mbox{$\tilde f(\b t) \coloneqq \sum_{\b k \in \Z^d} f(\b t+{2M}\b k)$} we extend it into a $({2M})$-periodic function.
	This periodized version then possesses the Fourier coefficients
	\begin{align*}
		c_{\b\ell}(\tilde f) 
		&=
		\frac{1}{|\I_{\b{2M}}|}
		\int\limits_{\left[-M,M\right)^d}
		f(\b t) \,\e^{-2\pi\i \b t\b y_{\b \ell}}\,\mathrm d\b t \\
		&=
		\frac{1}{|\I_{\b{2M}}|}
		\int\limits_{\left[-M,M\right)^d}
		\e^{2\pi\i \b t(\b x-\b y_{\b \ell})}\,\mathrm d\b t 
		=
		\sinc\left(2M\pi\left(\b x - \b y_{\b \ell}\right)\right),
		\quad \b\ell\in\Z^d ,
	\end{align*}
	cf.~\eqref{eq:Fourier_coeffs},
	i.\,e., the Fourier expansion of \new{\mbox{$\tilde f(\b t)$}}, \mbox{$\b t\in \left[-M,M\right)^d$}, for fixed $\b x$ is given by 
	\begin{align*}
		\e^{2\pi\i\b t \b x}
		=
		\sum_{\b\ell\in\Z^d} \e^{2\pi\i\b t \b y_{\b \ell}} \,\sinc\left(2M\pi\left(\b x - \b y_{\b \ell}\right)\right),
		\quad \b x\in\C^d,
	\end{align*}
	cf.~\eqref{eq:Fourier_series}.
	Since \new{\mbox{$\tilde f(\b t)$}} is continuous and piecewise differentiable, this Fourier series converges absolutely and uniformly, cf.~\cite[\new{Ex.~1.22}]{LuBo}.
	Thereby, we may consider the point evaluations at \mbox{$\b x = \b x_j$}, \mbox{$j=1,\dots,N$}, and \mbox{$\b t = \b k\in\I_{\b{2M}}$}, such that we obtain the representation
	\begin{align}
		\label{eq:exp_sinc_basis}
		\e^{2\pi\i\b k \b x_j}
		=
		\sum_{\b\ell\in\Z^d} \e^{2\pi\i\b k \b y_{\b \ell}} \,\sinc\left(2M\pi\left(\b x_j - \b y_{\b \ell}\right)\right).
	\end{align}
	Thus, we recognize that \eqref{eq:lgs_sinc_transformed} is a truncated version of \eqref{eq:exp_sinc_basis}. 
	In other words, this implies that the linear system \eqref{eq:claim_quadrature_coefficients_double} is equivalent to a discretization of \eqref{eq:claim_sinc_bandlim_double} incorporating infinitely many points \mbox{$\b y_{\b\ell}\in\R^d$} in \eqref{eq:lgs_claim_sinc_bandlim}.
	\ex
\end{Remark}

\begin{Remark}
	For bandlimited functions several fast evaluation methods including the $\sinc$ function are known.	
	The classical sampling theorem of Shannon-Whittaker-Kotelnikov, see \cite{Whittaker, Shannon49, Kotelnikov},
	states that any bandlimited function \mbox{$f\in L_1(\R^d)\cap C_0(\R^d)$} with maximum bandwidth $\b M$ can be recovered from its uniform samples \mbox{$f(\b L^{-1} \odot {\b \ell})$}, \mbox{$\b\ell\in\Z^d$}, with \mbox{$L \geq M$}, \mbox{$\b L \coloneqq L\cdot\b 1_d$}, and we have
	\begin{align}
		\label{eq:shannon}
		f(\b x) 
		= 
		\sum_{\b\ell\in\Z^d} f(\b L^{-1} \odot {\b \ell}) 
		\,\sinc\left(L\pi\left(\b x - \b L^{-1} \odot {\b \ell}\right)\right) ,
		\quad \b x\in\R^d .
	\end{align}
	Since the practical use of this sampling theorem is limited due to the infinite number of samples, which is impossible in practice, and the very slow decay of the sinc function, 
	various authors such as \cite{Q03, StTa06, MXZ09, LZ16, KiPoTa22} considered the regularized Shannon sampling formula with localized sampling
	\begin{align}
		\label{eq:shannon_regularized}
		f(\b x) 
		\approx
		\sum_{\b\ell\in\Z^d} f(\b L^{-1} \odot {\b \ell}) 
		\,\sinc\left(L\pi\left(\b x - \b L^{-1} \odot {\b \ell}\right)\right)
		\,\varphi_m\left(\b x - \b L^{-1} \odot {\b \ell}\right) ,
		\quad \b x\in\R^d ,
	\end{align}
	instead.
	Here \mbox{$\varphi_m \colon \R^d\to[0,1]$} is a compactly supported window function with truncation parameter \mbox{$m\in\N\setminus\{1\}$}, such that
	for \mbox{$\varphi_m$} with small support the direct evaluation of \eqref{eq:shannon_regularized} is efficient, see \cite{KiPoTa22} for the relation to the NFFT window functions.
	
	On the other hand, in the one-dimensional setting a fast $\sinc$ transform was introduced in \cite{KiPoTa21}, which is based on the Clenshaw-Curtis quadrature 
	\begin{align*}
		\sinc\left(L\pi\left(x - \tfrac{\ell}{L}\right)\right)
		\approx
		\sum_{k=0}^{n} w_{k} \,\e^{-\pi\i L(x-\frac{\ell}{L})z_k}
	\end{align*}
	using Chebyshev points \mbox{$z_k=\cos(\frac{k\pi}{n})\in[-1,1]$}, \mbox{$k=0,\dots,n$}, and corresponding Clenshaw-Curtis weights $w_k>0$.
	Thereby, sums of the form 
	\begin{align*}
		h(x) 
		= 
		\sum_{\ell\in\I_T} f\big(\tfrac{\ell}{L}\big) \,\sinc\left(L\pi\left(x - \tfrac{\ell}{L}\right)\right)
		\approx
		\sum_{k=0}^{n} w_{k} \,
		\bigg( \sum_{\ell\in\I_T} f\big(\tfrac{\ell}{L}\big) \,\e^{\pi\i \ell z_k} \bigg)
		\,\e^{-\pi\i Lxz_k}
	\end{align*}
	with uniform truncation parameter \mbox{$T\in 2\N$}, can efficiently be approximated by means of fast Fourier transforms.
	More precisely, for the term in brackets one may utilize an NFFT, cf.~\eqref{eq:nfft}.
	Then the resulting outer sum can be computed using an NNFFT, also referred to as NFFT of type III, see \cite{duro93, ElSt, LeGr05} or \cite[pp.~394–397]{PlPoStTa18}.
	\ex
\end{Remark}

\subsection{General error bound \label{sec:error_bound}}

In this section we summarize our previous findings by presenting a general error bound on density compensation factors computed by means of \eqref{eq:claim_quadrature_coefficients_double}, that applies to
trigonometric polynomials, 1-periodic functions \mbox{$f\in L_2(\T^d)\cap C(\T^d)$} and
band\-limited functions \mbox{$f\in L_1(\R^d)\cap C_0(\R^d)$} as well.

\begin{Theorem}
	\label{Thm:error_est_lgs}
	Let \mbox{$p,q \in \{1,2,\infty\}$} with \mbox{$\frac 1p + \frac 1q = 1$}.
	For given \mbox{$d, N\in\N$}, \mbox{$\b M\in(2\N)^d$} and nonequispaced points \mbox{$\b x_j \in \left[-\tfrac 12,\tfrac 12\right)^d$}, \mbox{$j=1,\dots,N$}, 
	let \mbox{$\b A \in \C^{N\times|\I_{\b{M}}|}$} be the nonequispaced Fourier matrix in \eqref{eq:matrix_A}.
	Further assume we can compute density compensation factors \mbox{$\b W = \diag\left(w_j\right)_{j=1}^N\in\C^{N\times N}$} by means of Algorithm~\ref{alg:precompute_density}, such that 
	\begin{align}
		\label{eq:residual_lgs}
		\sum_{j=1}^N w_j \,\e^{2\pi\i\b k\b x_j}
		=
		\delta_{\b 0,\b k} + \varepsilon_{\b k},
		\quad \b k\in\I_{\b{2M}},
	\end{align}
	with small \mbox{$\varepsilon_{\b k} \in \R$} for all \mbox{$\b k\in\I_{\b{2M}}$}. 
	
	Then there exists an \mbox{$\varepsilon \geq 0$} such that the corresponding density compensation procedure with \mbox{$\b W = \mathrm{diag} (w_j)_{j=1}^N$} satisfies the following error bounds.
	\begin{enumerate}[label=(\roman*)]
		\item 
		\label{error_est_trig_poly}
		For any trigonometric polynomial \mbox{$f\in L_2(\T^d)$} of degree~$\b M$ given in \eqref{eq:trig_poly_2d} we have
		\begin{align}
			\label{eq:error_est_trig_poly}
			\big\| \b{\hat f} - \b A^* \b W \b f \big\|_p
			&\leq
			|\I_{\b{M}}| \, \varepsilon  \cdot \big\| \b{\hat f} \big\|_p ,
		\end{align}
		where \mbox{$\b{\hat f}\coloneqq(\hat f_{\b k})_{\b k \in \I_{\b M}}$} are the coefficients given in \eqref{eq:trig_poly_2d}.
		\item 
		\label{error_est_periodic}
		For any 1-periodic function \mbox{$f\in L_2(\T^d)\cap C(\T^d)$} we have
		\begin{align}
			\label{eq:error_est_periodic}
			\big\| \b{\hat f} - \b A^* \b W \b f \big\|_p
			&\leq
			|\I_{\b{M}}| \, \varepsilon  \cdot \big\| \b{\hat f} \big\|_p
			+
			(N\,|\I_{\b{M}}|)^{1/p} \, \|\b w\|_{q} \cdot \|f-p_{\b M}\|_{C(\T^d)} ,
		\end{align}
		where \mbox{$\b{\hat f} \coloneqq (c_{\b k}(f))_{\b k\in\I_{\b M}}$} are the first $|\I_{\b{M}}|$ coefficients given in \eqref{eq:Fourier_series} and 
		$p_{\b{M}}$ is the best approximating trigonometric polynomial of degree $\b M$ of $f$.
		\item 
		\label{error_est_bandlim}
		For any bandlimited function \mbox{$f\in L_1(\R^d)\cap C_0(\R^d)$} with bandwidth $\b M$ we have
		\begin{align}
			\label{eq:error_est_bandlim}
			\big\| \b{\hat f} - \b A^* \b W \b f \big\|_p
			&\leq
			|\I_{\b{M}}| \, \varepsilon  \cdot \big\| \b{\hat f} \big\|_p
			+
			(N\,|\I_{\b{M}}|)^{1/p} \, \|\b w\|_{q} \cdot \|Q\|_{C(\T^d)} ,
		\end{align}
		where \mbox{$\b{\hat f} \coloneqq (\hat f(\b k))_{\b k\in\I_{\b M}}$} are the integer evaluations of \eqref{eq:inverse_integral} and 
		\mbox{$Q$}
		in \eqref{eq:quadrature_error} is the pointwise quadrature error of the equispaced quadrature rule \eqref{eq:quadrature_forward}.
	\end{enumerate}
\end{Theorem}

\begin{proof}
	We start with some general considerations that are independent of the function~$f$.
	By \eqref{eq:residual_lgs} we can find \new{\mbox{$\varepsilon \coloneqq \max_{\b k \in \I_{\b M}} |\varepsilon_{\b k}| \geq 0$}, such that \mbox{$|\varepsilon_{\b k}| \leq \varepsilon$}, \mbox{$\b k\in\I_{\b{2M}}$}, and thereby}
	\begin{align*}
		\left|\, \sum_{j=1}^N w_j \,\e^{2\pi\i\b k\b x_j} - \delta_{\b 0,\b k} \,\right|
		\leq \varepsilon,
		\quad \b k\in\I_{\b{2M}} .
	\end{align*}
	Then for all \mbox{$\b k, \b\ell\in\I_{\b{M}}$} with \mbox{$(\b\ell-\b k) \in\I_{\b{2M}}$} this yields
	\begin{align*}
		\left| \left[\b E_\mathrm{r}\right]_{\b k,\b\ell} \,\right|
		=
		\left|\, \sum_{j=1}^N w_j \,\e^{2\pi\i(\b\ell-\b k)\b x_j} - \delta_{\b k,\b\ell} \,\right|
		\leq \varepsilon,
	\end{align*}
	where \mbox{$\b E_\mathrm{r} \coloneqq \b A^* \b W \b A - \b I_{|\I_{\b M}|}$}.
	Hence, we have
	\begin{align}
		\label{eq:error_matrix_norm_1}
		\big\| \b A^* \b W \b A - \b I_{|\I_{\b M}|} \big\|_1 
		=
		\max_{\b\ell\in\I_{\b M}} \,\sum_{\b k \in \I_{\b M}} \left| \left[\b E_\mathrm{r}\right]_{\b k,\b\ell} \,\right|
		\leq
		\max_{\b\ell\in\I_{\b M}} \,\sum_{\b k \in \I_{\b M}} \varepsilon
		=
		|\I_{\b{M}}| \, \varepsilon ,
	\end{align}
	\begin{align}
		\label{eq:error_matrix_norm_infty}
		\big\| \b A^* \b W \b A - \b I_{|\I_{\b M}|} \big\|_\infty 
		=
		\max_{\b k \in \I_{\b M}} \,\sum_{\b\ell\in\I_{\b M}} \left| \left[\b E_\mathrm{r}\right]_{\b k,\b\ell} \,\right|
		\leq
		\max_{\b k \in \I_{\b M}} \,\sum_{\b\ell\in\I_{\b M}} \varepsilon
		=
		|\I_{\b{M}}| \, \varepsilon ,
	\end{align}
	and
	\begin{align}
		\label{eq:error_matrix_norm_fro}
		\big\| \b A^* \b W \b A - \b I_{|\I_{\b M}|} \big\|_{\mathrm F} 
		=
		\sqrt{\sum_{\b k \in \I_{\b M}} \,\sum_{\b\ell\in\I_{\b M}} \left| \left[\b E_\mathrm{r}\right]_{\b k,\b\ell} \,\right|^2}
		\leq
		\sqrt{\sum_{\b k \in \I_{\b M}} \,\sum_{\b\ell\in\I_{\b M}} \varepsilon^2} 
		=
		|\I_{\b{M}}| \, \varepsilon .
	\end{align}
	Considering the approximation error of \eqref{eq:Fourier_coeffs_nonequi}, it can be estimated by
	\begin{align}
		\big\| \b{\hat f} - \b A^* \b W \b f \big\|_p
		&\leq
		\big\| \b{\hat f} - \b A^* \b W \b A \b{\hat f} \big\|_p
		+
		\big\| \b A^* \b W \b A \b{\hat f} - \b A^* \b W \b f \big\|_p 
		\notag \\
		&=
		\big\| \big( \b A^* \b W \b A - \b I_{|\I_{\b M}|} \big) \b{\hat f} \big\|_p
		+
		\big\| \b A^* \b W \big( \b A \b{\hat f} - \b f \big) \big\|_p \notag \\
		&\leq
		\big\| \b A^* \b W \b A - \b I_{|\I_{\b M}|} \big\|_p  \cdot \big\| \b{\hat f} \big\|_p 
		+
		\big\| \b A^* \b W \big\|_p \cdot \big\| \b A \b{\hat f} - \b f \big\|_p .
		\label{eq:approx_error_general_second}
	\end{align}
	Using \mbox{$\b A \in \C^{N\times|\I_{\b{M}}|}$} from \eqref{eq:matrix_A} as well as \mbox{$\b W = \mathrm{diag} (w_j)_{j=1}^N = \diag(\b w)$} we have
	\begin{align*}
		\big\| \b A^* \b W \big\|_1 
		=
		\max_{j=1,\dots,N} \,\sum_{\b k \in \I_{\b M}} |w_j| \cdot \big|\e^{-2\pi\i\b k \b x_j}\big|
		\leq
		\max_{j=1,\dots,N} |w_j| \cdot \sum_{\b k \in \I_{\b M}} 1
		=
		|\I_{\b{M}}| \cdot \|\b w\|_{\infty} ,
	\end{align*}
	\begin{align*}
		\big\| \b A^* \b W \big\|_\infty 
		=
		\max_{\b k \in \I_{\b M}} \,\sum_{j=1}^N |w_j| \cdot \big|\e^{-2\pi\i\b k \b x_j}\big|
		\leq
		\sum_{j=1}^N |w_j| \cdot \max_{\b k \in \I_{\b M}} 1
		=
		\|\b w\|_{1} ,
	\end{align*}
	and
	\begin{align*}
		\big\| \b A^* \b W \big\|_{\mathrm F} 
		=
		\sqrt{\sum_{\b k \in \I_{\b M}} \,\sum_{j=1}^N |w_j|^2 \cdot \big|\e^{-2\pi\i\b k \b x_j}\big|^2}
		\leq
		\sqrt{\sum_{j=1}^N |w_j|^2 \cdot |\I_{\b{M}}|} 
		=
		\sqrt{|\I_{\b{M}}|} \cdot \|\b w\|_{2} .
	\end{align*}
	Hence, from \eqref{eq:error_matrix_norm_1} -- \eqref{eq:approx_error_general_second} and \mbox{$\|\cdot\|_2 \leq \|\cdot\|_{\mathrm F}$} it follows that
	\begin{align}
		\label{eq:approx_error_general}
		\big\| \b{\hat f} - \b A^* \b W \b f \big\|_p
		&\leq
		|\I_{\b{M}}| \, \varepsilon  \cdot \big\| \b{\hat f} \big\|_p
		+
		\big\| \b A \b{\hat f} - \b f \big\|_p \cdot |\I_{\b{M}}|^{1/p} \, \|\b w\|_{q}
	\end{align}
	for \mbox{$p \in \{1,2,\infty\}$} with \mbox{$\frac 1p + \frac 1q = 1$}.
	Now it merely remains to estimate 
	\mbox{$\big\| \b A \b{\hat f} - \b f \big\|_p$}
	for the specific choice of $f$.
	
	\ref{error_est_trig_poly}: 
	Since a trigonometric polynomial \eqref{eq:trig_poly_2d} of degree~$\b M$ satisfies \mbox{$\b A \b{\hat f} = \b f$}, the second error term in~\eqref{eq:approx_error_general} vanishes and we obtain the assertion~\eqref{eq:error_est_trig_poly}.
	
	\ref{error_est_periodic}: 
	When considering a general 1-periodic function \mbox{$f\in L_2(\T^d) \cap C(\T^d)$} in \eqref{eq:Fourier_series} we have
	\begin{align*}
		\Big| \big[ \b A \b{\hat f} - \b f \big]_j \Big|
		&=
		\Bigg| f(\b x_j) - \sum_{\b k\in\I_{\b{M}}} c_{\b k}(f) \,\e^{2\pi\i\b k\b x_j} \Bigg| \\
		&\leq 
		\max_{\b x\in\T^d} \Bigg| \sum_{\b k\in\Z^d\setminus \I_{\b{M}}} c_{\b k}(f) \,\e^{2\pi\i\b k\b x} \Bigg|
		=
		\| f - p_{\b{M}} \|_{C(\T^d)}
		\quad j=1,\dots, N,
	\end{align*}
	with the best approximating trigonometric polynomial $p_{\b{M}}$ of degree $\b M$ of $f$. 
	Thus, this yields \mbox{$\big\| \b A \b{\hat f} - \b f \big\|_p \leq N^{1/p} \, \| f - p_{\b{M}} \|_{C(\T^d)}$}
	and by \eqref{eq:approx_error_general} the assertion \eqref{eq:error_est_periodic}.
	
	\ref{error_est_bandlim}: 
	For a bandlimited function \mbox{$f\in L_1(\R^d)\cap C_0(\R^d)$} with bandwidth~$\b M$ 
	we may use the notation \mbox{$\b{\hat f} \coloneqq (\hat f(\b k))_{\b k\in\I_{\b M}}$} as well as the inverse Fourier transform~\eqref{eq:forward_integral} to estimate
	\begin{align*}
		\Big| \big[ \b A \b{\hat f} - \b f \big]_j \Big|
		&=
		\Bigg| f(\b x_j) - \sum_{\b k\in\I_{\b{M}}} \hat{f}(\b k) \,\e^{2\pi\i\b k\b x_j} \Bigg| 
		\leq
		\max_{\b x \in \T^d} | Q(\b x) |
		=
		\| Q \|_{C(\T^d)},
		\quad j=1,\dots,N ,
	\end{align*}
	with the pointwise quadrature error
	\begin{align}
		\label{eq:quadrature_error}
		Q(\b x) \coloneqq \int_{\left[-\frac M2,\frac M2\right)^d} \hat{f}(\b v) \,\e^{2\pi\i\b v \b x} \,\mathrm d\b v - \sum_{\b k\in \I_{\b{M}}} \hat{f}(\b k) \,\e^{2\pi\i\b k\b x}
	\end{align}
	of the uniform quadrature rule \eqref{eq:quadrature_forward}.
	For detailed investigations of quadrature errors for band\-limited functions we refer to \cite{KiPoTa21, GoRo22}.
	Hence, we obtain 
	\mbox{$\big\| \b A \b{\hat f} - \b f \big\|_p \leq N^{1/p} \, \| Q \|_{C(\T^d)}$}
	and by \eqref{eq:approx_error_general} the assertion \eqref{eq:error_est_bandlim}.

\end{proof}

By Corollary~\ref{Corollary:cond_exact} it is known that in the setting of trigonometric polynomials there is a linkage between an exact reconstruction \eqref{eq:exact_reconstr} and the matrix product
\mbox{$\b A^* \b W \b A$} being equal to identity \mbox{$\b I_{|\I_{\b{M}}|}$}.
The following theorem shows that the error of the reconstruction~\eqref{eq:Fourier_coeffs_nonequi} also affects the condition of the matrix \mbox{$\b A^* \b W \b A$}.

\begin{Theorem}
	\label{Thm:est_cond}
	Let \mbox{$\b A \in \C^{N\times|\I_{\b{M}}|}$} from \eqref{eq:matrix_A}, \mbox{$\b W = \mathrm{diag} (w_j)_{j=1}^N$} and \mbox{$\varepsilon \geq 0$} be given as in Theorem~\ref{Thm:error_est_lgs}.
	If additionally \mbox{$\varepsilon \,|\I_{\b{M}}| < 1$} is fulfilled, then we have
	\begin{align}
		\label{eq:est_cond}
		1 \le \kappa_2(\b A^* \b W \b A) \leq \frac{1+\varepsilon \,|\I_{\b{M}}|}{1-\varepsilon \,|\I_{\b{M}}|}
	\end{align}
	for the condition number \mbox{$\kappa_2(\b X) \coloneqq \|\b X\|_2 \|\b X^{-1}\|_2$}.
\end{Theorem}

\begin{proof}
	To estimate the condition number \mbox{$\kappa_2(\b A^* \b W \b A)$} we need to determine the norms \mbox{$\big\|\b A^* \b W \b A\big\|_2$} and \mbox{$\big\|(\b A^* \b W \b A)^{-1}\big\|_2$}.
	By \eqref{eq:residual_lgs} it is known that \mbox{$\b A^* \b W \b A = \b I_{|\I_{\b{M}}|} + \b{\mathcal E}$}, where 
	\mbox{$\b{\mathcal E} \coloneqq \left( \varepsilon_{\b\ell-\b k} \right)_{\b\ell, \b k\in\I_{\b{M}}}$}, and
	therefore we have
	\begin{align}   
		\label{eq:cond_norm_matrix}
		\big\|\b A^* \b W \b A\big\|_2 
		= \big\|\b I_{|\I_{\b{M}}|} + \b{\mathcal E}\big\|_2 
		\leq \big\|\b I_{|\I_{\b{M}}|}\big\|_2 + \|\b{\mathcal E}\|_2 .
	\end{align}
	Moreover, it is known by the theory of Neumann series, cf.~\cite[Thm.~4.20]{St98}, that if \mbox{$\big\|\b I_{|\I_{\b{M}}|} - \b T\big\|_2 < 1$} holds for a matrix \mbox{$\b T \in \C^{|\I_{\b{M}}|\times|\I_{\b{M}}|}$}, then \mbox{$\b T$} is invertible and its inverse is given by
	\begin{align*}
		\b T^{-1} = \sum_{n=0}^{\infty} \left( \b I_{|\I_{\b{M}}|} - \b T \right)^n. 
	\end{align*}
	Using this property for \mbox{$\b T = \b A^* \b W \b A$} we have
	\begin{align}
		\label{eq:cond_norm_inverse}
		\big\| ( \b A^* \b W \b A )^{-1} \big\|_2
		= \left\| \sum_{n=0}^{\infty} \left( \b I_{|\I_{\b{M}}|} - \b A^* \b W \b A \right)^n \right\|_2
		= \left\| \sum_{n=0}^{\infty} \b{\mathcal E}^n \right\|_2
		\leq \sum_{n=0}^{\infty} \left\| \b{\mathcal E}^n \right\|_2 ,
	\end{align}
	in case that \mbox{$\big\|\b I_{|\I_{\b{M}}|} - \b A^* \b W \b A\big\|_2 = \|\b{\mathcal E}\|_2 < 1$}.
	Hence, by \eqref{eq:cond_norm_matrix} and \eqref{eq:cond_norm_inverse} we obtain
	\begin{align}
		\label{eq:cond_first_est}
		\kappa_2(\b A^* \b W \b A)
		\leq 
		\left( 1 + \|\b{\mathcal E}\|_2 \right) \cdot
		\left( \sum_{n=0}^{\infty} \left\| \b{\mathcal E}^n \right\|_2 \right) .
	\end{align}
	Additionally, we know that \mbox{$|\varepsilon_{\b k}| \leq \varepsilon$}, \mbox{$\b k\in\I_{\b{2M}}$}, with some \mbox{$\varepsilon > 0$} and therefore
	\begin{align}
		\label{eq:norm_eps}
		\|\b{\mathcal E}\|_2 \leq \|\b{\mathcal E}\|_{\mathrm F}
		=
		\sqrt{\sum_{\b k \in \I_{\b M}} \,\sum_{\b\ell \in \I_{\b M}} |\varepsilon_{\b\ell-\b k}|^2}
		\leq 
		\sqrt{\sum_{\b k \in \I_{\b M}} \,\sum_{\b\ell \in \I_{\b M}} \varepsilon ^2}
		=
		\varepsilon \,|\I_{\b{M}}| .
	\end{align}
	In other words, the correctness of \eqref{eq:cond_norm_inverse} is ensured if \mbox{$\varepsilon \,|\I_{\b{M}}| < 1$}.
	Since the spectral norm is a sub-multiplicative norm, \eqref{eq:norm_eps} also implies
	\mbox{$\|\b{\mathcal E}^n\|_2 \leq \|\b{\mathcal E}\|_2^n \leq \left( \varepsilon \,|\I_{\b{M}}| \right)^n$}.
	Consequently, we have
	\begin{align}
		\label{eq:norm_eps_n}
		\sum_{n=0}^{\infty} \left\| \b{\mathcal E}^n \right\|_2
		\leq
		\sum_{n=0}^{\infty} \left( \varepsilon \,|\I_{\b{M}}| \right)^n
		=
		\frac{1}{1-\varepsilon \,|\I_{\b{M}}|} .
	\end{align}
	Thus, combining \eqref{eq:cond_first_est}, \eqref{eq:norm_eps} and \eqref{eq:norm_eps_n} yields the assertion \eqref{eq:est_cond}.
\end{proof}

\subsection{Connection to certain density compensation approaches from literature \label{sec:dcf_literature}}

In literature a variety of density compensation approaches can be found that are concerned with the setting of bandlimited functions and make use of a $\sinc$ transform
\begin{align}
	\label{eq:matrix_C}
	\b C
	\coloneqq 
	\bigg( |\I_{\b M}| \, \sinc\left(M\pi\left(\b x_j - \b M^{-1} \odot\b\ell\right)\right) \bigg)_{j=1,\, \b\ell \in \I_{\b M}}^{N} 
	\ \in \mathbb R^{N\times |\I_{\b M}|} 
\end{align}
instead of the Fourier transform~\eqref{eq:Fourier_coeffs_nonequi}. 
Namely, instead of directly using the quadrature~\eqref{eq:quadrature_inverse} for reconstruction, in these methods it is inserted into the inverse Fourier transform~\eqref{eq:forward_integral}, i.\,e.,
\begin{align}
	\label{eq:interpolation_sinc}
	f(\b x) 
	&= 
	\int\limits_{\left[-\frac M2,\frac M2\right)^d} \hat f(\b v)\,\e^{2\pi\i \b v\b x}\,\mathrm d\b v
	\approx
	\sum_{j=1}^{N} w_j\, f(\b x_j)
	\int\limits_{\left[-\frac M2,\frac M2\right)^d} \e^{-2\pi\i \b v (\b x_j-\b x)} \,\mathrm d\b v \notag \\
	&=
	\sum_{j=1}^{N} w_j\, f(\b x_j) \cdot |\I_{\b M}| \, \sinc(M\pi(\b x_j-\b x)),
	\quad \b x \in \R^d.
\end{align}
By using the $\sinc$ matrix \mbox{$\b C \in \mathbb R^{N\times |\I_{\b M}|}$} from \eqref{eq:matrix_C}, 
the weight matrix
\mbox{$\b W = \mathrm{diag} (w_j)_{j=1}^N$}
as well as the vectors
\mbox{$\b f = \left(f(\b x_j)\right)_{j=1}^N$} and
\mbox{$\b{\tilde f} = (f(\b M^{-1} \odot\b\ell))_{\b l \in \I_{\b M}}$},
the evaluation of \eqref{eq:interpolation_sinc} at equispaced points \mbox{$\b M^{-1} \odot\b\ell$}, \mbox{$\b\ell \in \I_{\b M}$},
can be denoted as \mbox{$\b{\tilde f} \approx \b C^* \b W \b f$}.
Using the equispaced quadrature rule in \eqref{eq:Fourier_coeffs_equi}, we find that evaluations \mbox{$\hat f(\b k)$} of \eqref{eq:inverse_integral} at the uniform grid points \mbox{$\b k \in \I_{\b{M}}$} can be approximated by \eqref{eq:approx_inverse_integral} by means of a simple FFT.
In matrix-vector notation this can be written as \mbox{$\b{\hat f} \approx \b{\tilde D}^* \b F_{|\I_{\b{M}}|}^* \b{\tilde f}$} where 
\mbox{$\b{\hat f} = (\hat f(\b k))_{\b k\in\I_{\b M}}$}, 
\mbox{$\b F_{|\I_{\b{M}}|} \coloneqq ( \e^{2\pi\i \b k (\b M^{-1}\odot\,\b\ell)} )_{\b\ell,\, \b k \in \I_{\b M}}$}, 
cf.~\eqref{eq:matrix_F}, 
and \mbox{$\b{\tilde D} \coloneqq \frac{1}{|\I_{\b M}|}\, \b I_{|\I_{\b M}|}$}. 
Thus, all in all one obtains an approximation of the form 
\mbox{$\b{\hat f} \approx \b{\tilde D}^* \b F_{|\I_{\b{M}}|}^* \b C^* \b W \b f$}.

Here some of these approaches, cf.~\cite{EgKiPo22}, shall be reconsidered in the context of the Fourier transform~\eqref{eq:Fourier_coeffs_nonequi}. 
We especially focus on the connection of the approaches among each other as well as to our new method introduced in Section~\ref{sec:trig_poly}.


\subsubsection{Density compensation using the pseudoinverse \label{sec:dcf}}
Since \eqref{eq:problem_infft} is in general not exactly solvable, we study the corresponding least squares problem, instead, i.\,e., we look for the approximant that minimizes the residual norm
\mbox{$\big\|\b f - \b A \b{\hat f}\big\|_2$}.
It is known (e.\,g. \cite[p.~15]{Bj96}) that this problem always has the unique solution
\begin{equation}
	\label{eq:solution_pseudoinversion_nfft}
	\b{\hat f} 
	\approx 
	\b{\tilde h}^{\mathrm{pinv}}
	\coloneqq
	\b A^\dagger \b f
\end{equation}
with the Moore-Penrose pseudoinverse $\b A^\dagger$.
Comparing \eqref{eq:solution_pseudoinversion_nfft} to the density compensation approach \eqref{eq:Fourier_coeffs_nonequi}, the weights $w_j$ should be chosen such that the matrix product \mbox{$\b A^* \b W$} approximates the pseudoinverse \mbox{$\b A^\dagger$} as best as possible,
i.\,e., we study the optimization problem
\begin{equation}
	\label{eq:Frobenius_norm_dcf}
	\underset{\b W = \mathrm{diag}(w_j)_{j=1}^N}{\text{Minimize }} \ 
	\big\| \b A^* \b W - \b A^\dagger \big\|_{\mathrm F}^2,
\end{equation}
where \mbox{$\|\cdot\|_{\mathrm F}$} denotes the Frobenius norm of a matrix.
It was shown in \cite{Sedarat00} that the solution to this least squares problem can be computed as
\begin{align}
	\label{eq:sol_dcf}
	w_{j}
	&=
	\frac{[\b A \b A^\dagger]_{j,j}}
	{[\b A \b A^*]_{j,j}}
	=
	\frac{1}{|\I_{\b M}|} \cdot [\b A \b A^\dagger]_{j,j},
	\quad j=1,\dots N.
\end{align}
However, since a singular value decomposition is necessary for the calculations in \eqref{eq:sol_dcf}, we obtain a high complexity of \mbox{$\mathcal O(N^2\,|\I_{\b M}|+|\I_{\b M}|^3)$}.
Therefore, we study some more sophisticated least squares approaches in the following.

\subsubsection{Density compensation using weighted normal equations of first kind \label{sec:wcf}}

It is known, that every least squares solution to \eqref{eq:problem_infft} satisfies the weighted normal equations of first kind 
\mbox{$\b A^* \b W \b A \b{\hat f} = \b A^* \b W \b f$},
see e.\,g. \cite[Thm.~1.1.2]{Bj96}.
As already mentioned in Corollary~\ref{Corollary:cond_exact}, we have an exact reconstruction formula \eqref{eq:Fourier_coeffs_nonequi} for all trigonometric polynomials \eqref{eq:trig_poly_2d} of degree $\b M$, if \mbox{$\b A^* \b W \b A = \b I_{|\I_{\b M}|}$} is fulfilled.
Thus, we aim to compute optional weights $w_j$, \mbox{$j=1,\dots,N$}, by considering the optimization problem
\begin{equation}
	\label{eq:Frobenius_norm_wcf}
	\underset{\b W = \mathrm{diag}(w_j)_{j=1}^N}{\text{Minimize }} \ 
	\big\| \b A^* \b W \b A - \b I_{|\I_{\b M}|} \big\|_{\mathrm F}^2.
\end{equation}
Analogous to \cite{Sedarat00} this could also be derived from \eqref{eq:Frobenius_norm_dcf} by introducing a right-hand scaling in the domain of measured data and minimizing the Frobenius norm of the weighted error matrix
\mbox{$\b E_\mathrm{r} \coloneqq \b E \cdot \b A$}, where
\mbox{$\b E \coloneqq \b A^* \b W - \b A^\dagger$} 
is the error matrix in \eqref{eq:Frobenius_norm_dcf}.

In \cite{Rosenfeld98} it was shown that a solution \mbox{$\b W = \diag(\b w)$} to \eqref{eq:Frobenius_norm_wcf} can be obtained by solving \mbox{$\b S \b w = \b b$} with
\begin{align}
	\label{eq:wcf_system_matrix}
	\b S 
	&\coloneqq
	\Big( \left| \left[ \b A \b A^* \right]_{j,h} \right|^2 \Big)_{j,h=1}^N
	\quad\text{ and }\quad
	\b b 
	= 
	|\I_{\b M}| \cdot \b 1_N .
\end{align}
However, since \mbox{$\b S \b w = \b b$} is not separable for single $w_j$, \mbox{$j=1,\dots,N$},
computing these weights is of complexity \mbox{$\mathcal O(N^3)$}.
This is why the authors in \cite{Sedarat00} restricted themselves to a maximal image size of \mbox{$64\times 64$} pixels, which corresponds to setting \mbox{$M=64$}.


\subsubsection{Density compensation using weighted normal equations of second kind \label{sec:pcf}}

Another approach for density compensation factors is based on
the weighted normal equations of second kind
\begin{align}
	\label{eq:normal_equations_second_kind}
	\b A \b A^* \b W \b y = \b f,\quad \b A^* \b W \b y = \b{\hat f}.
\end{align}
We recognize that by \eqref{eq:normal_equations_second_kind} we are given an exact approximation \mbox{$\b{\hat f} = \b A^* \b W \b f$} of the Fourier coefficients in \eqref{eq:Fourier_coeffs_nonequi}
in case \mbox{$\b y = \b f$}, and thereby \mbox{$\b A \b A^* \b W = \b I_{N}$}.
To this end, we consider the optimization problem
\begin{equation}
	\label{eq:Frobenius_norm_pcf}
	\underset{\b W = \mathrm{diag}(w_j)_{j=1}^N}{\text{Minimize }} \ 
	\| \b A \b A^* \b W - \b I_{N} \|_{\mathrm F}^2.
\end{equation}
As in Section~\ref{sec:wcf}, we remark that this optimization problem \eqref{eq:Frobenius_norm_pcf} could also be derived from \eqref{eq:Frobenius_norm_dcf} by introducing an additional left-hand scaling in the Fourier domain and minimizing the Frobenius norm of the weighted error matrix
\mbox{$\b E_\mathrm{l} \coloneqq \b A \cdot \b E$}.

\begin{Remark}
	An analogous approach considering the $\sinc$ transform~\eqref{eq:matrix_C} instead of the Fourier transform~\eqref{eq:Fourier_coeffs_nonequi} was already studied in \cite{PiMe99}. 
	Another version using a $\sinc$ transform evaluated at pointwise differences of the nonequispaced points instead of \eqref{eq:matrix_C} was studied in \cite{ChMu98, GrLeIn06}, where it was claimed that this approach coincides with the one in \cite{PiMe99}.
	However, we remark that due to the sampling theorem of Shannon-Whittaker-Kotelnikov, see~\eqref{eq:shannon}, applied to the function \mbox{$f(\b x)=\sinc(M\pi(\b x_j-\b x))$}, i.\,e.,
	\begin{align*}
		\sinc(M\pi(\b x_j-\b x)) 
		= 
		\sum_{\b\ell\in\Z^d} \sinc(M\pi(\b x_j-\b M^{-1} \odot {\b \ell})) 
		\,\sinc\left(M\pi\left(\b x - \b M^{-1} \odot {\b \ell}\right)\right) 
	\end{align*}
	and its evaluation at \mbox{$\b x = \b x_h$}, \mbox{$h=1,\dots,N$}, this claim only holds asymptotically for \mbox{$|\I_{\b M}|\to\infty$} in the setting of the $\sinc$ transform. 
	
	In contrast, when using the Fourier transform~\eqref{eq:Fourier_coeffs_nonequi} this equality can directly be seen.
	Then the analog to \cite{GrLeIn06} utilizes an approximation of the form $\b f \approx \b H \b W \b f$, where the matrix $\b H$ is defined as the system matrix of \eqref{eq:problem_infft} evaluated at pointwise differences of the nonequispaced points, i.\,e., 
	\begin{equation}
		\label{eq:matrix_H}	 
		\b H \coloneqq \Bigg( \sum_{\b k\in\I_{\b M}} \e^{2\pi\i \b k (\b x_j-\b x_h)} \Bigg)_{j,h=1}^N.
	\end{equation}
	Since by \eqref{eq:matrix_A} we have $\b H = \b A \b A^*$, minimizing the approximation error
	\begin{align*}
		\left\| \b H^* \b W \b f - \b f \right\|_2^2
		=
		\left\| \b A \b A^* \b W \b f - \b f \right\|_2^2
		&\leq
		\left\| \b A \b A^* \b W - \b I_N \right\|_{\mathrm F}^2 \cdot
		\left\| \b f \right\|_2^2,
	\end{align*}
	leads to the optimization problem \eqref{eq:Frobenius_norm_pcf} as well.
	\ex
\end{Remark}

It was shown in \cite{PiMe99} that the minimizer of \eqref{eq:Frobenius_norm_pcf} is given by
\begin{equation}
	\label{eq:weights_pcf}
	w_j
	=
	\frac{|\I_{\b M}|}
	{\sum_{h=1}^N \big|\! \left[ \b A \b A^* \right]_{j,h} \!\big|^2},
	\quad j=1,\dots,N.
\end{equation}
Since for fixed~$j$ the computation of \mbox{$\left[ \b A \b A^* \right]_{j,h}$}, \mbox{$h=1,\dots,N,$} is of complexity \mbox{$\mathcal O(N \, |\I_{\b M}|)$}, 
the weights \eqref{eq:weights_pcf} can be computed in \mbox{$\mathcal{O}(N^2 \, |\I_{\b M}|)$} arithmetic operations.
However, due to the explicit representation~\eqref{eq:matrix_H} the computation of \mbox{$\left[ \b A \b A^* \right]_{j,h}$}, \mbox{$h=1,\dots,N,$} for fixed~$j$ can be accelerated by means of the NFFT (see Algorithm~\ref{alg:nfft}).
Then this step takes \mbox{$\mathcal O(|\I_{\b M}|\log(|\I_{\b M}|) + N)$} arithmetic operations and the overall complexity is given by \mbox{$\mathcal O(N \cdot |\I_{\b M}|\log(|\I_{\b M}|) + N^2)$}.

As mentioned in \cite{PiMe99} one could also consider a simplified version of the optimization problem~\eqref{eq:Frobenius_norm_pcf}
%
by reducing the number of conditions, e.\,g. by summing the columns on both sides of \mbox{$\b A \b A^* \b W = \b I_{N}$} as
\begin{equation}
	\label{eq:relaxed_problem_long}
	\sum_{j=1}^N w_j \sum_{\b k\in\I_{\b M}} \e^{2\pi\i \b k (\b x_h-\b x_j)}
	=
	\sum_{j=1}^N \delta_{j,h} = 1,
	\quad h=1,\dots,N.
\end{equation}
By means of \eqref{eq:matrix_A} this can be written as
\mbox{$\b A \b A^* \b w = \b 1_N$}.
Since fast multiplication with $\b A$ and $\b A^*$ can be realized using the NFFT (see Algorithm~\ref{alg:nfft}) and the adjoint NFFT (see Algorithm~\ref{alg:nfft*}), respectively, a solution to the linear system of equations \mbox{$\b A \b A^* \b w = \b 1_N$} can be computed iteratively with arithmetic complexity \mbox{$\mathcal{O}(|\I_{\b M}|\log(|\I_{\b M}|)+N)$}.


Finally, we investigate the connection of this approach to our method introduced in Section~\ref{sec:trig_poly}.
To this end, suppose the linear system \eqref{eq:claim_quadrature_coefficients} is fulfilled for given \mbox{$\b w\in\C^N$}, i.\,e., by \mbox{$\b A^* = \overline{\b A^T}$} we have
\mbox{$\left( \delta_{\b 0,\b k} \right)_{\b k\in{\I_{\b M}}} = \overline{\b A^T \b w} = \b A^* \overline{\b w}$}.
Then multiplication with \mbox{$\b A \in \C^{N \times |\I_{\b{M}}|}$} in \eqref{eq:matrix_A} yields
\begin{equation*}
	\b A \b A^* \overline{\b w}
	=
	\b A \cdot \left( \delta_{\b 0,\b k} \right)_{\b k\in{\I_{\b M}}}
	=
	\left( \sum_{\b k\in\I_{\b M}} \delta_{\b 0,\b k} \cdot \e^{2\pi\i \b k \b x_j} \right)_{j=1}^N
	=
	\b 1_N.
\end{equation*}
In other words, an exact solution $\b w$ to the linear system \eqref{eq:claim_quadrature_coefficients}
implies that the conjugate complex weights $\b{\overline w}$ exactly solve the system \eqref{eq:relaxed_problem_long}.
However, the reversal does not hold true and therefore \eqref{eq:relaxed_problem_long} is not equivalent to \eqref{eq:claim_quadrature_coefficients}.
Moreover, we have seen in Corollary~\ref{Corollary:cond_exact} that an augmented variant of \eqref{eq:claim_quadrature_coefficients}, namely \eqref{eq:claim_quadrature_coefficients_double}, is necessary to obtain an exact reconstruction \mbox{$\hat f_{\b k} = h_{\b k}^{\mathrm{w}}$} in \eqref{eq:Fourier_coeffs_nonequi} for trigonometric polynomials \eqref{eq:trig_poly_2d} with maximum degree $\b M$.

\section{Direct inversion using matrix optimization \label{sec:opt_B}}

As seen in Remark~\ref{Remark:opt_density}, the previously considered density compensation techniques can be regarded as an optimization of the sparse matrix \mbox{$\b B\in\R^{N\times|\I_{\b M_{\b\sigma}}|}$} from the NFFT, cf. Section~\ref{subsec:nfft}.
Since density compensation allows only $N$ degrees of freedom, this limitation shall now be softened, i.\,e., 
instead of searching for optimal scaling factors for the rows of~$\b B$, we now study the optimization of each nonzero entry of the sparse matrix~$\b B$, cf.~\cite{KiPo20}.
To this end, we firstly have another look at the equispaced setting.
It is known by \eqref{eq:matrix_product_F*F} and \eqref{eq:matrix_product_FF*}, 
that for equispaced points and appropriately chosen parameters a one-sided inversion is given by composition of the Fourier matrix and its adjoint.
Hence, we aim to use this result to find a good approximation of the inverse in the general setting.

Considering problem \eqref{eq:problem_infft} we seek to find an appropriate matrix $\b X$ such that we have \mbox{$\b X \b A \approx \b I_{|\I_{\b M}|}$}, since then we can simply compute an approximation of the Fourier coefficients by means of
\mbox{$\b X \b f = \b X \b A \b{\hat f} \approx \b{\hat f}$}.
To find this left-inverse $\b X$, we utilize the fact that in the equispaced case it is known that \eqref{eq:matrix_product_F*F} holds in the overdetermined setting \mbox{$|\I_{\b M}| \leq N$}.
In \mbox{addition}, we also incorporate the approximate factorization \mbox{$\b A^* \approx \b D^* \b F^* \b B^*$} of the adjoint NFFT, cf.~Section~\ref{subsec:nfft*},
with the matrices $\b D\in \mathbb C ^{|\I_{\b M}|\times |\I_{\b M}|}$, $\b F\in \mathbb C ^{|\I_{\b M_{\b\sigma}}|\times |\I_{\b M}|}$ and $\b B\in \mathbb R ^{N\times |\I_{\b M_{\b\sigma}}|}$ defined in \eqref{eq:matrix_D}, \eqref{eq:matrix_F} and \eqref{eq:matrix_B}.
Combining both ingredients we aim for an approximation of the form \mbox{$\b D^* \b F^* \b B^* \b A \approx \b I_{|\I_{\b M}|}$}.
To achieve an approximation like this, we aim to modify the matrix~$\b B$ such that its sparse structure with at most \mbox{$(2m+1)^d$} entries per row and consequently the arithmetic complexity of its evaluation is preserved.
A matrix satisfying this property we call \mbox{\textit{${(2m+1)^d}$-sparse}}.

\begin{Remark}
	We remark that this approach can also be deduced from the density compensation method in Section~\ref{sec:inv_density_comp} as follows. 
	By Corollary~\ref{Corollary:cond_exact} it is known that an exact reconstruction needs to satisfy \mbox{$\b A^* \b W \b A = \b I_{|\I_{\b{M}}|}$}. 
	Since the reconstruction shall be realized efficiently by means of an adjoint NFFT, one rather studies \mbox{$\b D^* \b F^* \b B^* \b W \b A \approx \b I_{|\I_{\b{M}}|}$}. 
	Using the definition \mbox{$\b{\tilde B} \coloneqq \b W^* \b B$} as in Remark~\ref{Remark:opt_density}, we end up with an approximation of the form \mbox{$\b D^* \b F^* \b{\tilde B}^* \b A \approx \b I_{|\I_{\b{M}}|}$}.
	Thus, optimizing each nonzero entry of the sparse matrix~$\b B$ using this approximation is the natural generalization of the density compensation method from Section~\ref{sec:inv_density_comp}.
	\ex
\end{Remark}


Let $\b{\tilde B}$ denote such a modified matrix.
By defining \mbox{$\b{\tilde{h}} := \b D^* \b F^* \b{\tilde B}^* \b f$},
we recognize that the minimization of the approximation error
\begin{align}
	\label{eq:est_approx_error}
	\big\|\b{\tilde{h}} - \b{\hat f}\big\|_2
	&=
	\big\|\b D^* \b F^* \b{\tilde B}^* \b f - \b{\hat f}\big\|_2
	=
	\big\|\b D^* \b F^* \b{\tilde B}^* \b A \b{\hat f} - \b{\hat f}\big\|_2 \notag \\
	&=
	\big\|\big( \b D^* \b F^* \b{\tilde B}^* \b A - \b I_{|\I_{\b M}|} \big) \b{\hat f} \big\|_2
	\leq
	\big\| \b D^* \b F^* \b{\tilde B}^* \b A -\b I_{|\I_{\b M}|} \big\|_{\mathrm F} \,
	\big\|\b{\hat f}\big\|_2
\end{align}
implies the optimization problem
\begin{equation}
	\label{eq:opt_id}
	\underset{\b{\tilde B} \in \R^{N\times{|\I_{\b M_{\b\sigma}}|}} \colon \b{\tilde B}\, (2m+1)^d\text{-sparse }}{\text{Minimize }} \ 
	\big\| \b D^* \b F^* \b{\tilde B}^* \b A -\b I_{|\I_{\b M}|} \big\|_{\mathrm F}^2 .
\end{equation}
%
Note that a similar idea for the forward problem, i.\,e., the evaluation of \eqref{eq:nfft}, was already studied in \cite{st01}.
By the definition of the Frobenius norm we have $\|\b Z\|_{\mathrm F} = \|\b Z^*\|_{\mathrm F}$, such that \eqref{eq:opt_id} is equivalent to its adjoint
\begin{equation}
	\label{eq:opt_id_adj}
	\underset{\b{\tilde B} \in \R^{N\times{|\I_{\b M_{\b\sigma}}|}} \colon \b{\tilde B}\, (2m+1)^d\text{-sparse }}{\text{Minimize }} \ 
	\big\| \b A^* \b{\tilde B} \b F \b D - \b I_{|\I_{\b M}|} \big\|_{\mathrm F}^2.
\end{equation}
Since it is known by \eqref{eq:matrix_F} that \mbox{$\b F^* \b F = |\I_{\b M_{\b\sigma}}|\, \b I_{|\I_{\b M}|}$} and $\b D\in\R^{|\I_{\b{M}}|\times|\I_{\b{M}}|}$ is diagonal by \eqref{eq:matrix_D}, we have
\mbox{$\tfrac{1}{|\I_{\b M_{\b\sigma}}|} \b D^{-1} \b F^* \b F \b D = \b I_{|\I_{\b M}|}$}.
Thus, due to the fact that the Frobenius norm is a submultiplicative norm, we have
\begin{align}
	\label{eq:submult_fro}
	\big\| \b A^* \b{\tilde B} \b F \b D - \b I_{|\I_{\b M}|} \big\|_{\mathrm F}
	&=
	\big\| \big( \b A^* \b{\tilde B} - \tfrac{1}{|\I_{\b M_{\b\sigma}}|} \b D^{-1} \b F^* \big) \b F \b D \big\|_{\mathrm F}
	\notag \\
	&\leq
	\big\| \b A^* \b{\tilde B} - \tfrac{1}{|\I_{\b M_{\b\sigma}}|} \b D^{-1} \b F^* \big\|_{\mathrm F}
	\,\big\| \b F \b D \big\|_{\mathrm F}.
\end{align}
Hence, we consider the optimization problem
\begin{equation}
	\label{eq:opt_overdet}
	\underset{\b{\tilde B} \in \R^{N\times{|\I_{\b M_{\b\sigma}}|}} \colon \b{\tilde B}\, (2m+1)^d\text{-sparse }}{\text{Minimize }} \ 
	\big\| \b A^* \b{\tilde B} - \tfrac{1}{|\I_{\b M_{\b\sigma}}|} \b D^{-1} \b F^* \big\|_{\mathrm F}^2.
\end{equation}
Based on the definition of the Frobenius norm of a matrix~\mbox{$\b Z \in \R^{k\times n}$} and the definition of the Euclidean norm of a vector \mbox{$\b y\in\R^n$},
we obtain for $\b z_j$ being the columns of~\mbox{$\b Z \in \R^{k\times n}$} that 
\begin{equation}
	\label{eq:prop_Frobenius}
	\|\b Z\|_F^2 = \sum_{i=1}^{k} \sum_{j=1}^{n} |z_{ij}|^2 = \sum_{j=1}^{n} \|\b z_j\|_2^2.
\end{equation}
Since we aim to preserve the property that $\b B$ is a $(2m+1)^d$-sparse matrix, we rewrite the norm in \eqref{eq:opt_overdet} by \eqref{eq:prop_Frobenius} in terms of the columns of $\b{\tilde B}$ considering only the nonzero entries of each column.
To this end, analogously to~\eqref{eq:indexset_x} we define the index set
\begin{equation}
	\label{eq:indexset_l}
	\I_{\b M_{\b\sigma},m}(\b\ell) 
	:= 
	\left\{j \in \left\{1, \dots, N \right\}: 
	\exists\, \b z\in \mathbb Z^d\ \text{with} -m \b 1 \leq \b M_{\b\sigma} \odot \b x_j - \b\ell + \b z \leq m \b 1 \right\}
\end{equation}
of the nonzero entries of the $\b\ell$-th column of \mbox{$\b B\in \R^{N\times |\I_{\b M_{\b\sigma}}|}$}.
Thus, we have
\begin{equation}
	\label{eq:norm_columnwise}
	\big\| \b A^* \b{\tilde B} - \tfrac{1}{|\I_{\b M_{\b\sigma}}|} \b D^{-1} \b F^* \big\|_{\mathrm F}^2 
	=
	\sum_{\b\ell \in \I_{\b M_{\b\sigma}}} \big\|\b H_{\b\ell} \b{\tilde b}_{\b\ell} - \tfrac{1}{|\I_{\b M_{\b\sigma}}|} \b D^{-1} \b f_{\b\ell}\big\|_2^2,
\end{equation}
where \mbox{$\b{\tilde b}_{\b\ell} \in \R^{|\I_{\b M_{\b\sigma},m}(\b\ell)|}$} denotes the vectors of the nonzeros of each column of $\b{\tilde B}$,
\begin{equation}
	\label{eq:matrix_Hl}
	\b H_{\b\ell} \coloneqq \left( \e^{-2\pi\i \b k \b x_j} \right)_{\b k \in \I_{\b M},\, j\in \I_{\b M_{\b\sigma},m}(\b\ell)}
	\in \C^{|\I_{\b M}| \times |\I_{\b M_{\b\sigma},m}(\b\ell)|}
\end{equation}
are the corresponding submatrices of \mbox{$\b A^*\in\C^{|\I_{\b{M}}| \times N}$}, cf.~\eqref{eq:matrix_A}, 
and \mbox{$\b f_{\b\ell}\in\C^{|\I_{\b{M}}|}$} are the columns of \mbox{$\b F^* \in \C^{|\I_{\b{M}}| \times |\I_{\b M_{\b\sigma}}|}$}, cf.~\eqref{eq:matrix_F}.
Hence, we receive the equivalent optimization problems
\begin{equation}
	\label{eq:opt_overdet_columnwise}
	\underset{\b{\tilde b}_{\b\ell} \in \R^{|\I_{\b M_{\b\sigma},m}(\b\ell)|}}{\text{Minimize }}\ 
	\big\|\b H_{\b\ell} \b{\tilde b}_{\b\ell} - \tfrac{1}{|\I_{\b M_{\b\sigma}}|} \b D^{-1} \b f_{\b\ell}\big\|_2^2, 
	\quad \b\ell \in \I_{\b M_{\b\sigma}}.
\end{equation}
Thus, if the matrix \eqref{eq:matrix_Hl} has full column rank, the solution of the least squares problem~\eqref{eq:opt_overdet_columnwise} can be computed by means of the pseudoinverse $\b H_{\b\ell}^\dagger$ as 
\begin{equation}
	\label{eq:solution_overdet}
	\b b_{\b\ell}^{\mathrm{opt}}
	\coloneqq 
	\tfrac{1}{|\I_{\b M_{\b\sigma}}|} 
	\left(
	\b H_{\b\ell} ^* 
	\b H_{\b\ell} 
	\right)^{-1}
	\b H_{\b\ell} ^*
	\b D^{-1} \b f_{\b\ell}, 
	\quad \b\ell \in \I_{\b M_{\b\sigma}}.
\end{equation}
Having these vectors~$\b b_{\b\ell}^{\mathrm{opt}}$ we compose the optimized matrix~$\b B_{\mathrm{opt}}$, observing that $\b b_{\b\ell}^{\mathrm{opt}}$ only consist of the nonzero entries of $\b B_{\mathrm{opt}}$.
Then the approximation of the Fourier coefficients is given by
\begin{equation}
	\label{eq:approx}
	\b{\hat f} 
	\approx
	\b h_{\mathrm{opt}}
	\coloneqq
	\b D^* \b F^* \b B_{\mathrm{opt}}^* \b f.
\end{equation}
In other words, this approach yields an inverse NFFT by modifying the adjoint NFFT.

\begin{Remark}
	\label{Remark:matrix_normal_eq_overdet}
	To achieve an efficient algorithm we now have a closer look at the computation scheme \eqref{eq:solution_overdet}.
	We start with the computation of the matrix \mbox{$\b H_{\b\ell}^* \b H_{\b\ell}$}.
	By introducing the $d$-dimensional Dirichlet kernel
	\begin{equation*}
		D_{\b M}(\b x) 
		\coloneqq
		\sum_{k_1=-\frac{M}{2}+1}^{\frac{M}{2}-1} 
		\dots 
		\sum_{k_d=-\frac{M}{2}+1}^{\frac{M}{2}-1} 
		\e^{2\pi\i \b k \b x} 
		=
		\prod_{t=1}^{d} D_{\frac{M}{2}-1}(x_t)
		=
		\prod_{t=1}^{d} \frac{\sin((M-1)\pi x_t)}{\sin(\pi x_t)}, 
	\end{equation*}
	the matrix \mbox{$\b H_{\b\ell}^* \b H_{\b\ell}$} in \eqref{eq:solution_overdet} can explicitly be stated via
	\begin{align}
		\label{eq:matrix_normal_equations_overdet}
		\b H_{\b\ell}^* \b H_{\b\ell}
		&=
		\Bigg[
		\sum_{\b k \in \I_{\b M}}
		\e^{2\pi\i\b k (\b x_h-\b x_j)}
		\Bigg]_{h,j \in \I_{\b M_{\b\sigma},m}(\b\ell)} \notag\\
		&=
		\left[
		\prod_{t=1}^{d} 
		\left(
		D_{\frac{M}{2}-1}({x_h}_t-{x_j}_t)
		+
		\e^{-M\pi\i ({x_h}_t-{x_j}_t)}
		\right)
		\right]_{h,j \in \I_{\b M_{\b\sigma},m}(\b\ell)} ,
	\end{align}
	i.\,e., for given index set \mbox{$\I_{\b M_{\b\sigma},m}(\b\ell)$} the matrix \mbox{$\b H_{\b\ell}^* \b H_{\b\ell}$} can be determined in \mbox{$\mathcal O(|\I_{\b M_{\b\sigma},m}(\b\ell)|)$} operations.
	Considering the right hand sides of \eqref{eq:solution_overdet}, by definitions \eqref{eq:matrix_Hl}, \eqref{eq:matrix_D} and \eqref{eq:matrix_F} we have
	\begin{align}
		\label{eq:right_side_normal_equations_overdet}
		\b v_{\b\ell} 
		\coloneqq \tfrac{1}{|\I_{\b M_{\b\sigma}}|} \b H_{\b\ell}^* \b D^{-1} \b f_{\b\ell}
		&=
		\Bigg(
		\sum_{\b k \in \I_{\b M}}
		\hat w(\b k)\,
		\e^{2\pi\i\b k \left( \b x_j - \b M_{\b\sigma}^{-1}\odot\,\b\ell \right)}
		\Bigg)_{j \in \I_{\b M_{\b\sigma},m}(\b\ell)} , 
		\quad \b\ell \in \I_{\b M_{\b\sigma}}. 
	\end{align}
	Thus, since 
	\mbox{$\tfrac{1}{|\I_{\b M_{\b\sigma}}|} \b D^{-1} = \mathrm{diag}(\hat w(\b k))_{\b k\in \I_{\b M}}$},
	the computation of $\b v_{\b\ell}$ involves neither multiplication with nor division by the (possibly) huge number $|\I_{\b M_{\b\sigma}}|$ and is therefore numerically stable.
	\ex
\end{Remark}

This leads to the following algorithm.

\begin{algorithm}{Optimization of the sparse matrix $\b B$}
	\label{alg:opt_fast_overdet}
	For \mbox{$d,N \in \N$} let \mbox{$\b x_j \in \T^d$}, \mbox{$j=1,\dots,N,$} be given points.
	Further let \mbox{$\b M \coloneqq M \cdot \b 1_{d}$} with \mbox{$M\in 2\N$}, an over\-sampling factor \mbox{$\sigma\geq 1$} with \mbox{$2\N \ni  M_{\sigma} \coloneqq 2 \lceil \lceil \sigma M \rceil / 2 \rceil$} and \mbox{$\b M_{\b\sigma} \coloneqq M_{\sigma} \cdot \b 1_{d}$} as well as a truncation parameter \mbox{$m \ll M_\sigma$} be given.

	\begin{enumerate}	
		\item
		For $\b\ell \in \I_{\b M_{\b\sigma}}$: 
		\begin{itemize}
			\item[] \hspace{-0.6cm} Determine the index set $\I_{\b M_{\b\sigma},m}(\b\ell)$, cf. \eqref{eq:indexset_l}.
			\hfill $\mathcal O(|\I_{\b M_{\b\sigma},m}(\b\ell)|)$
			\item[] \hspace{-0.6cm} Compute the right side $\b v_{\b\ell}$ via \eqref{eq:right_side_normal_equations_overdet}.
			\hfill $\mathcal O(|\I_{\b M_{\b\sigma},m}(\b\ell)|)$
			\item[] \hspace{-0.6cm} Determine $\b H_{\b\ell}^* \b H_{\b\ell}$ via \eqref{eq:matrix_normal_equations_overdet}.
			\hfill $\mathcal O(|\I_{\b M_{\b\sigma},m}(\b\ell)|^2)$
			\item[] \hspace{-0.6cm} Solve \(\left(\b H_{\b\ell}^* \b H_{\b\ell}\right)\b b_{\b\ell}^{\mathrm{opt}}=\b v_{\b\ell}\), i.\,e., compute $\b b_{\b\ell}^{\mathrm{opt}}$, 
			cf. \eqref{eq:solution_overdet}.
			\hfill \mbox{$\mathcal O(|\I_{\b M_{\b\sigma},m}(\b\ell)|^3)$}
		\end{itemize}
		\item Compose $\b B_{\mathrm{opt}}\in \mathbb R^{N\times |\I_{\b M_{\b\sigma}}|}$ columnwise of the $\b b_{\b\ell}^{\mathrm{opt}} \in \R^{|\I_{\b M_{\b\sigma},m}(\b\ell)|}$. 
		\hfill $\mathcal O(|\I_{\b M}|)$
	\end{enumerate}
	\vspace{1ex} 
	\textnormal{\textbf{Output:}} optimized matrix $\b B_{\mathrm{opt}}$ \hfill
	\textnormal{\textbf{Complexity:}} $\mathcal O(|\I_{\b M}| \cdot |\I_{\b M_{\b\sigma},m}(\b\ell)|^3)$
\end{algorithm}

\begin{figure}[!h]
	\centering
	\captionsetup[subfigure]{justification=centering}
	\begin{subfigure}[t]{0.24\textwidth}
		\centering
		\includegraphics[width=0.8\textwidth]{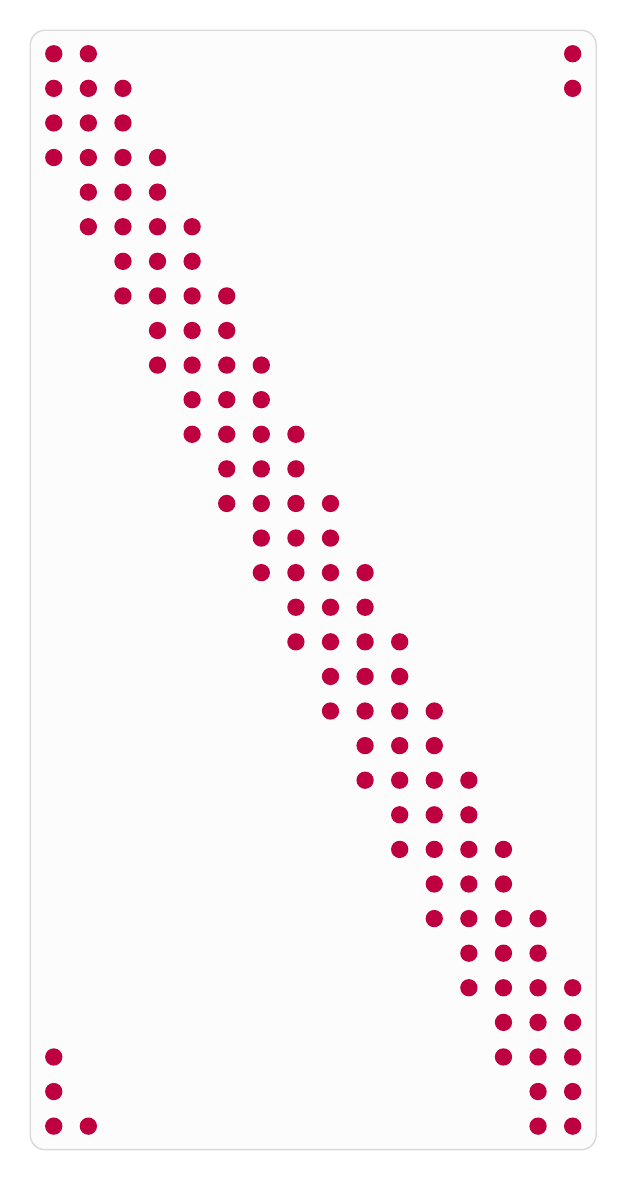}
		\caption{Equispaced points}
	\end{subfigure}
	\begin{subfigure}[t]{0.24\textwidth}
		\centering
		\includegraphics[width=0.8\textwidth]{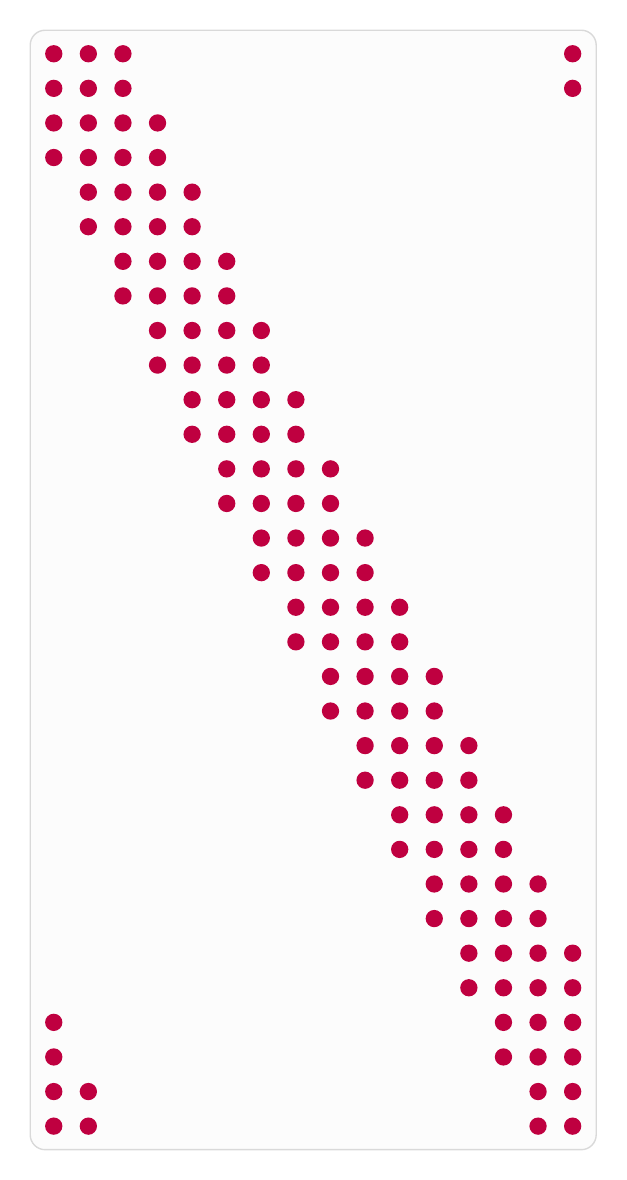}
		\caption{Jittered points}
	\end{subfigure}
	\begin{subfigure}[t]{0.24\textwidth}
		\centering
		\includegraphics[width=0.8\textwidth]{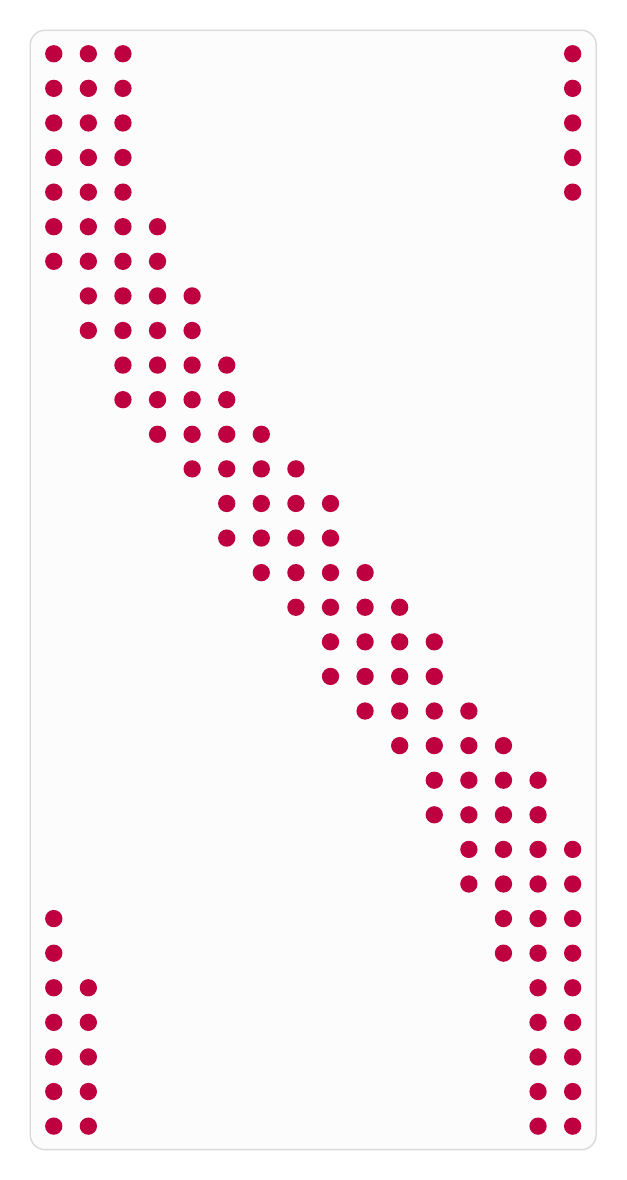}
		\caption{Chebyshev points}
	\end{subfigure}
	\begin{subfigure}[t]{0.24\textwidth}
		\centering
		\includegraphics[width=0.8\textwidth]{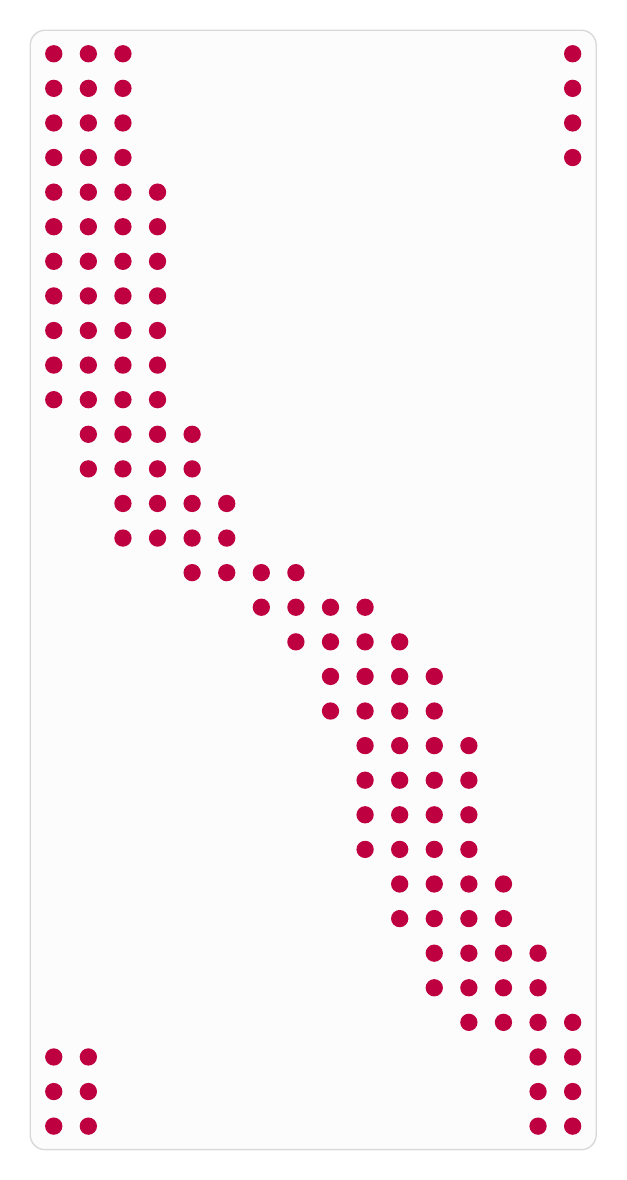}
		\caption{Random points}
	\end{subfigure}
	\caption{
		Nonzero entries of the matrix \mbox{$\b B \in \R^{N\times |\I_{\b M_{\b\sigma}}|}$} for several choices of the points \mbox{$\b x_j\in\T^d$}, \mbox{$j=1,\dots,N,$} with \mbox{$d=1$}, \mbox{$M_\sigma=M=16$}, \mbox{$N=2M$} and \mbox{$m=2$}.
		\label{fig:sparsity_B}}
\end{figure}

Note that a general statement about the dimensions of \mbox{$\b H_{\b\ell} \in \C^{|\I_{\b M}| \times |\I_{\b M_{\b\sigma},m}(\b\ell)|}$} is not possible, since the size of the set~\mbox{$\I_{\b M_{\b\sigma},m}(\b\ell)$} heavily depends on the distribution of the points.
To visualize this circumstance, we depicted some exemplary patterns of the nonzero entries of the original matrix \mbox{$\b B \in \R^{N\times |\I_{\b M_{\b\sigma}}|}$} in Figure~\ref{fig:sparsity_B}.
It can easily be seen that for all choices of the points each row contains the same number of nonzero entries, i.\,e., all index sets \eqref{eq:indexset_x} are of the same size of maximum \mbox{$(2m+1)^d$}.
However, when considering the columns instead, we recognize an evident mismatch in the number of nonzero entries.
We remark that because of the fact that each row of \mbox{$\b B \in \R^{N\times |\I_{\b M_{\b\sigma}}|}$} contains at most \mbox{$(2m+1)^d$} entries, each column contains \mbox{$\frac{N}{|\I_{\b M_{\b\sigma}}|}(2m+1)^d$} entries on average.
A general statement about the maximum size of the index sets~\eqref{eq:indexset_l} cannot be made.
Roughly speaking, the more irregular the distribution of the points is, the larger the index sets~\eqref{eq:indexset_l} can be.
Nevertheless, in general \mbox{$|\I_{\b M_{\b\sigma},m}(\b\ell)|$} is a small constant compared to \mbox{$|\I_{\b{M}}|$}, such that Algorithm~\ref{alg:opt_fast_overdet} ends up with total arithmetic costs of approximately \mbox{$\mathcal O(|\I_{\b M}|)$}.

In conclusion, our approach for an inverse NFFT can be summarized as follows.

\begin{algorithm}{iNFFT -- optimization approach}
	\label{alg:infft_overdet}
	For \mbox{$d,N \in \N$} let \mbox{$\b x_j \in \T^d$}, \mbox{$j=1,\dots,N,$} be given points as well as \mbox{$\b f \in \C^N$}.
	Further let \mbox{$\b M \coloneqq M \cdot \b 1_{d}$} with \mbox{$M\in 2\N$}, an over\-sampling factor \mbox{$\sigma\geq 1$} with \mbox{$2\N \ni  M_{\sigma} \coloneqq 2 \lceil \lceil \sigma M \rceil / 2 \rceil$} and \mbox{$\b M_{\b\sigma} \coloneqq M_{\sigma} \cdot \b 1_{d}$} as well as a truncation parameter \mbox{$m \ll M_\sigma$} be given.
	\begin{enumerate}
		\setcounter{enumi}{-1}
		\item Precompute the optimal sparse matrix $\b B_{\mathrm{opt}}$ using Algorithm~\ref{alg:opt_fast_overdet}.
		\item Compute 
		\mbox{$\b h_{\mathrm{opt}} \coloneqq \b D^* \b F^* \b B_{\mathrm{opt}}^* \b f$}, cf. \eqref{eq:approx},
		by means of a modified adjoint NFFT.
	\end{enumerate}
	\vspace{1ex} 
	\textnormal{\textbf{Output:}} $\b h_{\mathrm{opt}} \approx \b{\hat f}\in\C^{|\I_{\b{M}}|}$, cf. \eqref{eq:problem_infft}. \hfill
	\textnormal{\textbf{Complexity:}} $\mathcal O(|\I_{\b M}|\log(|\I_{\b M}|) + N)$
\end{algorithm}

\begin{Theorem}
	\label{Thm:error_est_opt}
	Let \mbox{$\b B_{\mathrm{opt}}\in\R^{N \times |\I_{\b M_{\b\sigma}}|}$} be the optimized matrix computed by means of Algorithm~\ref{alg:opt_fast_overdet} and let 
	\mbox{$\b h_{\mathrm{opt}} = \b D^* \b F^* \b B_{\mathrm{opt}}^* \b f \in \C^{|\I_{\b{M}}|}$} 
	be the corresponding approximation of~\mbox{$\b{\hat f}$} computed by means of Algorithm~\ref{alg:infft_overdet}.
	Further assume that each column \mbox{$\b b_{\b\ell}^{\mathrm{opt}}\in\R^{|\I_{\b M_{\b\sigma},m}|}$} of \mbox{$\b B_{\mathrm{opt}}\in\R^{N \times |\I_{\b M_{\b\sigma}}|}$} as solution to \eqref{eq:opt_overdet_columnwise} possesses a small residual 
	\begin{align}
		\label{eq:res_opt_columnwise}
		\big\|\b H_{\b\ell} \b b_{\b\ell}^{\mathrm{opt}} - \tfrac{1}{|\I_{\b M_{\b\sigma}}|} \b D^{-1} \b f_{\b\ell}\big\|_2^2 
		= 
		\varepsilon_{\b\ell} \geq 0 ,
		\quad \b\ell\in\I_{\b M_{\b\sigma}} .
	\end{align}
	Then there exists an $\varepsilon \geq 0$ such that
	\begin{align}
		\label{eq:approx_error_opt}
		\big\|\b h_{\mathrm{opt}} - \b{\hat f}\big\|_2^2
		&\leq
		\varepsilon \sum_{\b k \in \I_{\b M}} \frac{1}{\hat w(\b k)^2} \cdot
		\big\|\b{\hat f}\big\|_2^2 .
	\end{align}
	Moreover, the (asymmetric) Dirichlet kernel 
	\begin{align}
		\label{eq:kernel_Dirichlet_asym}
		w_{\mathrm{D}} 
		\coloneqq
		\sum_{\b k\in\I_{\b M}} \e^{2\pi\i \b k \b x} 
		=
		\prod_{t=1}^{d} \left( D_{\frac{M}{2}-1}(x_t) + \e^{-M\pi\i x_t} \right)
	\end{align}
	is the optimal window function for the inverse NFFT in Algorithm~\ref{alg:infft_overdet}.
\end{Theorem}

\begin{proof}
	As in \eqref{eq:est_approx_error} the approximation error can be estimated by
	\begin{align}
		\label{eq:est_approx_error_opt}
		\big\|\b h_{\mathrm{opt}} - \b{\hat f}\big\|_2^2
		&=
		\big\|\big( \b D^* \b F^* \b B_{\mathrm{opt}}^* \b A - \b I_{|\I_{\b M}|} \big) \b{\hat f} \big\|_2^2
		\leq
		\big\| \b D^* \b F^* \b B_{\mathrm{opt}}^* \b A -\b I_{|\I_{\b M}|} \big\|_{\mathrm F}^2 \,
		\big\|\b{\hat f}\big\|_2^2.
	\end{align}
	Using the same arguments as for \eqref{eq:opt_id_adj} and \eqref{eq:submult_fro} we proceed with
	\begin{align}
		\label{eq:submult_fro_opt}
		\big\| \b D^* \b F^* \b B_{\mathrm{opt}}^* \b A -\b I_{|\I_{\b M}|} \big\|_{\mathrm F}^2
		&=
		\big\| \b A^* \b B_{\mathrm{opt}} \b F \b D - \b I_{|\I_{\b M}|} \big\|_{\mathrm F}^2
		\notag \\
		&\leq
		\big\| \b A^* \b B_{\mathrm{opt}} - \tfrac{1}{|\I_{\b M_{\b\sigma}}|} \b D^{-1} \b F^* \big\|_{\mathrm F}^2
		\,\big\| \b F \b D \big\|_{\mathrm F}^2.
	\end{align}

	To estimate the first Frobenius norm in \eqref{eq:submult_fro_opt}, we rewrite it analogously to \eqref{eq:norm_columnwise} columnwise as
	\begin{align*}
		\big\| \b A^* \b B_{\mathrm{opt}} - \tfrac{1}{|\I_{\b M_{\b\sigma}}|} \b D^{-1} \b F^* \big\|_{\mathrm F}^2
		=
		\sum_{\b\ell \in \I_{\b M_{\b\sigma}}} \big\|\b H_{\b\ell} \b b_{\b\ell}^{\mathrm{opt}} - \tfrac{1}{|\I_{\b M_{\b\sigma}}|} \b D^{-1} \b f_{\b\ell}\big\|_2^2 ,
	\end{align*}
	where \mbox{$\b b_{\b\ell}^{\mathrm{opt}} \in \R^{|\I_{\b M_{\b\sigma},m}(\b\ell)|}$} are the nonzeros of the columns of $\b B_{\mathrm{opt}}$, 
	\mbox{$\b H_{\b\ell} \in \C^{|\I_{\b M}| \times |\I_{\b M_{\b\sigma},m}(\b\ell)|}$} in \eqref{eq:matrix_Hl} are the corresponding submatrices of \mbox{$\b A^*\in\C^{|\I_{\b{M}}| \times N}$}, cf.~\eqref{eq:matrix_A}, 
	and \mbox{$\b f_{\b\ell}\in\C^{|\I_{\b{M}}|}$} are the columns of \mbox{$\b F^* \in \C^{|\I_{\b{M}}| \times |\I_{\b M_{\b\sigma}}|}$}, cf.~\eqref{eq:matrix_F}.
	Since \mbox{$\b b_{\b\ell}^{\mathrm{opt}}\in\R^{|\I_{\b M_{\b\sigma},m}|}$} as solutions to the least squares problems \eqref{eq:opt_overdet_columnwise} satisfy \eqref{eq:res_opt_columnwise}, we can find \new{\mbox{$\varepsilon \coloneqq \max_{\b\ell\in\I_{\b M_{\b\sigma}}} \varepsilon_{\b\ell} \geq 0$}, such that \mbox{$\varepsilon_{\b\ell} \leq \varepsilon$}, \mbox{$\b\ell\in\I_{\b M_{\b\sigma}}$}, and thereby}
	\begin{align*}
		\sum_{\b\ell \in \I_{\b M_{\b\sigma}}} \big\|\b H_{\b\ell} \b b_{\b\ell}^{\mathrm{opt}} - \tfrac{1}{|\I_{\b M_{\b\sigma}}|} \b D^{-1} \b f_{\b\ell}\big\|_2^2 
		\leq \sum_{\b\ell \in \I_{\b M_{\b\sigma}}} \varepsilon_{\b\ell}
		\leq \varepsilon \, |\I_{\b M_{\b\sigma}}| .
	\end{align*}
	Therefore, we may write \eqref{eq:submult_fro_opt} as
	\begin{align}
		\label{eq:submult_fro_opt_1}
		\big\| \b D^* \b F^* \b B_{\mathrm{opt}}^* \b A -\b I_{|\I_{\b M}|} \big\|_{\mathrm F}^2
		&\leq
		\varepsilon \,|\I_{\b M_{\b\sigma}}| \cdot 
		\,\big\| \b F \b D \big\|_{\mathrm F}^2.
	\end{align}
	
	Thus, it remains to estimate the Frobenius norm \mbox{$\big\| \b F \b D \big\|_{\mathrm F}^2$}.
	By the definitions of the Frobenius norm and the trace $\mathrm{tr}(\b Z)$ of a matrix $\b Z$, it is clear that 
	\mbox{$\|\b Z\|_{\mathrm F}^2 = \mathrm{tr}(\b Z^* \b Z)$}.
	Since by \eqref{eq:matrix_F} we have that \mbox{$\b F^* \b F = |\I_{\b M_{\b\sigma}}|\, \b I_{|\I_{\b M}|}$}, this yields
	\begin{align}
		\label{eq:norm_FD}
		\big\| \b F \b D \big\|_{\mathrm F}^2
		=
		\mathrm{tr}(\b D^* \b F^* \b F \b D)
		=
		|\I_{\b M_{\b\sigma}}|\cdot \mathrm{tr}(\b D^* \b D)
		=
		|\I_{\b M_{\b\sigma}}|\cdot \big\| \b D \big\|_{\mathrm F}^2 .
	\end{align}
	Using the definition \eqref{eq:matrix_D} of the diagonal matrix \mbox{$\b D\in\R^{|\I_{\b{M}}| \times |\I_{\b{M}}|}$}, we obtain 
	\begin{align}
		\label{eq:norm_D}
		\big\| \b D \big\|_{\mathrm F}^2
		=
		\frac 1{|\I_{\b M_{\b\sigma}}|^2} \,\sum_{\b k \in \I_{\b M}} \frac 1{\hat{w}(\b k)^2} .
	\end{align}
	Then combination of \eqref{eq:submult_fro_opt_1}, \eqref{eq:norm_FD} and \eqref{eq:norm_D} implies
	\begin{align*}
		\big\| \b D^* \b F^* \b B_{\mathrm{opt}}^* \b A -\b I_{|\I_{\b M}|} \big\|_{\mathrm F}^2
		&\leq
		\varepsilon \,\sum_{\b k \in \I_{\b M}} \frac 1{\hat{w}(\b k)^2} ,
	\end{align*}
	such that \eqref{eq:est_approx_error_opt} yields the assertion \eqref{eq:approx_error_opt}.
	
	Since it is known that $0\leq \hat w(\b k) \leq 1$, $\b k\in\I_{\b M}$, for suitable window functions of the NFFT, cf.~\cite{PoTa21a}, we have
	\begin{align*}
		1 \leq \frac 1{\hat{w}(\b k)} \leq \frac 1{\hat{w}(\b k)^2} 
	\end{align*}
	and therefore
	\begin{align*}
		\sum_{\b k \in \I_{\b M}} \frac 1{\hat{w}(\b k)^2}
		\geq
		\sum_{\b k \in \I_{\b M}} 1
		=
		|\I_{\b M}|.
	\end{align*}
	Hence, the smallest constant is achieved in \eqref{eq:approx_error_opt} when \mbox{$\hat w(\b k) = 1$}, \mbox{$\b k\in\I_{\b M}$}, i.\,e., the (asymmetric) Dirichlet kernel \eqref{eq:kernel_Dirichlet_asym} is the optimal window function for the inverse NFFT in Algorithm~\ref{alg:infft_overdet}.
\end{proof}

Note that for trigonometric polynomials \eqref{eq:trig_poly_2d} 
the error bound of Theorem~\ref{Thm:error_est_opt} with the optimal window function \eqref{eq:kernel_Dirichlet_asym} is the same as the error bound \eqref{eq:error_est_trig_poly} from \autoref{Thm:error_est_lgs}.

\begin{Remark}
	\label{Remark:overdet_adj}
	Up to now, we only focused on the problem~\eqref{eq:problem_infft}.
	Finally, considering the inverse adjoint NFFT in~\eqref{eq:problem_infft*}, we remark that this problem can also be solved by means of the optimization procedure in Algorithm~\ref{alg:opt_fast_overdet}.
	Assuming again \mbox{$|\I_{\b M}| \leq N$}, this is the underdetermined setting for the adjoint problem \eqref{eq:problem_infft*}.
	Therefore, the minimum norm solution of \eqref{eq:problem_infft*} is given by the normal equations of second kind 
	$$\b A^* \b A \b y = \b h, \quad \b f = \b A \b y.$$
	Incorporating the matrix decomposition of the NFFT, cf.~Section~\ref{subsec:nfft}, we recognize that a modification of the matrix~\mbox{$\b B\in\R^{N\times |\I_{\b M_{\b\sigma}}|}$} such that \mbox{$\b A^* \b B \b F \b D \approx \b I_{|\I_{\b M}|}$} implies \mbox{$\b y \approx \b h$} and hence \mbox{$\b f \approx \b B \b F \b D \b h$}.
	Thus, the optimization problem~\eqref{eq:opt_id_adj} is also the one to consider for \eqref{eq:problem_infft*}. 
	In other words, our approach provides both, an inverse NFFT as well as an inverse adjoint NFFT.
	\ex
\end{Remark}

\section{Numerics \label{sec:numerics}}

Concluding, we have a look at some numerical examples.
Besides comparing the density compensation approach from Section~\ref{sec:inv_density_comp} to the optimization approach from Section~\ref{sec:opt_B},
for both trigonometric polynomials \eqref{eq:trig_poly_2d} and bandlimited functions,
we also demonstrate the accuracy of these approaches.

\begin{Remark}
	\label{Remark:grids}
	At first we introduce some exemplary grids.
	For visualization we restrict ourselves to the two-dimensional setting \mbox{$d=2$}.
	\begin{enumerate}[label=\emph{(\roman*)}]
		\item
		We start with a sampling scheme that is somehow ``close'' to the Cartesian grid, but also possesses a random part.
		To this end, we consider the two-dimensional Cartesian grid and add a two-dimensional perturbation, i.\,e.,
		\begin{equation}
			\label{eq:jittered_grid}
			x_{t,j} \coloneqq \left( -\frac 12 + \frac{2t-1}{N_1}, -\frac 12 + \frac{2j-1}{N_2} \right)^T 
			+ \left(\frac{1}{N_1}\,\eta_1, \frac{1}{N_2}\,\eta_2\right)^T,
		\end{equation}
		\mbox{$t=1,\dots,N_1$}, \mbox{$j=1,\dots,N_2$}, and \mbox{$\eta_1,\eta_2 \sim U(-1,1)$}, where $U(-1,1)$ denotes the uniform distribution on the interval~$(-1,1)$.
		A visualization of this jittered grid can be found in Figure~\ref{fig:jittered_grids_jitter}.
		Additionally, we also consider the random grid
		\mbox{$x_{t,j} \coloneqq \frac 12 \left( \eta_1, \eta_2 \right)^T,$}
		see Figure~\ref{fig:jittered_grids_rand}.
		\begin{figure}[!h]
			\centering
			\captionsetup[subfigure]{justification=centering}
			\begin{subfigure}[t]{0.3\textwidth}
				\centering
				\includegraphics[width=\textwidth]{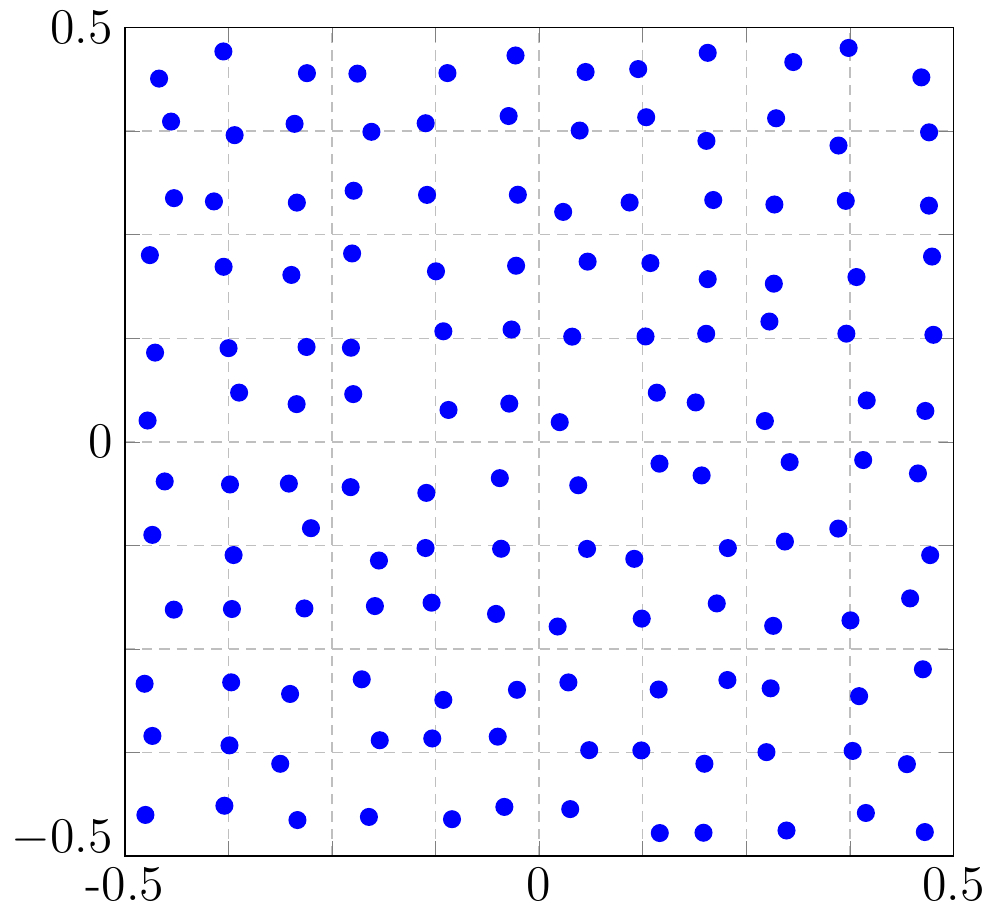}
				\caption{Jittered grid}
				\label{fig:jittered_grids_jitter}
			\end{subfigure}
			\begin{subfigure}[t]{0.3\textwidth}
				\centering
				\includegraphics[width=\textwidth]{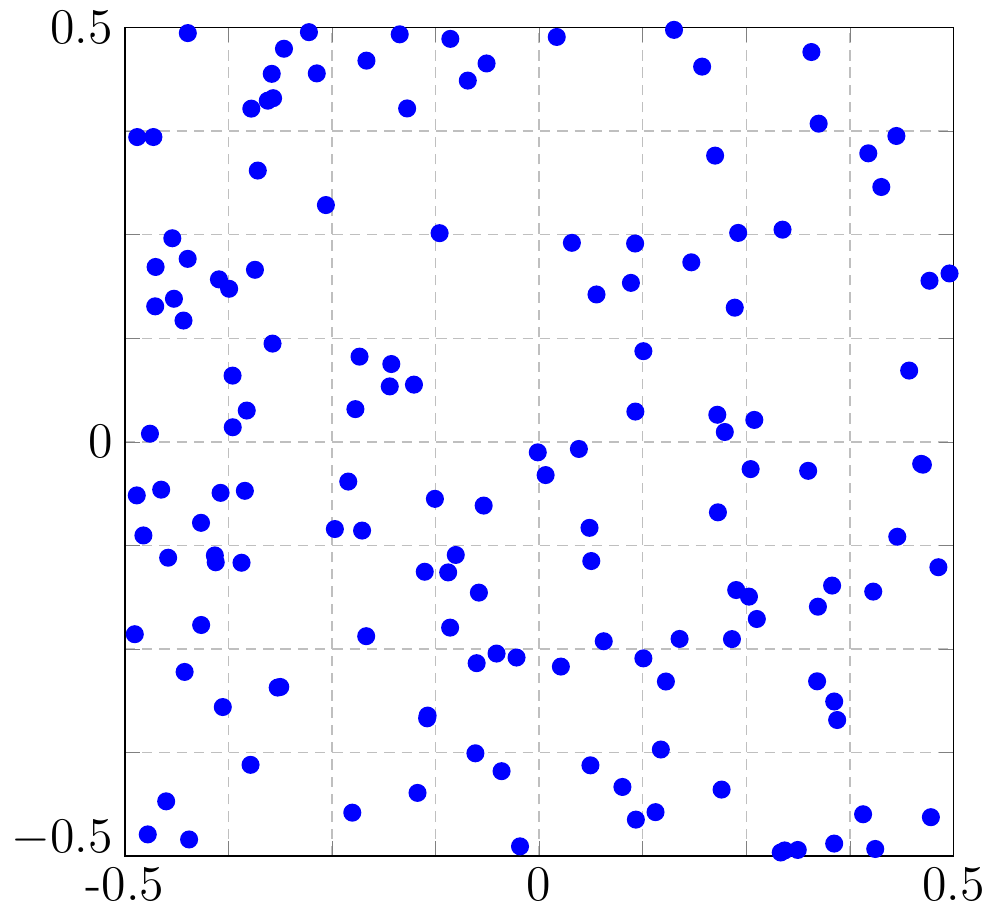}
				\caption{Random grid}
				\label{fig:jittered_grids_rand}
			\end{subfigure}
			\caption{
				Exemplary randomized grids of size $N_1=N_2=12$.
				\label{fig:jittered_grids}}
		\end{figure}
		\item 
		Moreover, we examine grids of polar kind, as mentioned in \cite{FeKuPo06}.
		For \mbox{$R,T\in 2\N$} the points of the \textit{polar grid} are given by a signed radius 
		\mbox{$r_j \coloneqq \frac jR \in \left[-\frac 12,\frac 12\right)$} 
		and an angle 
		\mbox{$\theta_t \coloneqq \frac{\pi t}{T} \in \left[-\frac{\pi}{2},\frac{\pi}{2}\right)$} as
		\begin{equation}
			\label{eq:polar_grid}
			x_{t,j} \coloneqq r_j \left(\cos \theta_t, \sin \theta_t\right)^T, \quad (j,t)^T \in \I_R \times \I_T.
		\end{equation}
		Since it is known that the inversion problem is ill-conditioned for this grid we consider a modification, the \textit{modified polar grid}
		\begin{equation}
			\label{eq:mpolar_grid}
			x_{t,j} \coloneqq r_j \left(\cos \theta_t, \sin \theta_t\right)^T, \quad (j,t)^T \in \I_{\sqrt 2R} \times \I_T,
		\end{equation}
		i.\,e., we added more concentric circles and excluded the points outside the unit square, see Figure~\ref{fig:polar_grids_mpolar}.
		Another sampling scheme which is known to lead to more stable results than the polar grid is the \textit{linogram} or \textit{pseudo-polar grid}, where the points lie on concentric squares instead of concentric circles, see Figure~\ref{fig:polar_grids_linogram}.
		There we distinguish two sets of points, i.\,e., 
		\begin{equation}
			\label{eq:linogram_grid}
			x_{t,j}^{BH} \coloneqq \left( \frac jR, \frac{4t}{T} \frac jR \right)^T,
			\quad 
			x_{t,j}^{BV} \coloneqq \left( -\frac{4t}{T} \frac jR, \frac jR \right)^T,
			\quad 
			(j,t)^T \in \I_R \times \I_{\frac T2}.
		\end{equation}
		\begin{figure}[!h]
			\centering
			\captionsetup[subfigure]{justification=centering}
			\begin{subfigure}[t]{0.3\textwidth}
				\centering
				\includegraphics[width=\textwidth]{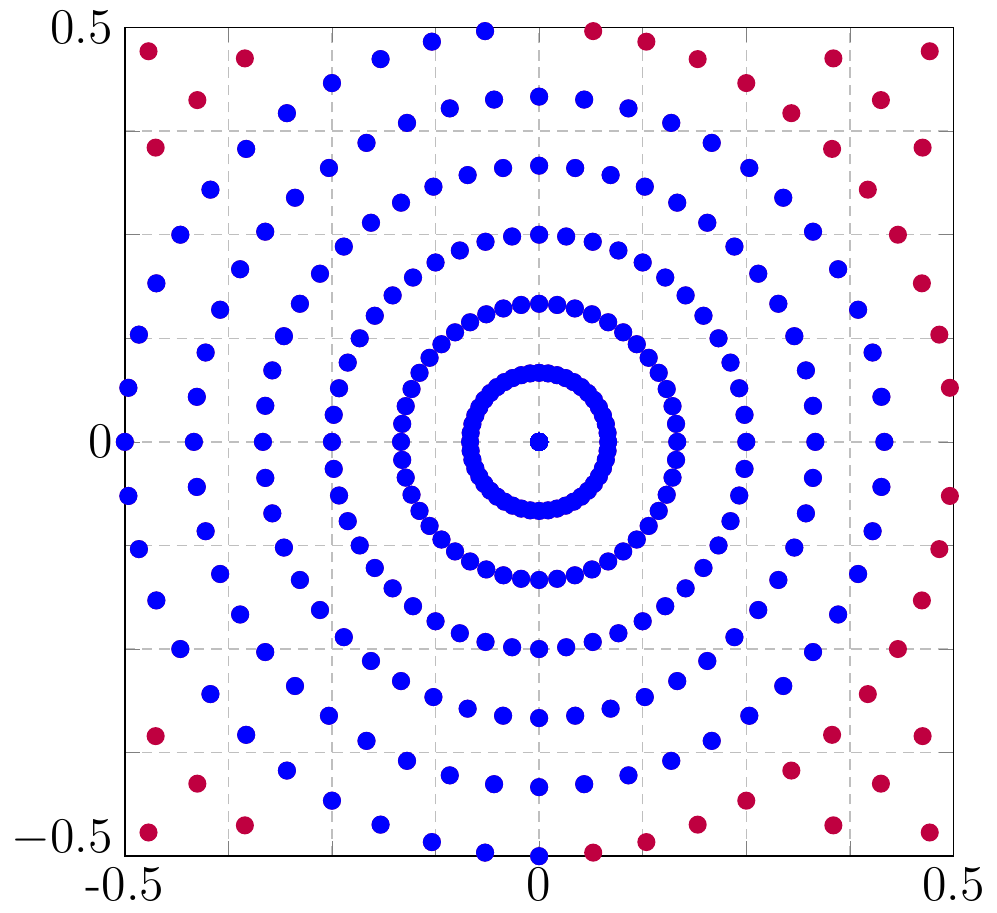}
				\caption{Polar (blue) and modified polar (red) grid}
				\label{fig:polar_grids_mpolar}
			\end{subfigure}
			\begin{subfigure}[t]{0.3\textwidth}
				\centering
				\includegraphics[width=\textwidth]{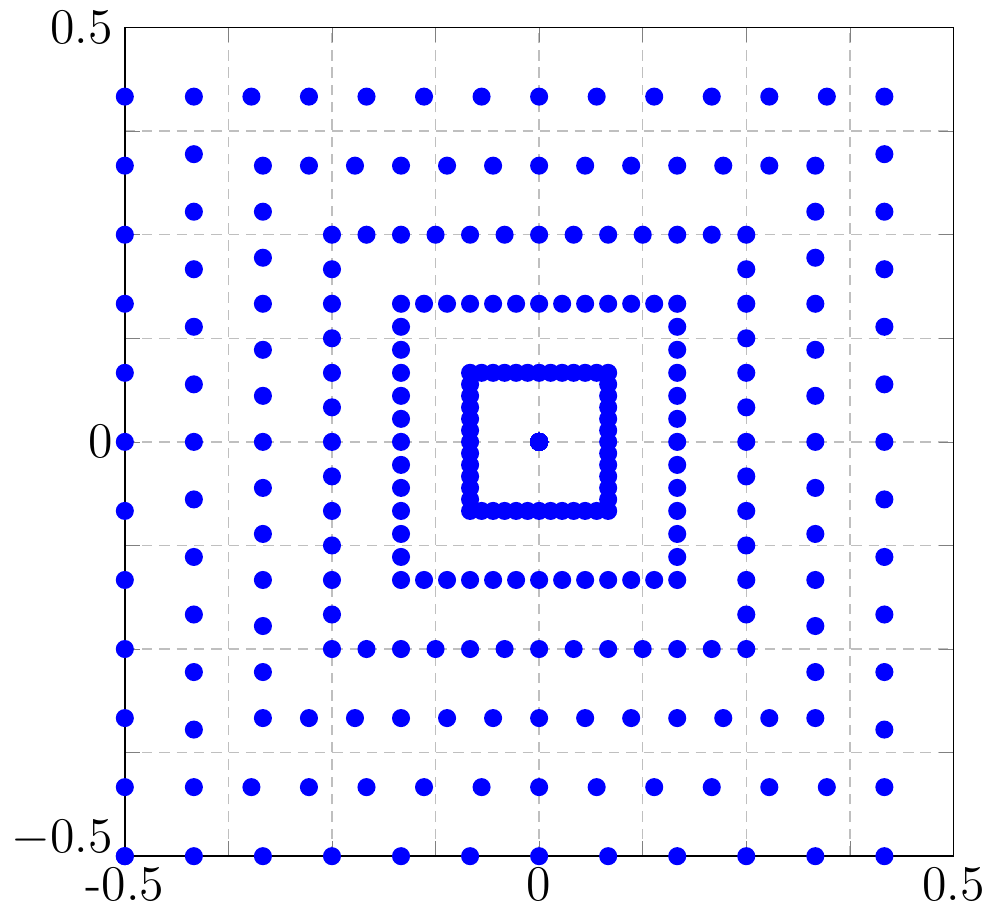}
				\caption{Linogram / pseudo-polar grid}
				\label{fig:polar_grids_linogram}
			\end{subfigure}
			\begin{subfigure}[t]{0.3\textwidth}
				\centering
				\includegraphics[width=\textwidth]{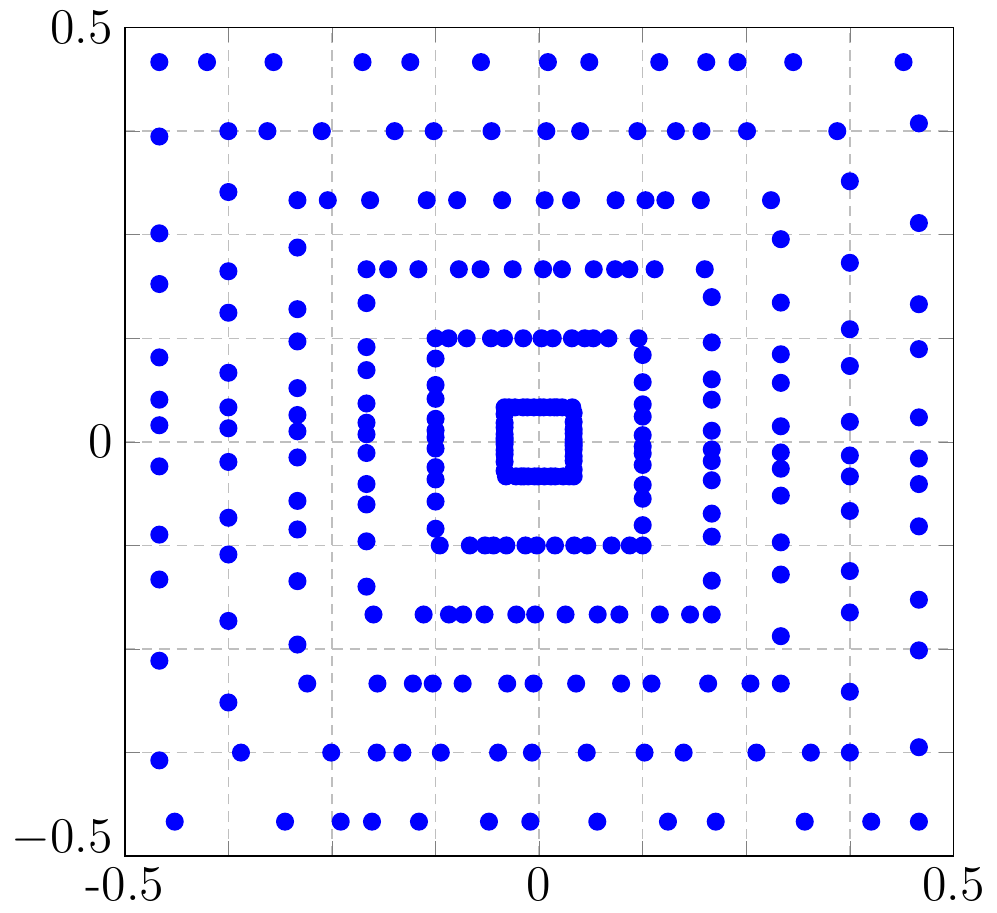}
				\caption{Golden angle linogram grid}
				\label{fig:polar_grids_golden}
			\end{subfigure}
			\caption{
				Polar grids of size $R=12$ and $T=2R$.
				\label{fig:polar_grids}}
		\end{figure}
		\item 
		Another modification of these polar grids was introduced in \cite{He19}, where the angles~$\theta_t$ are not chosen equidistantly but are obtained by golden angle increments.
		For the \textit{golden angle polar grid} we only exchange the equispaced angles of the polar grid to 
		\begin{equation}
			\label{eq:golden_angles}
			\theta_t = \mod \left( \frac{\pi}{2} + t\,\frac{2\pi}{1+\sqrt{5}} , \pi \right) - \frac{\pi}{2},\quad t=0,\dots,T-1.
		\end{equation}
		The \textit{golden angle linogram grid} is given by
		\begin{equation*}
			x_{t,j} \coloneqq \left\{ \begin{matrix}
				\left( \frac{2j+1}{2R}, \,\frac{2j+1}{2R} \tan\left(\theta_t-\frac{\pi}{4}\right) \right)^T & \colon \theta_t \in \left[0,\frac{\pi}{2}\right) \\
				\left( -\frac{2j+1}{2R} \cot\left(\theta_t-\frac{\pi}{4}\right), \,\frac{2j+1}{2R} \right)^T & \colon \theta_t \in \left[-\frac{\pi}{2},0\right) \\
			\end{matrix} \right\}, 
			\quad j \in \I_R,
		\end{equation*}
		with $\theta_t$ in \eqref{eq:golden_angles}, as illustrated in Figure~\ref{fig:polar_grids_golden}.
		\ex
	\end{enumerate}
\end{Remark}

Before comparing the different approaches from Sections~\ref{sec:inv_density_comp} and \ref{sec:opt_B}, we study the quality of our methods for the grids mentioned in Remark~\ref{Remark:grids}.
More specifically, in Example~\ref{ex:err_trig_poly_2d} we investigate the accuracy of the density compensation method from Algorithm~\ref{alg:infft_density} with the weights introduced in Section~\ref{sec:trig_poly}, 
and in Example~\ref{ex:min_norm_overdet_2d} we check if the norm minimization targeted in Section~\ref{sec:opt_B} is successful.

\begin{Example}
	\label{ex:err_trig_poly_2d}
	Firstly, we examine the quality of our density compensation method in Algorithm~\ref{alg:infft_density} 
	for a trigonometric polynomial $f$ as in \eqref{eq:trig_poly_2d} 
	with given Fourier coefficients \mbox{$\hat f_{\b k} \in [1,10]$}, \mbox{$\b k\in\I_{\b{M}}$}. 
	\new{In this test we consider several $d\in \{1,2,3\}$}.
	For the corresponding function evaluations of \eqref{eq:trig_poly_2d} at given points \new{\mbox{$\b x_j \in \left[-\frac 12,\frac 12\right)^d$}}, \mbox{$j=1,\dots,N$},
	we test how well these Fourier coefficients can be approximated.
	More precisely, we consider the estimate
	\mbox{$\b{\tilde h}^{\mathrm{w}} = \b D^* \b F^* \b B^* \b W \b f$}, cf.~\eqref{eq:reconstr_density_nfft},
	with the matrix \mbox{$\b W = \diag\left(w_j\right)_{j=1}^N$} of density compensation factors computed by means of \new{Algorithm~\ref{alg:precompute_density}}, i.\,e., by \eqref{eq:normal_equations_second_kind_double}, in case \mbox{$|\I_{\b {2M}}| \leq N$}, or by \eqref{eq:normal_equations_first_kind_double}, if \mbox{$|\I_{\b {2M}}| > N$}, 
	and compute the relative errors
	\begin{equation}
		\label{eq:errors}
		e_p \coloneqq \frac{\|\b{\tilde h}^{\mathrm{w}}-\b{\hat f}\|_{p}}{\|\b{\hat f}\|_{p}} ,
		\quad p\in\{2,\infty\} .
	\end{equation}
	\new{By \eqref{eq:error_est_trig_poly} it is known that
		\begin{align*}
			\frac{\big\| \b{\hat f} - \b A^* \b W \b f \big\|_p}{\big\| \b{\hat f} \big\|_p}
			&\leq
			|\I_{\b{M}}| \, \varepsilon ,
			\quad p \in \{2,\infty\} , 
		\end{align*}
		with the residual \mbox{$\varepsilon = \big\| \b A_{|\I_{\b {2M}}|}^T \,\b w - \b e_{\b 0} \big\|_{\infty} \geq 0$}, cf.~\eqref{eq:residual_lgs}.}
	
	In our experiment we use random points $\b x_j$ with \new{\mbox{$N_t = 2^{9-d}$}, \mbox{$t=1,\dots,d$}}, cf.~Figure~\ref{fig:jittered_grids_rand}, and, for several problem sizes \new{\mbox{$\b M=M \cdot \b 1_d$}}, \mbox{$M=2^c$} with \mbox{$c=1,\dots,\new{11-d},$}
	we choose random Fourier coefficients \mbox{$\hat f_{\b k} \in [1,10]$}, \mbox{$\b k\in\I_{\b{M}}$}.
	Afterwards, we compute the evaluations of the trigonometric polynomial \eqref{eq:trig_poly_2d} by means of an NFFT
	and use the resulting vector~$\b f$ as input for the reconstruction.
	Due to the randomness we repeat this 10 times and then consider the \new{maximum} error over all runs. 
	The corresponding results are displayed in Figure~\ref{fig:err_trig_poly_2d}.
	It can clearly be seen that \mbox{$|\I_{\b{2M}}|<N$}, i.\,e., as long as \new{\mbox{$M<\frac{N_1}{2}=2^{8-d}$}}, the weights computed by means of \eqref{eq:normal_equations_second_kind_double} lead to an exact reconstruction of the given Fourier coefficients.
	However, as soon as we are in the setting \mbox{$|\I_{\b{2M}}|>N$} the least squares approximation via \eqref{eq:normal_equations_first_kind_double} does not yield good results anymore.
	\begin{figure}[ht]
		\centering
		\captionsetup[subfigure]{justification=centering}
		\begin{subfigure}[t]{0.32\textwidth}
			\includegraphics[width=\textwidth]{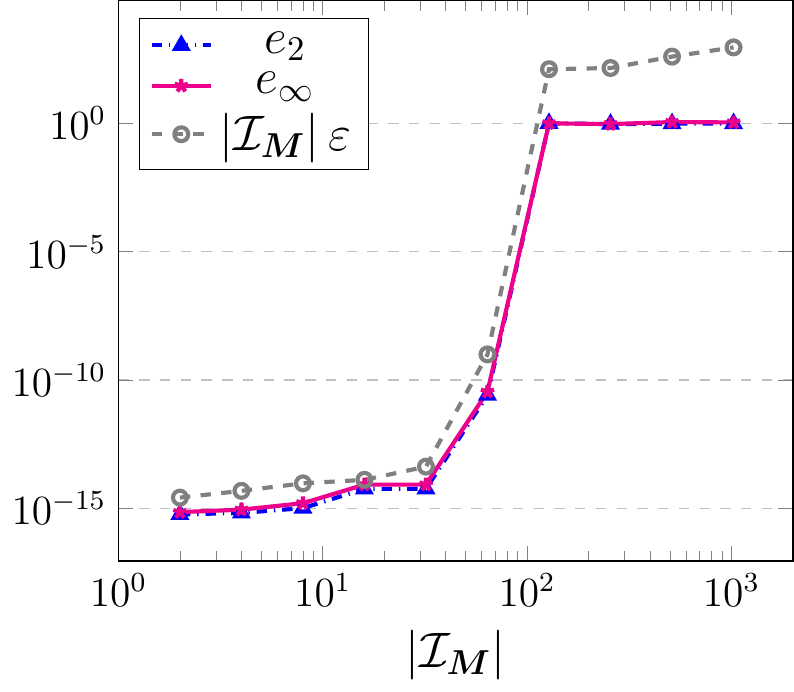}
			\caption{$d=1$}
		\end{subfigure}
		\begin{subfigure}[t]{0.32\textwidth}
			\includegraphics[width=\textwidth]{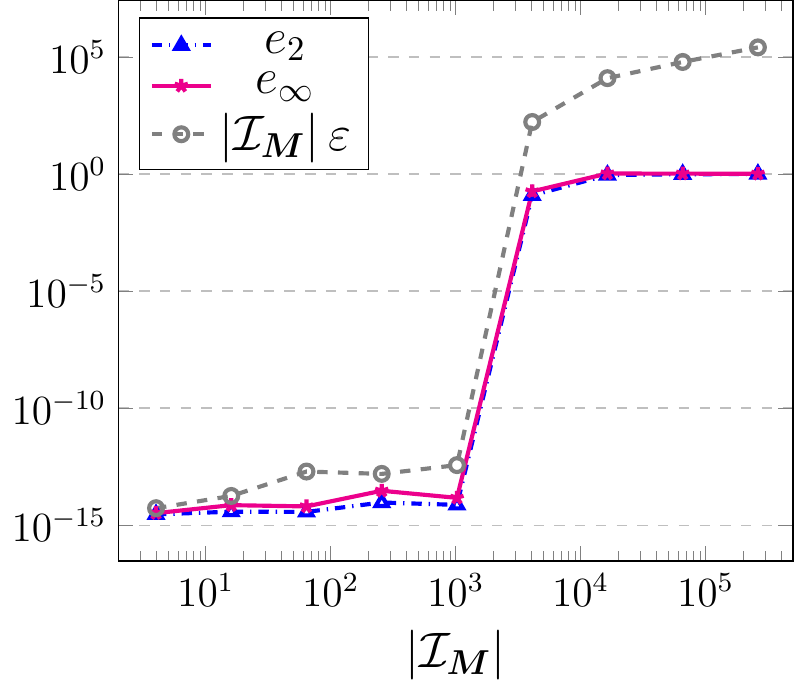}
			\caption{$d=2$}
		\end{subfigure}
		\begin{subfigure}[t]{0.32\textwidth}
			\includegraphics[width=\textwidth]{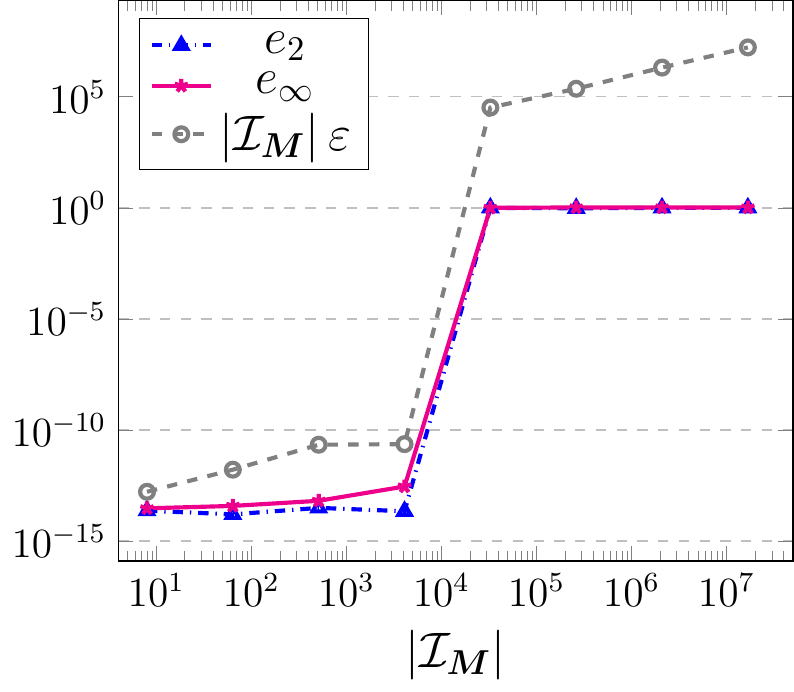}
			\caption{$d=3$}
		\end{subfigure}
		\caption{
			Relative errors~\eqref{eq:errors} of the reconstruction of the Fourier coefficients of a trigonometric polynomial \eqref{eq:trig_poly_2d} with given \mbox{$\hat f_{\b k} \in [1,10]$}, \mbox{$\b k\in\I_{\b{M}}$}, computed via the density compensation method from Algorithm~\ref{alg:infft_density}, for random grids with \new{\mbox{$N_t = 2^{9-d}$}, \mbox{$t=1,\dots,d$}}, and \new{\mbox{$\b M=M \cdot \b 1_d$}}, \mbox{$M=2^c$} with \mbox{$c=1,\dots,\new{11-d}$}.
			\label{fig:err_trig_poly_2d}}
	\end{figure}
	\ex
\end{Example}

\begin{Example}
	\label{ex:min_norm_overdet_2d}
	In order to study the quality of our optimization method in Section~\ref{sec:opt_B}, we consider the Frobenius norms 
	\begin{equation}
		\label{eq:matrixnorms_2d}
		n(w,m,\sigma) \coloneqq \left\| \b A^* \b B \b F \b D - \b I_{|\I_{\b M}|} \right\|_{\mathrm F},
		\quad
		n^{\textrm{opt}}(w,m,\sigma) \coloneqq \left\| \b A^* \b B_{\textrm{opt}} \b F \b D - \b I_{|\I_{\b M}|} \right\|_{\mathrm F},
	\end{equation}
	where $\b B$ denotes the original matrix from the NFFT in \eqref{eq:matrix_B} and $\b B_{\mathrm{opt}}$ the optimized matrix generated by Algorithm~\ref{alg:opt_fast_overdet}.
	For the original matrix $\b B$ we utilize the common \mbox{B-Spline} window function
	\begin{equation}
		\label{eq:window_bspline}
		w_{\mathrm{B}} \coloneqq B_{2m}(\b M_{\b\sigma} \b x)
	\end{equation}
	with the centered \mbox{B-Spline} of order $2m$, cf. \cite[p.~388]{PlPoStTa18}.
	The optimized matrix $\b B_{\mathrm{opt}}$ shall be computed by means of the \mbox{B-Spline} \eqref{eq:window_bspline} as well as the Dirichlet window function~\eqref{eq:kernel_Dirichlet_asym}, which is the optimal window by Theorem~\ref{Thm:error_est_opt}.
	
	\new{Due to memory limitations in the computation of the Frobenius norms \eqref{eq:matrixnorms_2d}, we have to settle for very small problems, which however show the functionality of Algorithm~\ref{alg:opt_fast_overdet}.
		For this reason we consider $d=2$ and choose \mbox{$\b M=(12,12)^T$}} as well as \mbox{$N_1=N_2=R=2^\mu$}, \mbox{$\mu\in\{2,\dots,7\}$}, and \mbox{$T=2R$} for the grids mentioned in Remark~\ref{Remark:grids}.
	In other words, we test Algorithm~\ref{alg:opt_fast_overdet} in the underdetermined setting \mbox{$|\I_{\b M}| > N$} as well as for the overdetermined setting \mbox{$|\I_{\b M}| \leq N$}.
	
	Having a look at the results for the grids in Remark~\ref{Remark:grids}, it becomes apparent that they separate into two groups.
	Figure~\ref{fig:min_norm_overdet_polar} displays the results for the polar grid~\eqref{eq:polar_grid}, which are the same as for the golden angle polar grid, cf.~\eqref{eq:golden_angles}.
	In these cases there is only a slight improvement by the optimization.
	However, for all other mentioned grids the minimization procedure in Algorithm~\ref{alg:opt_fast_overdet} is very effective.
	The results for these grids are depicted in Figure~\ref{fig:min_norm_overdet_mpolar} exemplarily for the modified polar grid~\eqref{eq:mpolar_grid}.
	Moreover, it can be seen that our optimization procedure in Algorithm~\ref{alg:opt_fast_overdet} is most effective in the overdetermined setting \mbox{$|\I_{\b M}| \leq N$}. 
	%
	\begin{figure}[ht]
		\centering
		\captionsetup[subfigure]{justification=centering}
		\begin{subfigure}{0.496\textwidth}
			\includegraphics[width=\textwidth]{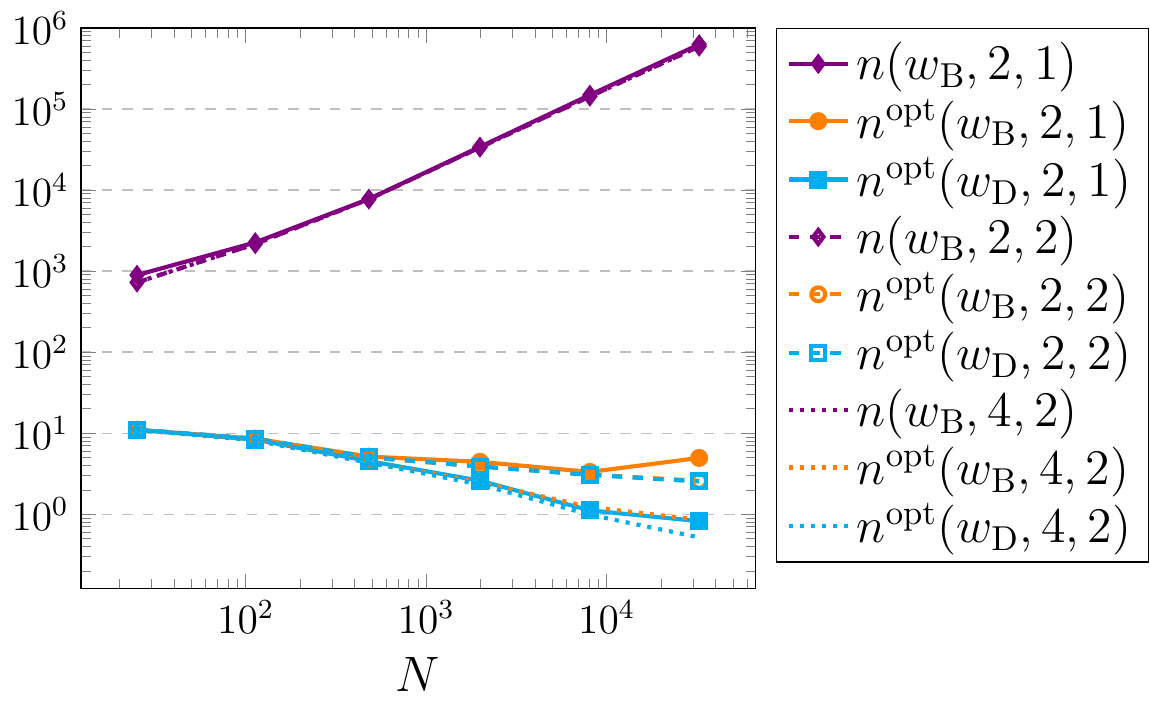}
			\subcaption{Polar grid~\eqref{eq:polar_grid}}
			\label{fig:min_norm_overdet_polar}
		\end{subfigure}
		\begin{subfigure}{0.496\textwidth}
			\includegraphics[width=\textwidth]{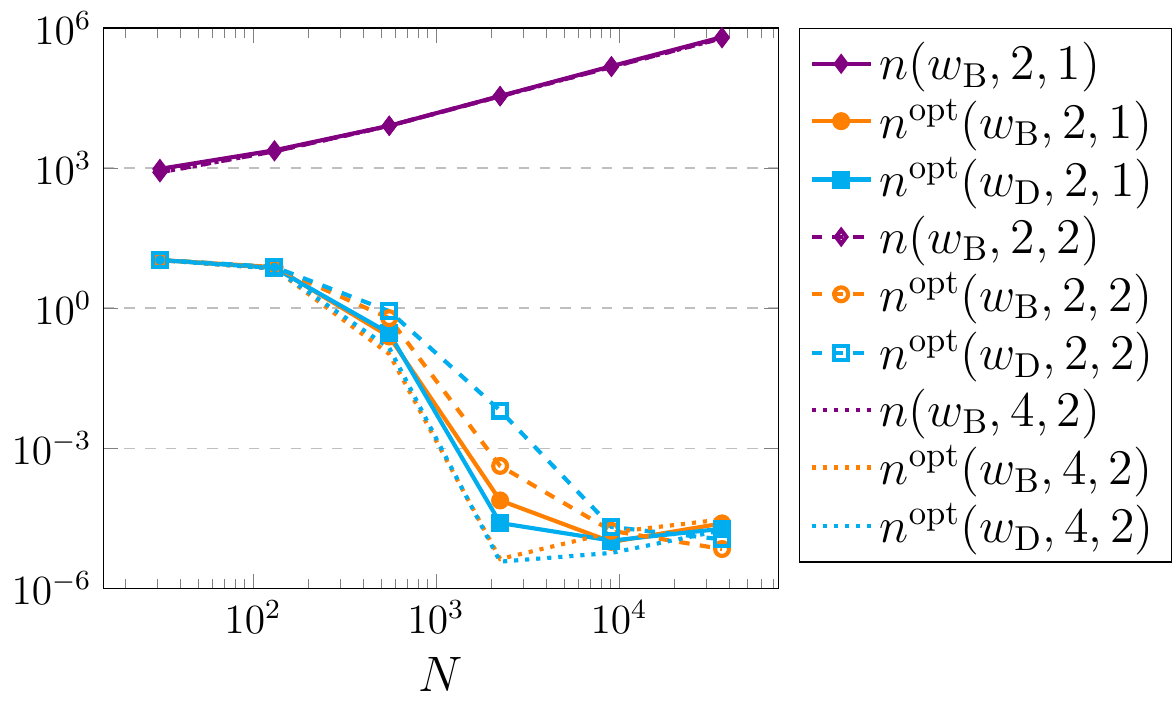}
			\subcaption{Modified polar grid~\eqref{eq:mpolar_grid}}
			\label{fig:min_norm_overdet_mpolar}
		\end{subfigure}
		\caption{
			Frobenius norms~\eqref{eq:matrixnorms_2d} of the original matrix~$\b B$ (violet) and the optimized matrix~$\b B_{\textrm{opt}}$ generated by Algorithm~\ref{alg:opt_fast_overdet} using the {B-Spline} $w_{\mathrm{B}}$ (orange) 
			as well as the Dirichlet window $w_{\mathrm{D}}$ (cyan) 
			with \mbox{$R=2^\mu$}, \mbox{$\mu\in\{2,\dots,7\}$}, and \mbox{$T=2R$} as well as \mbox{$\b M=(12,12)^T$}, \mbox{$m\in\{2,4\}$} and \mbox{$\sigma\in\{1,2\}$}.
			\label{fig:min_norm_overdet}}
	\end{figure}
	
	One reason for the different behavior of polar and modified polar grid could be the ill-posedness of the inversion problem for the polar grid, which becomes evident in huge condition numbers of $\b H_{\b\ell}^* \b H_{\b\ell}$, whereas the problem for modified polar grids is well-posed.
	Another reason can be found in the optimization procedure itself. 
	Having a closer look at the polar grid, see Figure~\ref{fig:polar_grids_mpolar}, there are no grid points in the corners of the unit square.
	Therefore, some of the index sets~
	\mbox{$\I_{\b M_{\b\sigma},m}(\b\ell)$}, cf.~\eqref{eq:indexset_l}, 
	are empty and no optimization can be done for the corresponding matrix columns.
	This could also cause the worsened minimization properties of the polar grid.
	\ex
\end{Example}

Next, we proceed with comparing the density compensation approach from Section~\ref{sec:inv_density_comp} using the weights \mbox{$w_j$} introduced in Section~\ref{sec:trig_poly} to the optimization approach for modifying the matrix $\b B$ from Section~\ref{sec:opt_B}.
To this end, we show an example concerning trigonometric polynomials \eqref{eq:trig_poly_2d} of degree~$\b M$ and a second one that deals with bandlimited functions of bandwidth~$\b M$.
\new{Here we restrict ourselves to the two-dimensional setting \mbox{$d=2$} for better visualization of the results.}

\begin{Example}
	\label{ex:reconstr_phantom}
	Similar to \cite{AvCoDoIsSh08, KiPo20} we have a look at the reconstruction of the Shepp-Logan phantom, see Figure~\ref{fig:reconstr_phantom_linogram_original}.
	Here we treat the phantom data as given Fourier coefficients \mbox{$\b{\hat f}\coloneqq(\hat f_{\b k})_{\b k \in \I_{\b M}}$} of a trigonometric polynomial \eqref{eq:trig_poly_2d}.
	For given points \mbox{$\b x_j \in \left[-\frac 12,\frac 12\right)^2$}, \mbox{$j=1,\dots,N$}, we then compute the evaluations of the trigonometric polynomial \eqref{eq:trig_poly_2d} by means of an NFFT
	and use the resulting vector as input for the reconstruction.
	
	In a first experiment, we test the inversion methods from Sections~\ref{sec:inv_density_comp} and \ref{sec:opt_B} as in \cite{AvCoDoIsSh08} for increasing input sizes.
	To this end, we choose \mbox{$\b M = (M,M)^T$}, \mbox{$M=2^c$} with \mbox{$c=3,\dots,\new{10},$}
	and linogram grids~\eqref{eq:linogram_grid} of size \mbox{$R=2M$}, \mbox{$T=2R$}, i.\,e., we consider the setting \mbox{$|\I_{\b{2M}}|<N$}.
	For using Algorithm~\ref{alg:infft_overdet} we choose the oversampling factor \mbox{$\sigma = 1.0$} and the truncation parameter \mbox{$m = 4$}.
	For each input size we measure the computation time of the precomputational steps, i.\,e., the computation of the weight matrix $\b W$ or the computation of the optimized sparse matrix~$\b B_{\textrm{opt}}$, as well as the time needed for the reconstruction, i.\,e., the corresponding adjoint NFFT, see Algorithms~\ref{alg:infft_density} and \ref{alg:infft_overdet}. 
	Moreover, for the reconstruction \mbox{$\b{\tilde h} \in \{ \b{\tilde h}^{\mathrm{w}}, \b h_{\mathrm{opt}} \}$}, cf.~\eqref{eq:reconstr_density_nfft} and \eqref{eq:approx}, we consider the relative errors 
	\begin{equation}
		\label{eq:relative_errors}
		e_{2} \coloneqq \frac{\|\b{\tilde h}-\b{\hat f}\|_{2}}{\|\b{\hat f}\|_{2}}.
	\end{equation}
	The corresponding results can be found in Table~\ref{table:reconstr_phantom_linogram_2M}.
	We remark that since we are in the setting \mbox{$|\I_{\b{2M}}|<N$}, the density compensation method in Algorithm~\ref{alg:infft_density} with weights computed by \eqref{eq:normal_equations_second_kind_double} indeed produces nearly exact results.
	Although, our optimization procedure from Algorithm~\ref{alg:infft_overdet} achieves small errors as well, this reconstruction is not as good as the one by means of our density compensation method.
	
	Note that in comparison to \cite{AvCoDoIsSh08} our method in Algorithm~\ref{alg:infft_density} using density compensation produces errors of the same order, but is much more effective for solving several problems using the same points $\b x_j$ for different input values $\b f$.
	Since our precomputations have to be done only once in this setting, we strongly profit from the fact that we only need to perform an adjoint NFFT as reconstruction, which is very fast,
	whereas in \cite{AvCoDoIsSh08} they would need to execute their whole routine each time again.
	\begin{table}[ht]
		\centering
		\begin{tabular}{|c|c|c|c|c|c|c|}
			\hline
			& \multicolumn{2}{c|}{Relative error $e_2$} & \multicolumn{2}{c|}{Precomputation time} & \multicolumn{2}{c|}{Reconstruction time} \\
			\hline
			$M$ & Alg.~\ref{alg:infft_density} & Alg.~\ref{alg:infft_overdet} & Alg.~\ref{alg:infft_density} & Alg.~\ref{alg:infft_overdet} & Alg.~\ref{alg:infft_density} & Alg.~\ref{alg:infft_overdet} \\
			\hline
			8 & 1.3332e-15 & 6.8606e-14 & 9.8254e-02 & 1.9220e+00 & 6.2500e-04 & 2.2360e-03 \\ 
			16 & 7.2315e-15 & 1.5718e-07 & 1.6157e-01 & 8.3276e+00 & 2.5100e-03 & 3.4760e-03 \\ 
			32 & 2.3383e-14 & 4.5778e-07 & 3.3032e-01 & 4.3169e+01 & 3.1860e-03 & 7.4790e-03 \\ 
			64 & 2.5859e-14 & 4.7505e-07 & 3.4324e+00 & 2.4103e+02 & 5.0420e-03 & 4.9310e-03 \\ 
			128 & 7.9006e-14 & 5.9962e-07 & 9.4725e+00 & 1.2045e+03 & 2.9860e-02 & 5.5123e-02 \\ 
			256 & 2.6386e-13 & 4.0943e-06 & 3.8365e+01 & 5.8347e+03 & 6.6443e-02 & 6.7810e-01 \\ 
			512 & 1.0917e-12 & 2.0184e-06 & 1.4020e+02 & 2.9235e+04 & 2.1674e-01 & 3.2674e+00 \\ 
			\hspace{-3pt}\new{1024}\hspace{-3pt} & 4.2563e-12 & 1.3491e-05 & 7.2153e+02 & 1.4342e+05 & 7.4896e-01 & 1.6114e+01 \\
			\hline
		\end{tabular}
		\caption{
			Relative errors~\eqref{eq:relative_errors} of the reconstruction of the Shepp-Logan phantom of size~$M$ as well as the runtime in seconds 
			for the density compensation method from Algorithm~\ref{alg:infft_density} compared to Algorithm~\ref{alg:infft_overdet} with \mbox{$\sigma = 1.0$} and \mbox{$m = 4$}, 
			using linogram grids~\eqref{eq:linogram_grid} of size \mbox{$R=2M$},\, \mbox{$T=2R$}. 
			\label{table:reconstr_phantom_linogram_2M}}
	\end{table}

	As a second experiment we aim to decrease the amount of overdetermination, i.\,e., we want to keep the size \mbox{$|\I_{\b{M}}|$} of the phantom, but reduce the number $N$ of the points $\b x_j$, \mbox{$j=1,\dots,N$}.
	To this end, we now consider linogram grids~\eqref{eq:linogram_grid} of the smaller size \mbox{$R=M$}, \mbox{$T=2R$}, i.\,e., now we have \mbox{$|\I_{\b{2M}}|>N$}.
	The reconstruction of the phantom of size \mbox{$1024\times 1024$} is presented in Figure~\ref{fig:reconstr_phantom_linogram}~(top) including a detailed view of the 832nd row of this reconstruction (bottom).
	Despite the reconstruction via Algorithm~\ref{alg:infft_overdet} as well as the density compensation method using weights computed by means of \eqref{eq:normal_equations_first_kind_double}, we also considered the result using Voronoi weights.
	For all approaches we added the corresponding relative errors \eqref{eq:relative_errors} to Figure~\ref{fig:reconstr_phantom_linogram} as well.

	Due to the fact that the exactness condition \mbox{$|\I_{\b{2M}}|<N$} (cf.~Section~\ref{sec:underdet}) is violated, 
	it can be seen in Figure~\ref{fig:reconstr_phantom_linogram_exact} that the density compensation method using weights computed by means of \eqref{eq:normal_equations_first_kind_double} does not yield an exact reconstruction is this setting.
	On the contrary, we recognize that our optimization method, see Figure~\ref{fig:reconstr_phantom_linogram_opt}, achieves a huge improvement in comparison to the density compensation techniques in Figure~\ref{fig:reconstr_phantom_linogram_voronoi} and \ref{fig:reconstr_phantom_linogram_exact} since no artifacts are visible.
	Presumably, this arises because there are more degrees of freedom in the optimization of the matrix $\b B$ from Section~\ref{sec:opt_B} than with the density compensation techniques from Section~\ref{sec:inv_density_comp}, cf.~Remark~\ref{Remark:opt_density}.
	We remark that although the errors are not as small as in Table~\ref{table:reconstr_phantom_linogram_2M}, by comparing Figure~\ref{fig:reconstr_phantom_linogram_original} and \ref{fig:reconstr_phantom_linogram_opt} it becomes apparent that the differences are not even visible anymore.
	Note that for this result the number $N$ of points is ca. 4 times lower as for the results in depicted in Table~\ref{table:reconstr_phantom_linogram_2M}, i.\,e.,
	we only needed twice as much function values as Fourier coefficients, whereas e.\,g. in \cite{AvCoDoIsSh08} they worked with a factor of more than 4.
	\begin{figure}[ht]
		\centering
		\captionsetup[subfigure]{justification=centering}
		\begin{subfigure}[t]{0.24\textwidth}
			\centering 
			\includegraphics[width=\textwidth]{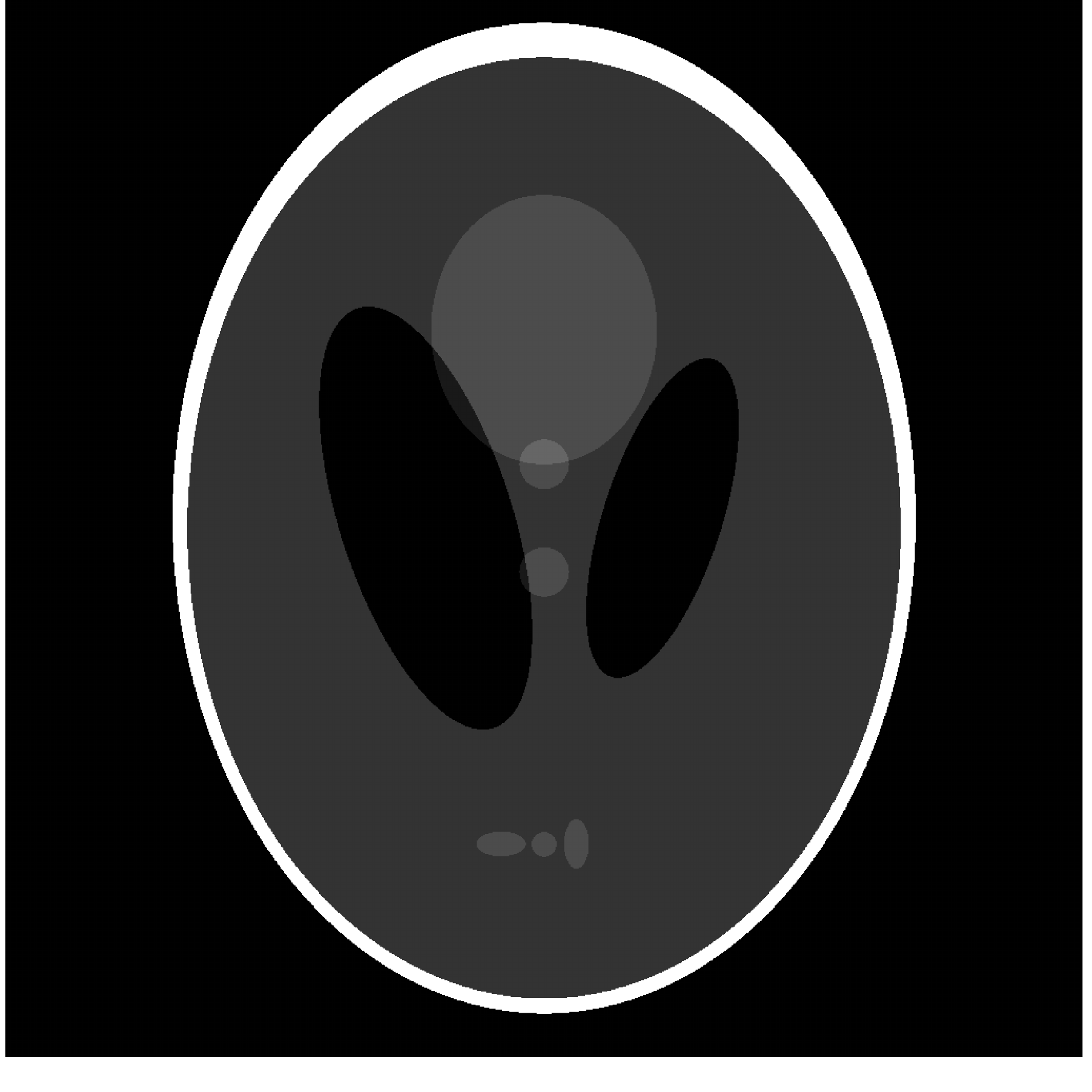}
		\end{subfigure}
		\begin{subfigure}[t]{0.24\textwidth}
			\centering
			\includegraphics[width=\textwidth]{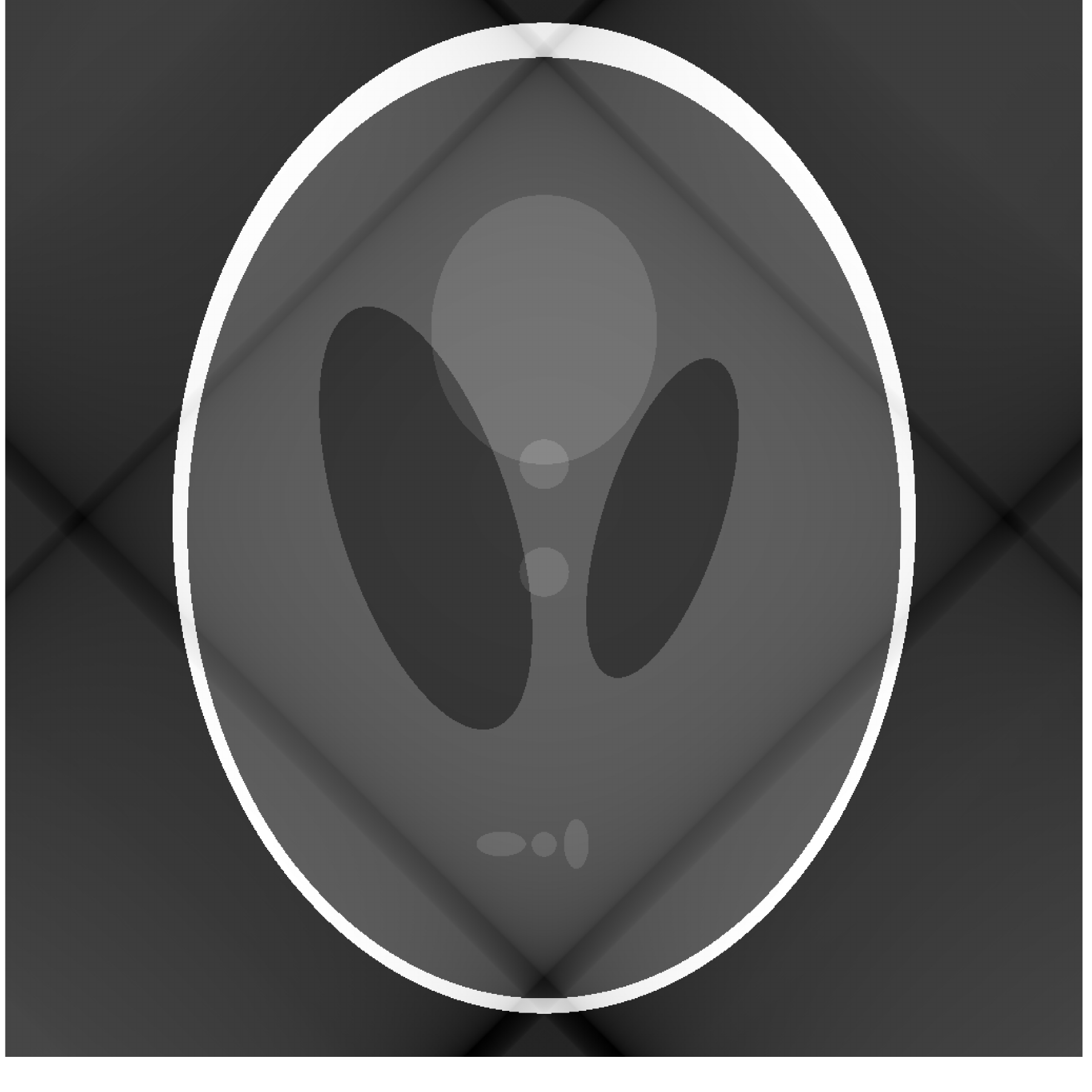}
		\end{subfigure}
		\begin{subfigure}[t]{0.24\textwidth}
			\centering
			\includegraphics[width=\textwidth]{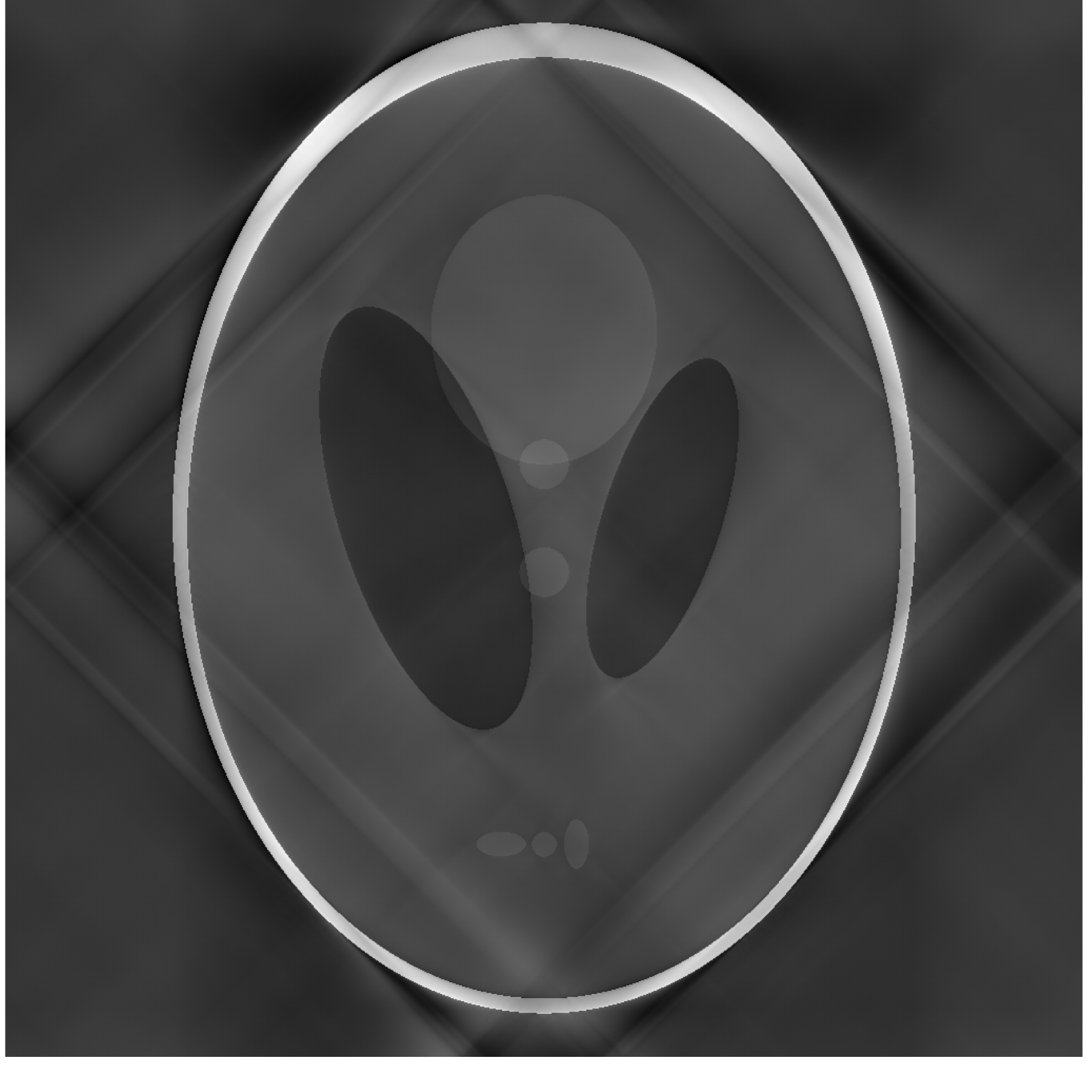}
		\end{subfigure}
		\begin{subfigure}[t]{0.24\textwidth}
			\centering
			\includegraphics[width=\textwidth]{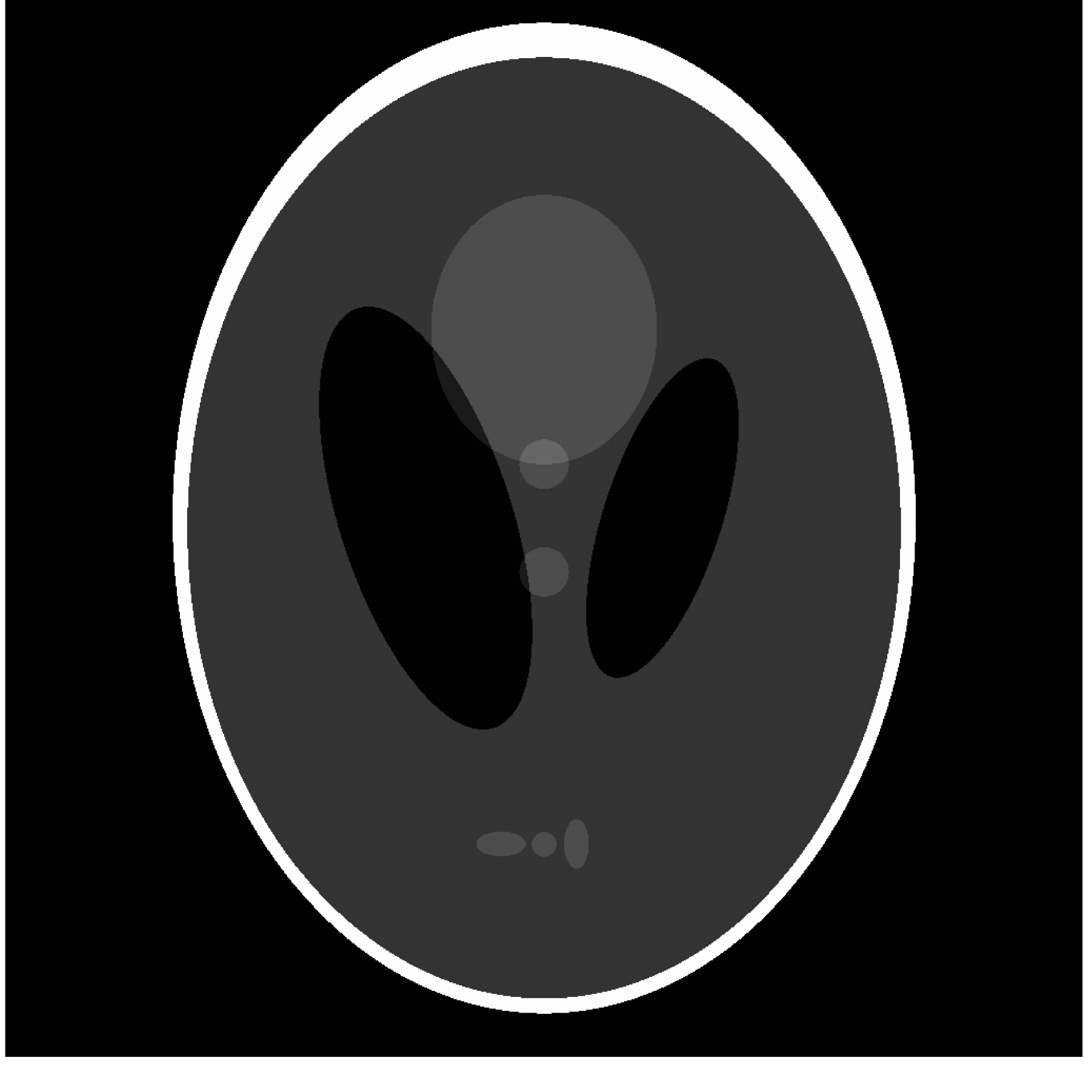}
		\end{subfigure}
		\begin{subfigure}[t]{0.24\textwidth}
			\centering
			\includegraphics[width=\textwidth]{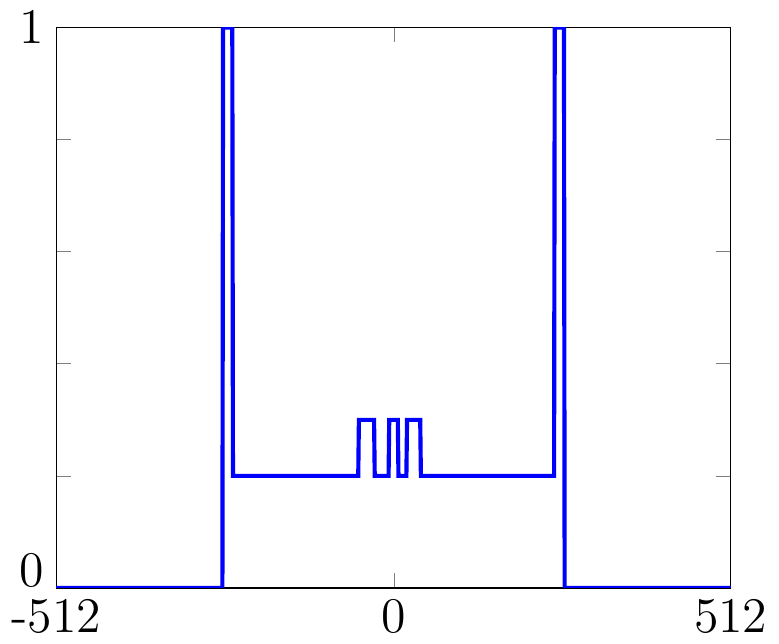}
			\caption{Original phantom}
			\label{fig:reconstr_phantom_linogram_original}
		\end{subfigure}
		\begin{subfigure}[t]{0.24\textwidth}
			\centering
			\includegraphics[width=\textwidth]{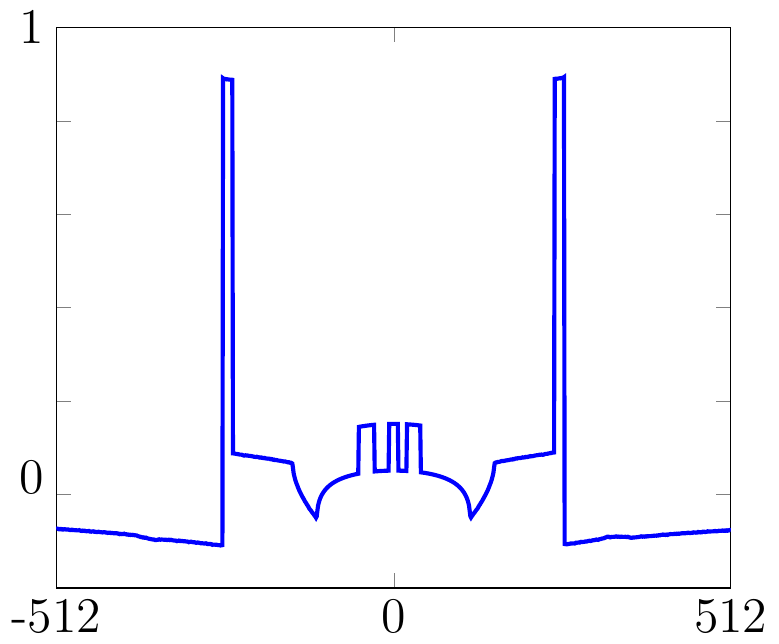}
			\caption{Voronoi weights \\ with \mbox{$e_2$=5.3040e-01}}
			\label{fig:reconstr_phantom_linogram_voronoi}
		\end{subfigure}
		\begin{subfigure}[t]{0.24\textwidth}
			\centering
			\includegraphics[width=\textwidth]{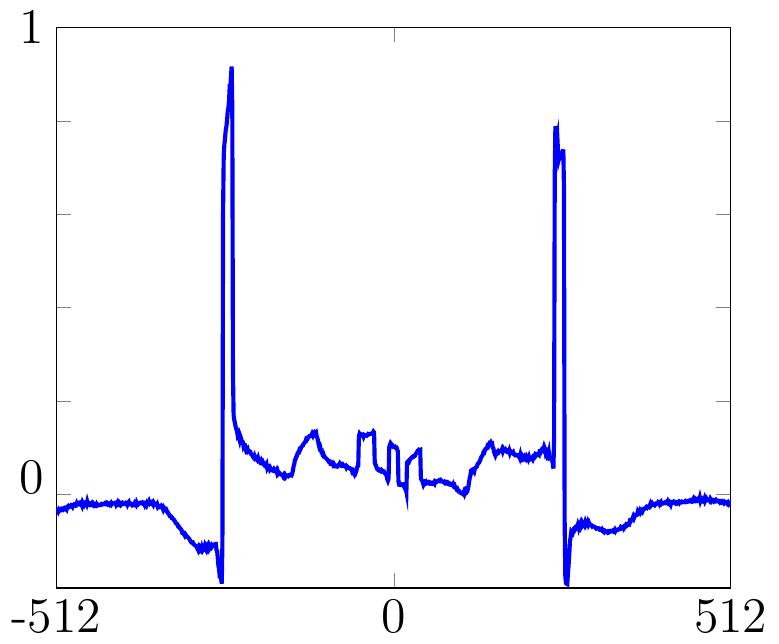}
			\caption{Algorithm~\ref{alg:infft_density} \\ with \mbox{$e_2$=5.0585e-01}}
			\label{fig:reconstr_phantom_linogram_exact}
		\end{subfigure}
		\begin{subfigure}[t]{0.24\textwidth}
			\centering
			\includegraphics[width=\textwidth]{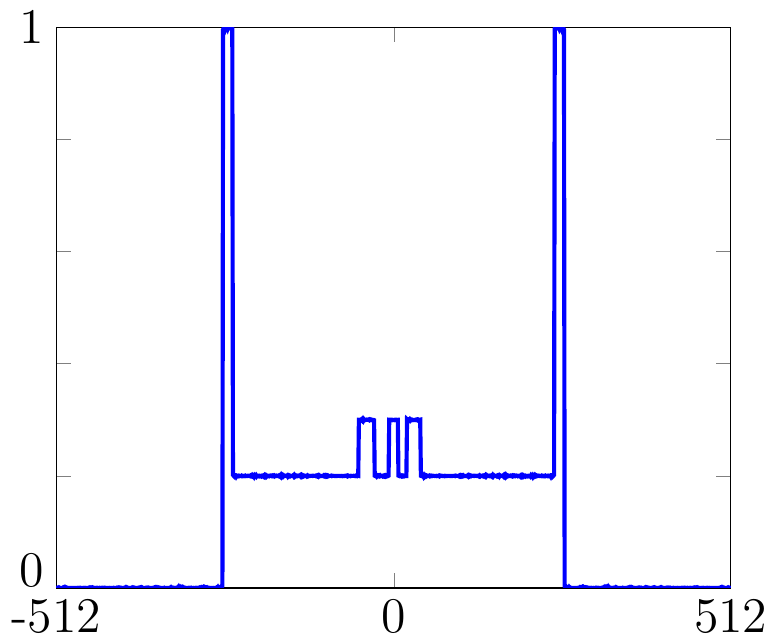}
			\caption{Algorithm~\ref{alg:infft_overdet} \\ with \mbox{$e_2$=2.2737e-03}}
			\label{fig:reconstr_phantom_linogram_opt}
		\end{subfigure}
		\caption{
			Reconstruction of the Shepp-Logan phantom of size \mbox{$1024\times 1024$} (top) via the density compensation method from Section~\ref{sec:trig_poly} using Voronoi weights and Algorithm~\ref{alg:infft_density} with weights computed by \eqref{eq:normal_equations_first_kind_double} compared to Algorithm~\ref{alg:infft_overdet} for the linogram grid~\eqref{eq:linogram_grid} of size \mbox{$R=M=1024$},\, \mbox{$T=2R$}; as well as a detailed view of the 832nd row each (bottom).
			\label{fig:reconstr_phantom_linogram}}
	\end{figure}
	\ex
\end{Example}

\begin{Example}
	\label{ex:reconstr_bandlim}
	Finally, we examine the reconstruction properties for bandlimited functions \mbox{$f\in L_1(\R^d)\cap C_0(\R^d)$} with maximum bandwidth~$\b M$.
	To this end, we firstly specify a compactly supported function $\hat f$ and consequently compute its inverse Fourier transform $f$, such that its samples $f(\b x_j)$ for given \mbox{$\b x_j \in \left[-\frac 12,\frac 12\right)^2$}, \mbox{$j=1,\dots,N$}, can be used for the reconstruction of the samples \mbox{$\hat f(\b k)$}, \mbox{$\b k\in\I_{\b{M}}$}.
	Here we consider the tensorized function \mbox{$\hat f(\b v) = g(v_1)\cdot g(v_2)$},
	where $g(v)$ is the one-dimensional triangular pulse 
	\mbox{$g(v) \coloneqq ( 1-\left|\tfrac{v}{b}\right| ) \cdot \chi_{[-b, b]}(v)$}.
	Then for all \mbox{$b\in\N$} with \mbox{$b\leq \frac M2$} 
	the associated inverse Fourier transform 
	\begin{align*}
		f(\b x)
		=
		\int_{\R^2} \hat f(\b v) \,\mathrm e^{2\pi\mathrm i \b v\b x} \,\mathrm{d}\b v 
		=
		b^2 \,\mathrm{sinc}^2(b\pi \b x)
		=
		b^2 \,\mathrm{sinc}^2(b\pi x_1) \,\mathrm{sinc}^2(b\pi x_2) ,
		\quad \b x\in\R^2 ,
	\end{align*}
	is bandlimited with bandwidth $\b M$.
	In this case, we consider \mbox{$M=64$} and \mbox{$b=24$} as well as the jittered grid~\eqref{eq:jittered_grid} of size \mbox{$N_1=N_2=144$}, i.\,e., we study the setting \mbox{$|\I_{\b {2M}}| \leq N$}.
	
	Now the aim is comparing the different density compensation methods considered in Section~\ref{sec:inv_density_comp} and the optimization approach from Section~\ref{sec:opt_B}.
	More precisely, we consider the reconstruction using Voronoi weights, the weights computed via \eqref{eq:wcf_system_matrix}, the weights in \eqref{eq:weights_pcf} and Algorithm~\ref{alg:infft_density} with weights computed via \eqref{eq:normal_equations_second_kind_double}, as well as Algorithm~\ref{alg:infft_overdet}.
	For the reconstruction \mbox{$\b{\tilde h} \in \{ \b{\tilde h}^{\mathrm{w}}, \b h_{\mathrm{opt}} \}$}, cf.~\eqref{eq:reconstr_density_nfft} and \eqref{eq:approx}, we then compute the pointwise absolute errors \mbox{$\big|\b{\tilde h} - \b{\hat f}\big|$}.
	The corresponding results are displayed in Figure~\ref{fig:reconstr_bandlim_jittered}.
	It can easily be seen that Voronoi weights, see Figure~\ref{fig:reconstr_bandlim_jittered_voronoi}, and the weights in \eqref{eq:weights_pcf}, see Figure~\ref{fig:reconstr_bandlim_jittered_pcf}, do not yield a good reconstruction, as expected.
	The other three approaches produce nearly the same reconstruction error, which is also obtained by reconstruction on an equispaced grid and therefore is the best possible.
	In other words, in case of band\-limited functions the truncation error in \eqref{eq:approx_inverse_integral} is  dominating and thus reconstruction errors smaller than the ones shown in Figure~\ref{fig:reconstr_bandlim_jittered} cannot be expected.
	
	Note that the comparatively small choice of \mbox{$M=64$} was made in order that the computation of the
	weights via \eqref{eq:wcf_system_matrix}, see Figure~\ref{fig:reconstr_bandlim_jittered_wcf},
	as well as
	the weights in \eqref{eq:weights_pcf}, see Figure~\ref{fig:reconstr_bandlim_jittered_pcf},
	is affordable, cf.~Section~\ref{sec:wcf}.
	In contrast, our new methods using Algorithm~\ref{alg:infft_density}, see Figure~\ref{fig:reconstr_bandlim_jittered_exact}, or Algorithm~\ref{alg:infft_overdet}, see Figure~\ref{fig:reconstr_bandlim_jittered_opt}, are much more effective and therefore better suited for the given problem.
	\ex
	\begin{figure}[ht]
		\centering
		\captionsetup[subfigure]{justification=centering}
		\begin{subfigure}[t]{0.32\textwidth}
			\centering 
			\includegraphics[width=\textwidth,trim={2cm 5cm 1cm 0},clip]{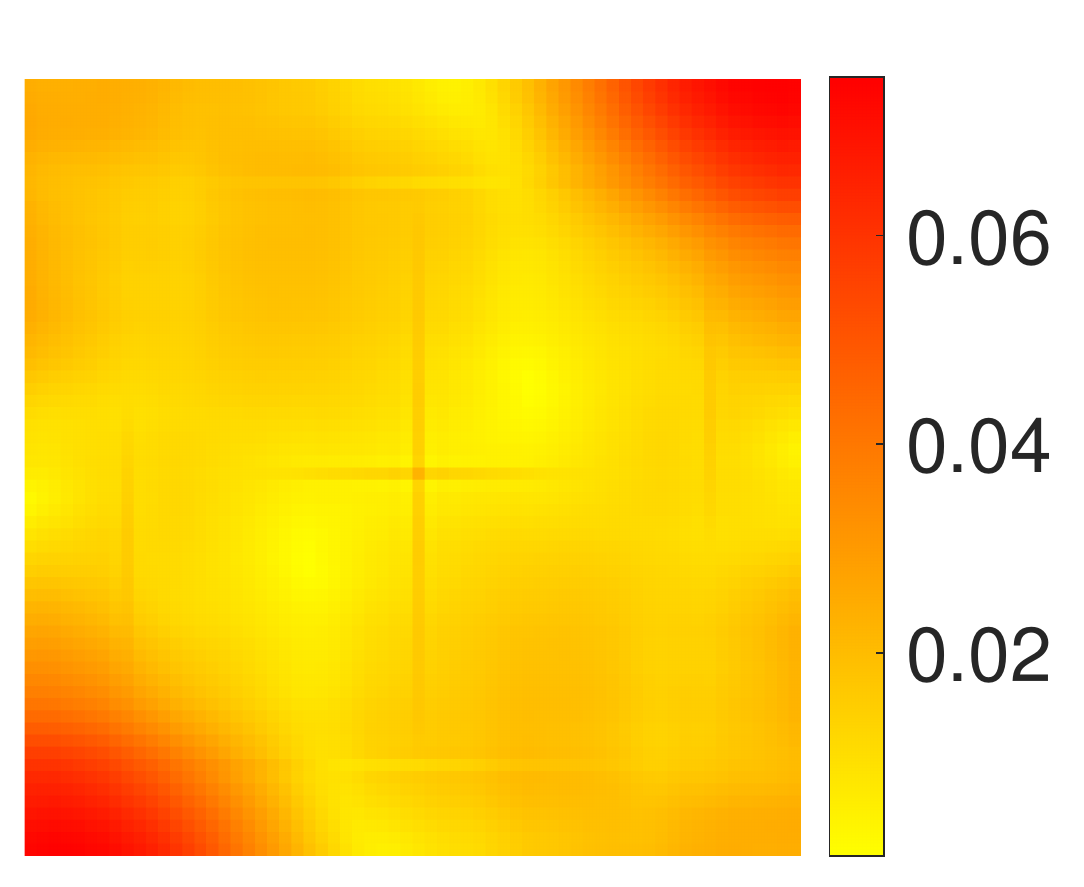}
			\caption{Voronoi weights}
			\label{fig:reconstr_bandlim_jittered_voronoi}
		\end{subfigure}
		\begin{subfigure}[t]{0.33\textwidth}
			\centering
			\includegraphics[width=\textwidth,trim={2cm 5.6cm 1cm 0},clip]{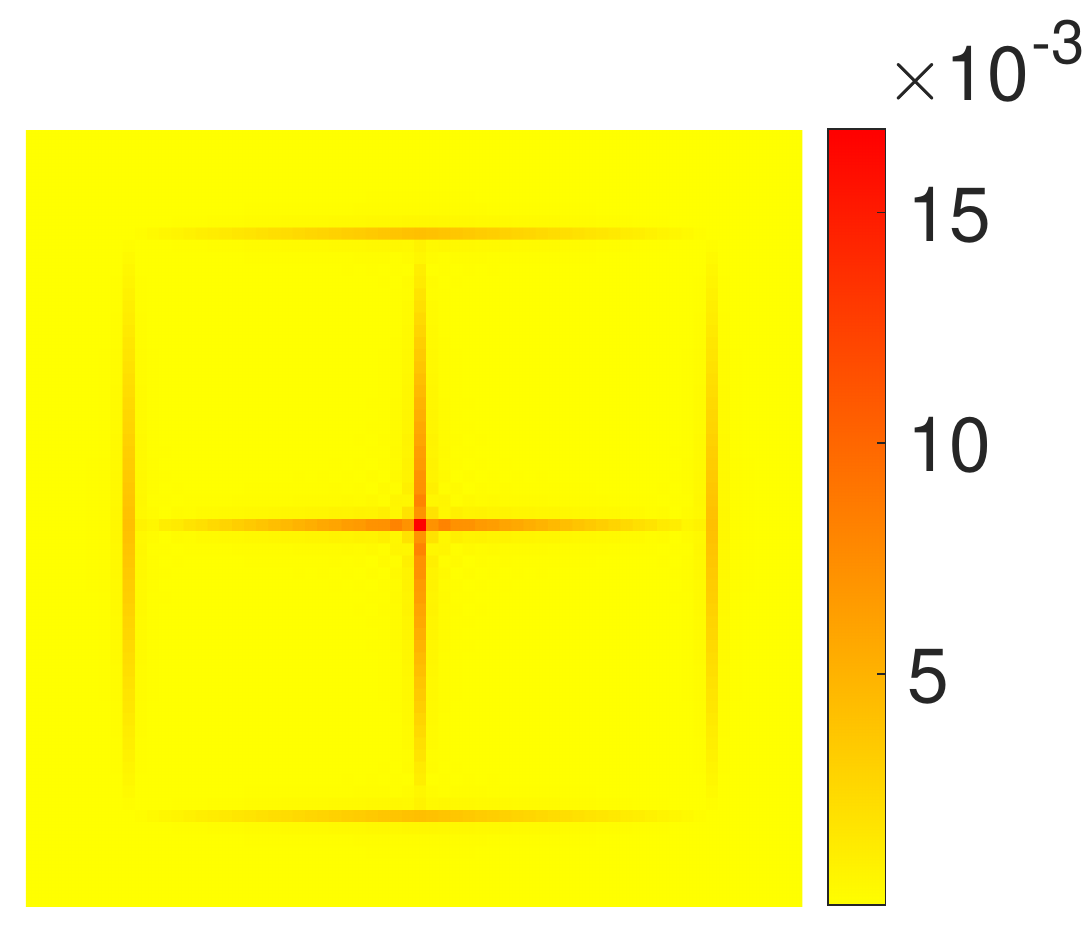}
			\caption{Algorithm~\ref{alg:infft_density}}
			\label{fig:reconstr_bandlim_jittered_exact}
		\end{subfigure}
		\begin{subfigure}[t]{0.33\textwidth}
			\centering
			\includegraphics[width=\textwidth,trim={2cm 5.6cm 1cm 0},clip]{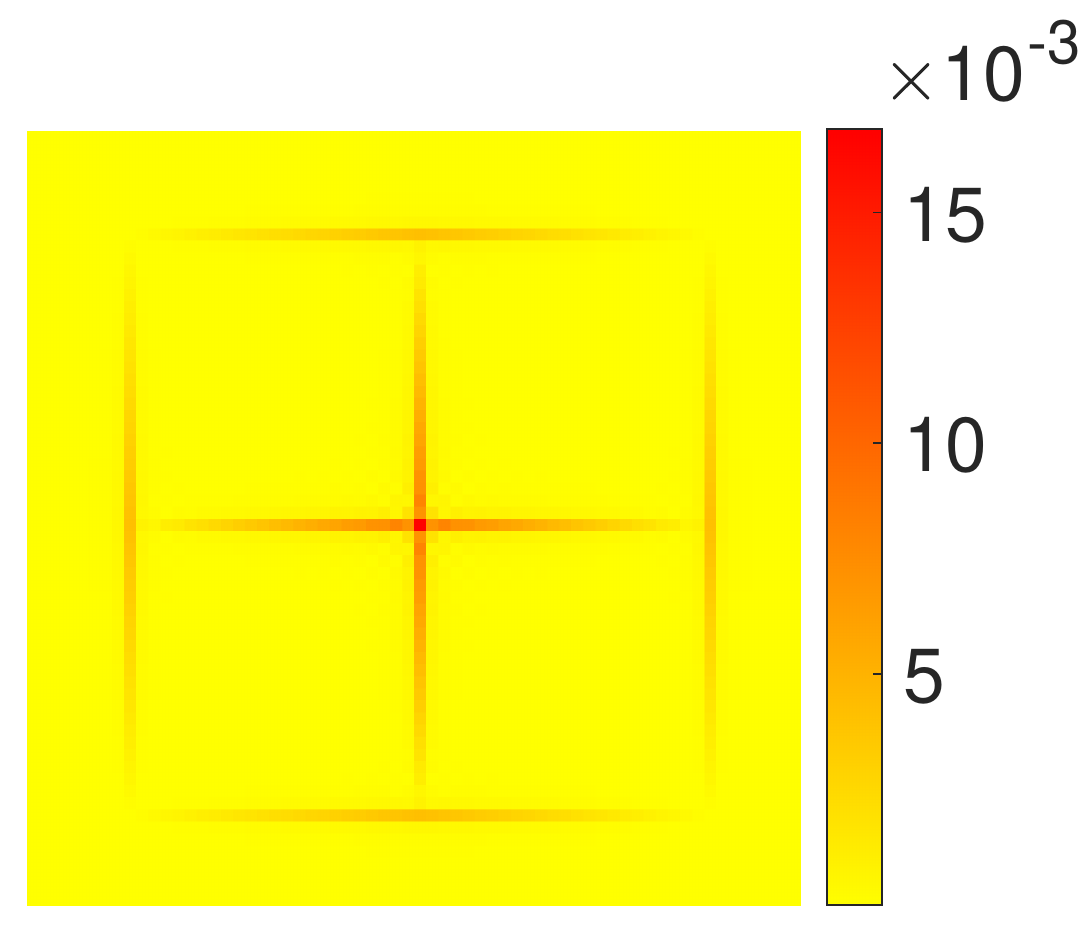}
			\caption{Weights via \eqref{eq:wcf_system_matrix}}
			\label{fig:reconstr_bandlim_jittered_wcf}
		\end{subfigure}
		\begin{subfigure}[t]{0.32\textwidth}
			\centering
			\includegraphics[width=\textwidth,trim={2cm 5cm 1cm 0},clip]{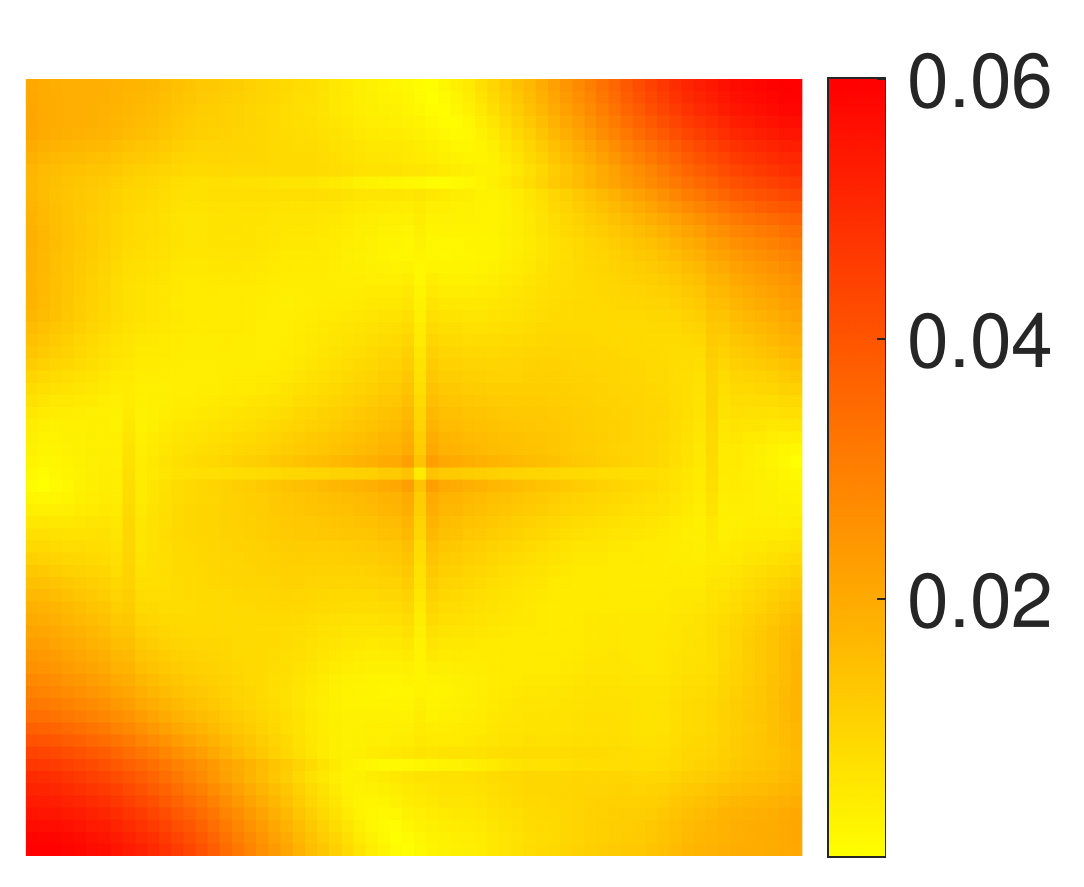}
			\caption{Weights \eqref{eq:weights_pcf}}
			\label{fig:reconstr_bandlim_jittered_pcf}
		\end{subfigure}
		\begin{subfigure}[t]{0.33\textwidth}
			\centering
			\includegraphics[width=\textwidth,trim={2cm 5.6cm 1cm 0},clip]{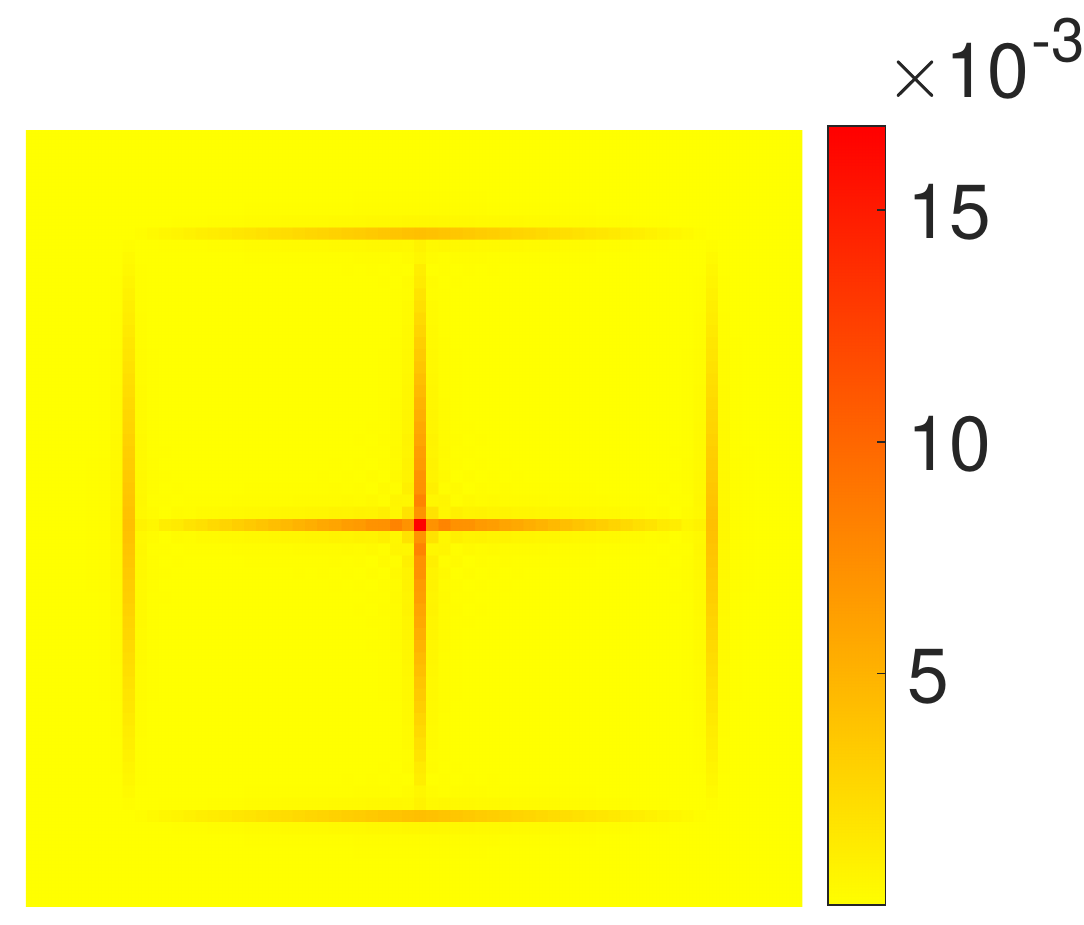}
			\caption{Algorithm~\ref{alg:infft_overdet}}
			\label{fig:reconstr_bandlim_jittered_opt}
		\end{subfigure}
		\caption{
			Pointwise absolute error \mbox{$\big|\b{\tilde h} - \b{\hat f}\big|$} of the reconstruction of samples \mbox{$\hat f(\b k)$} of the tensorized triangular pulse function with \mbox{$M=64$} and \mbox{$b=24$}, using the density compensation methods considered in Section~\ref{sec:inv_density_comp} as well as the optimization approach from Algorithm~\ref{alg:infft_overdet} for the jittered grid~\eqref{eq:jittered_grid} of size \mbox{$N_1=N_2=144$}.
			\label{fig:reconstr_bandlim_jittered}}
	\end{figure}
\end{Example}

\section{Conclusion}

In the present paper we considered several direct methods for computing an inverse NFFT, i.\,e., reconstructing the Fourier coefficients~$\hat f_{\b k}$, \mbox{$\b k\in\I_{\b{M}}$}, from given nonequispaced data~$f(\b x_j)$, \mbox{$j=1,\dots,N$}.
Being a direct method here means, that for a fixed set of points~$\b x_j$, \mbox{$j=1,\dots, N$}, the reconstruction can be realized with the same number of arithmetic operations as a single application of an adjoint NFFT (see Algorithm~\ref{alg:nfft*}).
As we have seen in \eqref{eq:matrix_product_nonequi}, a certain precomputational step is compulsory, since the adjoint NFFT does not yield an inversion by itself. 
Although this precomputations might be rather costly, they need to be done only once for a given set of points 
$\b x_j$, \mbox{$j=1,\dots,N$}.
Therefore, direct methods are especially beneficial in case of fixed points.
For this reason, we studied two different approaches of this kind and especially focused on methods for the multidimensional setting \mbox{$d \geq 1$} that are applicable for general sampling patterns.

Firstly, in Section~\ref{sec:inv_density_comp} we examined the well known approach of sampling density compensation,
where suitable weights \mbox{$w_j\in \C$}, \mbox{$j=1,\dots,N,$} are precomputed, such that the reconstruction can be realized by means of an adjoint NFFT applied to the scaled data \mbox{$w_j f(\b x_j)$}.
We started our investigations with trigonometric polynomials in Section~\ref{sec:trig_poly}.
In Corollary~\ref{Corollary:cond_exact} we introduced the main formula \eqref{eq:claim_quadrature_coefficients_double}, that yields exact reconstruction for all trigonometric polynomials of degree $\b M$.
In addition to this theoretical consider\-ations, we also discussed practical computation schemes for the overdetermined as well as the underdetermined setting, as summarized in Algorithm~\ref{alg:precompute_density}.
Afterwards, in Section~\ref{sec:bandlim} we studied the case of bandlimited functions, which often occurs in the context of MRI, and discussed that the same numerical procedures as in Section~\ref{sec:trig_poly} can be used in this setting as well.
In Section~\ref{sec:error_bound} we then summarized the previous findings by presenting a general error bound on density compensation factors computed by means of Algorithm~\ref{alg:precompute_density} in \autoref{Thm:error_est_lgs}.
In addition, this also yields an estimate on the condition number of the matrix product \mbox{$\b A^* \b W \b A$}, as shown in \autoref{Thm:est_cond}.
In Section~\ref{sec:dcf_literature} we surveyed certain approaches from literature and commented on their connection among each other as well as to the method presented in Section~\ref{sec:trig_poly}.

Subsequently, in Section~\ref{sec:opt_B} we studied another direct inversion method, 
where the matrix representation \mbox{$\b A \approx \b B \b F \b D$} of the NFFT is used to modify the sparse matrix $\b B$, such that a reconstruction is given by \mbox{$\b{\hat f} \approx \b D^* \b F^* \b B_{\mathrm{opt}}^* \b f$}.
In other words, the inversion is done by a modified adjoint NFFT, while the optimization of the matrix $\b B$ can be realized in a precomputational step, see Algorithm~\ref{alg:infft_overdet}.

Finally, in Section~\ref{sec:numerics} we had a look at some numerical examples to investigate the accuracy of the previously introduced methods.
We have seen that our approaches are best-suited for the overdetermined setting \mbox{$|\I_{\b{M}}|\leq N$} and
work for many different sampling patterns.
More specifically, in the highly overdetermined case \mbox{$|\I_{\b{2M}}|\leq N$} we have theoretically proven as well as numerically verified in several examples that the density compensation technique in Algorithm~\ref{alg:infft_density} leads to an exact reconstruction for trigonometric polynomials.
In case not that much data is available and we have to reduce the amount of overdetermination such that \mbox{$|\I_{\b{2M}}| > N$}, we have shown that the optimization approach from Algorithm~\ref{alg:infft_overdet} is preferable, since the
higher number of degrees of freedom in the optimization (see~Remark~\ref{Remark:opt_density}) yields better results.
In addition, also for the setting of bandlimited functions we demonstrated that our methods are much more efficient than existing ones.

\section*{Acknowledgments}
Melanie Kircheis gratefully acknowledges the support from the BMBF grant 01$\mid$S20053A (project SA$\ell$E).
Daniel Potts acknowledges the funding by Deutsche Forschungsgemeinschaft (German Research Foundation) -- Project--ID 416228727 -- SFB 1410.

Moreover, the authors thank the referees and the editor for their very useful suggestions for improvements.

\bibliographystyle{abbrv}

\begin{thebibliography}{10}

\bibitem{Gelb16}
R.~Archibald, A.~Gelb, and R.~B. Platte.
\newblock Image reconstruction from undersampled {F}ourier data using the
  polynomial annihilation transform.
\newblock {\em J. Sci. Comput.}, 67(2):432--452, 2016.

\bibitem{Gelb05}
R.~Archibald, A.~Gelb, and J.~Yoon.
\newblock Polynomial fitting for edge detection in irregularly sampled signals
  and images.
\newblock {\em SINUM}, 43(1):259--279, 2005.

\bibitem{AvCoDoIsSh08}
A.~Averbuch, R.~Coifman, D.~Donoho, M.~Israeli, and Y.~Shkolnisky.
\newblock A framework for discrete integral transformations {I} – the
  pseudopolar {F}ourier transform.
\newblock {\em SIAM J. Sci. Comput.}, 30:764--784, 2008.

\bibitem{AvCoDoElIs}
A.~Averbuch, R.~Coifman, D.~L. Donoho, M.~Elad, and M.~Israeli.
\newblock Fast and accurate polar {F}ourier transform.
\newblock {\em Appl. Comput. Harmon. Anal.}, 21:145--167, 2006.

\bibitem{AvShSh16}
A.~Averbuch, G.~Shabat, and Y.~Shkolnisky.
\newblock Direct inversion of the three-dimensional pseudo-polar {F}ourier
  transform.
\newblock {\em SIAM J. Sci. Comput.}, 38(2):A1100--A1120, 2016.

\bibitem{BaGr03}
R.~F. Bass and K.~Gr{{\"o}}chenig.
\newblock Random sampling of multivariate trigonometric polynomials.
\newblock {\em \textrm{SIAM} J. Math. Anal.}, 36:773--795, 2004.

\bibitem{bey95}
G.~Beylkin.
\newblock On the fast {F}ourier transform of functions with singularities.
\newblock {\em Appl. Comput. Harmon. Anal.}, 2:363--381, 1995.

\bibitem{Bj96}
{\AA}.~Bj{\"o}rck.
\newblock {\em Numerical Methods for Least Squares Problems}.
\newblock \textrm{SIAM}, Philadelphia, PA, USA, 1996.

\bibitem{BoPo06}
A.~B{\"o}ttcher and D.~Potts.
\newblock Probability against condition number and sampling of multivariate
  trigonometric random polynomials.
\newblock {\em Electron. Trans. Numer. Anal.}, 26:178--189, 2007.

\bibitem{CaWaBo08}
E.~J. Candès, M.~B. Wakin, and S.~P. Boyd.
\newblock Enhancing sparsity by reweighted $\ell_1$ minimization.
\newblock {\em J. Fourier Anal. Appl.}, 14(5):877--905, 2008.

\bibitem{ChYi08}
R.~Chartrand and W.~Yin.
\newblock Iteratively reweighted algorithms for compressive sensing.
\newblock In {\em IEEE International Conference on Acoustics, Speech and Signal
  Processing}, pages 3869--3872. IEEE, 2008.

\bibitem{ChMu98}
H.~Choi and D.~C. Munson.
\newblock Analysis and design of minimax-optimal interpolators.
\newblock {\em IEEE Trans. Signal Process.}, 46(6):1571--1579, 1998.

\bibitem{Gelb18}
V.~Churchill, R.~Archibald, and A.~Gelb.
\newblock Edge-adaptive $\ell_2 $ regularization image reconstruction from
  non-uniform {F}ourier data.
\newblock {\em Inverse Probl. Imaging}, 13(5):931--958, 2019.

\bibitem{DaDeFoG10}
I.~Daubechies, R.~DeVore, M.~Fornasier, and C.~S. Güntürk.
\newblock Iteratively reweighted least squares minimization for sparse
  recovery.
\newblock {\em Commun. Pure Appl. Math}, 63(1):1--38, 2010.

\bibitem{MRI22}
M.~Doneva, M.~Akcakaya, and C.~Prieto, editors.
\newblock {\em Magnetic Resonance Image Reconstruction: Theory, Methods and
  Applications}, volume~6.
\newblock Academic Press, 2022.

\bibitem{DuSc}
A.~J.~W. Duijndam and M.~A. Schonewille.
\newblock Nonuniform fast {F}ourier transform.
\newblock {\em Geophysics}, 64:539--551, 1999.

\bibitem{duro93}
A.~Dutt and V.~Rokhlin.
\newblock Fast {F}ourier transforms for nonequispaced data.
\newblock {\em \textrm{SIAM} J. Sci. Stat. Comput.}, 14:1368--1393, 1993.

\bibitem{duro95}
A.~Dutt and V.~Rokhlin.
\newblock Fast {F}ourier transforms for nonequispaced data {II}.
\newblock {\em Appl. Comput. Harmon. Anal.}, 2:85--100, 1995.

\bibitem{EgKiPo22}
H.~Eggers, M.~Kircheis, and D.~Potts.
\newblock {Non-Cartesian MRI reconstruction}.
\newblock In M.~Doneva, M.~Akcakaya, and C.~Prieto, editors, {\em Magnetic
  Resonance Image Reconstruction: Theory, Methods and Applications}, volume~6.
  Academic Press, 2022.

\bibitem{ElSt}
B.~Elbel and G.~Steidl.
\newblock Fast {F}ourier transform for nonequispaced data.
\newblock In C.~K. Chui and L.~L. Schumaker, editors, {\em Approximation Theory
  IX}, pages 39--46, Nashville, 1998. Vanderbilt University Press.

\bibitem{fa07}
G.~E. Fasshauer.
\newblock {\em Meshfree approximation methods with {MATLAB}}.
\newblock World Scientific Publishers, 2007.

\bibitem{FeGrSt95}
H.~G. Feichtinger, K.~Gr{{\"o}}chenig, and T.~Strohmer.
\newblock Efficient numerical methods in non-uniform sampling theory.
\newblock {\em Numer. Math.}, 69:423--440, 1995.

\bibitem{FeKuPo06}
M.~Fenn, S.~Kunis, and D.~Potts.
\newblock On the computation of the polar {FFT}.
\newblock {\em Appl. Comput. Harmon. Anal.}, 22:257--263, 2007.

\bibitem{fesu02}
J.~A. Fessler and B.~P. Sutton.
\newblock Nonuniform fast {F}ourier transforms using min-max interpolation.
\newblock {\em {\textrm{IEEE}} Trans. Signal Process.}, 51:560--574, 2003.

\bibitem{Fou02}
K.~Fourmont.
\newblock Non equispaced fast {F}ourier transforms with applications to
  tomography.
\newblock {\em J. Fourier Anal. Appl.}, 9:431--450, 2003.

\bibitem{GeSo14}
A.~Gelb and G.~Song.
\newblock A frame theoretic approach to the nonuniform fast {F}ourier
  transform.
\newblock {\em SIAM J. Numer. Anal.}, 52(3):1222--1242, 2014.

\bibitem{GoRo22}
A.~Gopal and V.~Rokhlin.
\newblock A fast procedure for the construction of quadrature formulas for
  bandlimited functions.
\newblock {\em Technical Report}, YALEU/DCS/TR-1563, 2022.

\bibitem{GrKuPo09}
M.~Gr{\"a}f, S.~Kunis, and D.~Potts.
\newblock On the computation of nonnegative quadrature weights on the sphere.
\newblock {\em Appl. Comput. Harmon. Anal.}, 27:124--132, 2009.

\bibitem{GrLe04}
L.~Greengard and J.-Y. Lee.
\newblock Accelerating the nonuniform fast {F}ourier transform.
\newblock {\em \textrm{SIAM} Rev.}, 46:443--454, 2004.

\bibitem{GrLeIn06}
L.~Greengard, J.-Y. Lee, and S.~Inati.
\newblock The fast sinc transform and image reconstruction from nonuniform
  samples in {$k$}-space.
\newblock {\em Commun. Appl. Math. Comput. Sci.}, 1:121--131, 2006.

\bibitem{Groechenig20}
K.~Gr{\"o}chenig.
\newblock Sampling, {Marcinkiewicz--Zygmund} inequalities, approximation, and
  quadrature rules.
\newblock {\em J. Approx. Theory}, 257:105455, 2020.

\bibitem{HeRo84}
G.~Heinig and K.~Rost.
\newblock {\em Algebraic methods for {T}oeplitz-like matrices and operators},
  volume~19 of {\em Mathematical Research}.
\newblock Akademie-Verlag, Berlin, 1984.

\bibitem{He19}
E.~S. Helou et~al.
\newblock The discrete {F}ourier transform for golden angle linogram sampling.
\newblock {\em Inverse Problems}, 35(125004), 2019.

\bibitem{Hu09}
D.~Huybrechs.
\newblock Stable high-order quadrature rules with equidistant points.
\newblock {\em J. Comput. Appl. Math.}, 231(2):933--947, 2009.

\bibitem{nfft3}
J.~Keiner, S.~Kunis, and D.~Potts.
\newblock {NFFT 3.5, C subroutine library}.
\newblock \url{http://www.tu-chemnitz.de/~potts/nfft}.
\newblock Contributors: F.~Bartel, M.~Fenn, T.~G\"orner, M.~Kircheis, T.~Knopp,
  M.~Quellmalz, M.~Schmischke, T.~Volkmer, A.~Vollrath.

\bibitem{KeKuPo09}
J.~Keiner, S.~Kunis, and D.~Potts.
\newblock Using {NFFT3} - a software library for various nonequispaced fast
  {Fourier} transforms.
\newblock {\em {ACM} Trans. Math. Software}, 36:Article 19, 1--30, 2009.

\bibitem{KiPo19}
M.~Kircheis and D.~Potts.
\newblock {D}irect inversion of the nonequispaced fast {F}ourier transform.
\newblock {\em Linear Algebra Appl.}, 575:106--140, 2019.

\bibitem{KiPo20}
M.~Kircheis and D.~Potts.
\newblock Efficient multivariate inversion of the nonequispaced fast {F}ourier
  transform.
\newblock {\em PAMM}, 20(1):e202000120, 2021.

\bibitem{KiPoTa21}
M.~Kircheis, D.~Potts, and M.~Tasche.
\newblock {Nonuniform fast Fourier transforms with nonequispaced spatial and
  frequency data and fast sinc transforms}.
\newblock {\em Numer. Algor.}, 2022.

\bibitem{KiPoTa22}
M.~Kircheis, D.~Potts, and M.~Tasche.
\newblock {On regularized Shannon sampling formulas with localized sampling}.
\newblock {\em Sampl. Theory Signal Process. Data Anal.}, 20(20), 2022.

\bibitem{KnKuPo}
T.~Knopp, S.~Kunis, and D.~Potts.
\newblock A note on the iterative {MRI} reconstruction from nonuniform k-space
  data.
\newblock {\em Int. J. Biomed. Imag.}, 2007, 2007.
\newblock ID 24727.

\bibitem{Kotelnikov}
V.~A. Kotelnikov.
\newblock On the transmission capacity of the “ether” and wire in
  electrocommunications.
\newblock In {\em Modern Sampling Theory: Mathematics and Application}, pages
  27--45. Birkh\"auser, Boston, 2001.
\newblock Translated from Russian.

\bibitem{KuNa18}
S.~Kunis and D.~Nagel.
\newblock On the condition number of {V}andermonde matrices with pairs of
  nearly-colliding nodes.
\newblock {\em Numer. Algor.}, 87:473–496, 2021.

\bibitem{kupo04}
S.~Kunis and D.~Potts.
\newblock Stability results for scattered data interpolation by trigonometric
  polynomials.
\newblock {\em \textrm{SIAM} J. Sci. Comput.}, 29:1403--1419, 2007.

\bibitem{KuPo06}
S.~Kunis and D.~Potts.
\newblock Time and memory requirements of the nonequispaced {FFT}.
\newblock {\em Sampl. Theory Signal Image Process.}, 7:77--100, 2008.

\bibitem{LeGr05}
J.-Y. Lee and L.~Greengard.
\newblock The type 3 nonuniform {FFT} and its applications.
\newblock {\em J. Comput. Physics}, 206:1--5, 2005.

\bibitem{LZ16}
R.~Lin and H.~Zhang.
\newblock {Convergence analysis of the Gaussian regularized Shannon sampling
  formula}.
\newblock {\em Numer. Funct. Anal. Optim.}, 38(2):224--247, 2017.

\bibitem{LiMa12}
Y.~Liu, J.~Ma, Y.~Fan, and Z.~Liang.
\newblock Adaptive-weighted total variation minimization for sparse data toward
  low-dose x-ray computed tomography image reconstruction.
\newblock {\em Phys. Med. Biol.}, 57(23):7923, 2012.

\bibitem{LuBo}
J.~Lund and K.~L. Bowers.
\newblock {\em {Sinc Methods for Quadrature and Differential Equations}}.
\newblock Society for Industrial and Applied Mathematics, 1992.

\bibitem{MXZ09}
C.~Micchelli, Y.~Xu, and H.~Zhang.
\newblock Optimal learning of bandlimited functions from localized sampling.
\newblock {\em J. Complexity}, 25(2):85--114, 2009.

\bibitem{st01}
A.~Nieslony and G.~Steidl.
\newblock Approximate factorizations of {F}ourier matrices with nonequispaced
  knots.
\newblock {\em Linear Algebra Appl.}, 266:337--351, 2003.

\bibitem{PiMe99}
J.~G. Pipe and P.~Menon.
\newblock Sampling density compensation in {MRI}: rationale and an iterative
  numerical solution.
\newblock {\em Magn. Reson. Med.}, 41:179--186, 1999.

\bibitem{PlPoStTa18}
G.~Plonka, D.~Potts, G.~Steidl, and M.~Tasche.
\newblock {\em Numerical Fourier Analysis}.
\newblock Applied and Numerical Harmonic Analysis. Birkh\"auser, 2018.

\bibitem{PoTa21b}
D.~Potts and M.~Tasche.
\newblock {Continuous window functions for NFFT}.
\newblock {\em Adv. Comput. Math.}, 47(53):1--34, 2021.

\bibitem{PoTa21a}
D.~Potts and M.~Tasche.
\newblock {Uniform error estimates for nonequispaced fast Fourier transforms}.
\newblock {\em Sampl. Theory Signal Process. Data Anal.}, 19(17):1--42, 2021.

\bibitem{PrWa01}
K.~P. Pruessmann and F.~T. A.~W. Wayer.
\newblock Major speedup of reconstruction for sensitivity encoding with
  arbitrary trajectories.
\newblock In {\em Proc. Intl. Soc. Mag. Reson. Med. 9}, page 767. Glasgow,
  Scotland, 2001.

\bibitem{Q03}
L.~Qian.
\newblock {On the regularized Whittaker--Kotelnikov--Shannon sampling formula}.
\newblock {\em Proc. Amer. Math. Soc.}, 131(4):1169--1176, 2003.

\bibitem{RaPrSiBoEg99}
V.~Rasche, R.~Proksa, R.~Sinkus, P.~B{\"o}rnert, and H.~Eggers.
\newblock Resampling of data between arbitrary grids using convolution
  interpolation.
\newblock {\em {\textrm{IEEE}} Trans. Med. Imag.}, 18:385--392, 1999.

\bibitem{Rosenfeld98}
D.~Rosenfeld.
\newblock An optimal and efficient new gridding algorithm using singular value
  decomposition.
\newblock {\em Magn. Reson. Med.}, 40(1):14--23, 1998.

\bibitem{RuTo18}
D.~Ruiz-Antolin and A.~Townsend.
\newblock A nonuniform fast {F}ourier transform based on low rank
  approximation.
\newblock {\em SIAM J. Sci. Comput.}, 40(1):A529--A547, 2018.

\bibitem{Gelb19}
T.~Scarnati and A.~Gelb.
\newblock Accelerated variance based joint sparsity recovery of images from
  {F}ourier data.
\newblock arXiv, 2019.

\bibitem{Sedarat00}
H.~Sedarat and D.~G. Nishimura.
\newblock On the optimality of the gridding reconstruction algorithm.
\newblock {\em IEEE Trans. Med. Imaging}, 19(4):306--317, 2000.

\bibitem{Selva18}
J.~Selva.
\newblock Efficient type-4 and type-5 non-uniform {FFT} methods in the
  one-dimensional case.
\newblock {\em IET Signal Processing}, 12(1):74--81, 2018.

\bibitem{Shannon49}
C.~E. Shannon.
\newblock Communication in the presence of noise.
\newblock {\em Proc. I.R.E.}, 37:10--21, 1949.

\bibitem{Gelb10}
W.~Stefan, R.~A. Renaut, and A.~Gelb.
\newblock Improved total variation-type regularization using higher order edge
  detectors.
\newblock {\em SIIMS}, 3(2):232--251, 2010.

\bibitem{st97}
G.~Steidl.
\newblock A note on fast {F}ourier transforms for nonequispaced grids.
\newblock {\em Adv. Comput. Math.}, 9:337--353, 1998.

\bibitem{St98}
G.~Stewart.
\newblock {\em Matrix Algorithms: Basic decompositions}.
\newblock SIAM, Philadelphia, 1998.

\bibitem{StTa06}
T.~Strohmer and J.~Tanner.
\newblock Fast reconstruction methods for bandlimited functions from periodic
  nonuniform sampling.
\newblock {\em SIAM J. Numer. Anal.}, 44(3):1071--1094, 2006.

\bibitem{SuFeNo01}
B.~P. Sutton, J.~A. Fessler, and D.~C. Noll.
\newblock A min-max approach to the nonuniform {$N$}-dimensional {FFT} for
  rapid iterative reconstruction of {MR} images.
\newblock In {\em Proc. ISMRM 9th Scientific Meeting}, page 763, 2001.

\bibitem{Tch57}
V.~Tchakaloff.
\newblock Formules de cubature mécaniques à coefficients non négatifs.
\newblock {\em Bull. Sci. Math.}, 81:123--134, 1957.

\bibitem{Gelb15}
G.~Wasserman, R.~Archibald, and A.~Gelb.
\newblock Image reconstruction from {F}ourier data using sparsity of edges.
\newblock {\em J. Sci. Comput.}, 65(2):533--552, 2015.

\bibitem{Whittaker}
E.~T. Whittaker.
\newblock On the functions which are represented by the expansions of the
  interpolation theory.
\newblock {\em Proc. R. Soc. Edinb.}, 35:181--194, 1915.

\end{thebibliography}

\end{document}